\documentclass[a4paper,12pt,twoside]{article}
\usepackage[latin1]{inputenc}
\usepackage{amsfonts,amssymb,amsmath,exscale}
\usepackage[dvips]{graphicx,color,psfrag}
\usepackage{fancyhdr}
\usepackage{rotating}
\usepackage{pifont}
\usepackage[hyperref]{xcolor}
\definecolor{myblue}{RGB}{0,0,139}
\usepackage[colorlinks=true,citecolor=blue,urlcolor=blue]{hyperref}
\hypersetup{bookmarksnumbered}
\usepackage[square]{natbib}
\usepackage[font=small,labelfont=bf]{caption}
\usepackage{setspace}
\fancypagestyle{plain}{%
  \fancyhead{}%
  \fancyfoot[c]{\sffamily\thepage}%
}
\makeatletter
\def\cleardoublepage{\clearpage\if@twoside \ifodd\c@page\else
  \hbox{}
  \vspace*{\fill}
  \thispagestyle{empty}
  \newpage
  \if@twocolumn\hbox{}\newpage\fi\fi\fi}
\usepackage{algorithm2e}
\newcommand{\increment}{\text{\reflectbox{$\Delta$}}}
\usepackage{amsmath,amssymb}
\renewcommand{\citep}[1]{(\cite{#1}{})}

\newcommand{\ce}{\centerline}
\newcommand{\B}{{\cal B}}

\DeclareRobustCommand{\BB}{{\boldsymbol{\mathnormal B}}}

\DeclareRobustCommand{\BK}{{\boldsymbol{\mathnormal K}}}

\DeclareRobustCommand{\Bn}{{\boldsymbol{\mathnormal n}}}

\DeclareRobustCommand{\Bs}{{\boldsymbol{\mathnormal s}}}
\DeclareRobustCommand{\Bt}{{\boldsymbol{\mathnormal t}}}
\DeclareRobustCommand{\Bu}{{\boldsymbol{\mathnormal u}}}

\DeclareRobustCommand{\Bx}{{\boldsymbol{\mathnormal x}}}

\DeclareRobustCommand{\calB}{{\mathcal B}}
\DeclareRobustCommand{\calH}{{\mathcal H}}
\DeclareRobustCommand{\calR}{{\mathcal R}}
\DeclareRobustCommand{\calT}{{\mathcal T}}
\DeclareRobustCommand{\calW}{{\mathcal W}}
\DeclareRobustCommand{\Bve}{{\boldsymbol{\varepsilon}}}
\DeclareMathAlphabet{\Inbb}{U}{bbmss}{m}{n}
\DeclareRobustCommand{\nIE}{{\Inbb E}}
\DeclareRobustCommand{\nIh}{{\Inbb h}}
\DeclareRobustCommand{\Bsigma}{{\boldsymbol{\sigma}}}
\DeclareMathAlphabet{\gothic}{U}{euf}{m}{n}
\DeclareRobustCommand{\frakD}{{\gothic D}}
\DeclareRobustCommand{\frakk}{{\gothic k}}
\DeclareRobustCommand{\frakW}{{\gothic W}}
\DeclareMathAlphabet{\Bgothic}{U}{euf}{b}{n}
\DeclareRobustCommand{\BfrakC}{{\Bgothic C}}
\DeclareRobustCommand{\BfrakF}{{\Bgothic F}}
\DeclareRobustCommand{\BfrakU}{{\Bgothic U}}
\newcommand{\tr}{\mathop{\operator@font tr}}
\newcommand{\dev}{\mathop{\operator@font dev}}
\newcommand{\half}{{\frac{1}{2}}}
\DeclareMathAlphabet{\Ibb}{U}{msb}{m}{n}
\DeclareRobustCommand{\IE}{{\Ibb E}}
\DeclareRobustCommand{\IF}{{\Ibb F}}
\DeclareRobustCommand{\IX}{{\Ibb X}}
\newcommand{\bn}{{\mbox{\boldmath$n$}}}
\newcommand{\br}{{\mbox{\boldmath$r$}}}
\newcommand{\ebn}{\[}
\newcommand{\een}{\]\@ignoretrue}
\DeclareRobustCommand{\Bone}{{\boldsymbol{\mathit 1}}}
\DeclareRobustCommand{\Bzero}{{\boldsymbol{\mathit 0}}}
\DeclareRobustCommand{\Bnzero}{{\boldsymbol 0}}
\newcommand{\bzero}{\Bnzero}
\DeclareRobustCommand{\Bvarepsilon}{{\boldsymbol{\varepsilon}}}
\renewcommand{\div}{\mathop{\operator@font div}}
\newcommand{\WITH}{\quad\mbox{with}\quad}
\newcommand{\AND}{\quad\mbox{and}\quad}

\begin{document}

\thispagestyle{empty}

\ce{\textbf{\large A variational minimization formulation for hydraulically}}

\medskip

\ce{\textbf{\large induced fracturing  in elastic-plastic solids}}

\bigskip

\ce{Daniel Kienle\footnote{%
\footnotesize Corresponding author.
\textit{E-mail address}: \href{mailto:kienle@mechbau.uni-stuttgart.de}{kienle@mechbau.uni-stuttgart.de}}, Marc-Andr\'e Keip
}

\bigskip

\ce{Institute of Applied Mechanics}
\ce{Department of Civil and Environmental Engineering}
\ce{University of Stuttgart, Stuttgart, Germany}

\bigskip

\noindent \textbf{Abstract.} A variational modeling framework for hydraulically induced fracturing of elastic-plastic solids is developed in the present work. The developed variational structure provides a global minimization problem. While fracture propagation is modeled by means of a phase-field approach to fracture, plastic effects are taken into account by using a Drucker--Prager-type yield-criterion function. This yield-criterion function governs the plastic evolution of the fluid-solid mixture. Fluid storage and transport are described by a Darcy--Biot-type formulation. Thereby the fluid storage is decomposed into a contribution due to the elastic deformations and one due to the plastic deformations. A local return mapping scheme is used for the update of the plastic quantities. The global minimization structure demands a $H(\div)$-conforming finite-element formulation. Furthermore this is combined with an enhanced-assumed-strain formulation in order to overcome locking phenomena arising from the plastic deformations. The robustness and capabilities of the presented framework will be shown in a sequence of numerical examples.


\bigskip

\section{Introduction}
\label{Section1}

One of the main applications for a model of hydraulically induced fractures can be found in the recently utilized oil-production technique called fracking. During this process a highly pressurized fluid is injected into a perforated bore hole. The goal is to induce fractures in layers of the earth's crust that store large amounts of oil and gas. The fractures are created to increase the fluid permeability of the present soil or rock, therewith leading to a higher flow of gas and oil into the bore hole where it is collected. The main criticisms of this technique are related to the inducement of fractures that could lead to environmental issues such as seismic activities or contamination of drinking water. Therefore, in recent years various attempts have been made to better capture and understand the underlying physical processes, for example by the development of models that can be used in numerical simulations. The idea is that an accurate and robust model can help in forecasting both risks and potential of such a technique. To be successful, such a model has to capture three distinct mechanisms: First, it has to describe the mechanical deformation of the present porous medium (which could be soil or rock); second, it has to describe the fluid transport throughout the intact and fractured porous medium; third, the model needs to be capable of describing fracture initiation and propagation, for example driven by fluid pressure.

To embed the present work into the literature, we comment on recent development in the field. For the modeling of hydraulic fracturing, it is of substantial importance to describe the mechanical deformation of the underlying porous media. Powerful techniques in this direction are given on the one hand by model formulations within the context of the so-called Theory of Porous Media, see
\cite{deboer00},
\cite{ehlers02}
and
\cite{bluhm1997volume}.
An alternative approach that is suitable for the description of fully saturated porous media with one fluid phase is formulated in the seminal works of
\cite{terzaghi25}
and
\cite{biot41}.
Associated models reside in the area of the Biot theory of consolidation, see
\cite{detournay+cheng93}
and
\cite{bear72}
for a general overview.

When it comes to the modeling of fluid flow within a porous medium, it is often made use of the law of
\cite{darcy56}.
Darcy's law provides a phenomenological approach to the modeling of fluid flux that is driven by gradients of fluid pressure or, more precisely, by gradients of the chemical potential. Darcy's law is often referred to as a macroscopic approach since it describes the fluid flow through the porous medium's pore scale in an averaged or homogenized sense. Fluid flow within fractures, i.e.\ regions without solid content, could directly be modeled based on the Navier--Stokes-equation. Associated simplifications can be based on the lubrication theory and lead to so-called Poiseuille-type flow models. While the first approach is often used in combination with models related to the Theory of Porous Media, the latter one can nicely be embedded in formulations that relate to Biot's theory of consolidation.

Since the present work aims at the description of fracturing processes in porous media, we briefly refer to some related literature. Fundamental concepts for the description of fracture were proposed by
\cite{griffith21}
and
\cite{irwin58}.
These works provide energy-based fracture criteria for brittle materials and build the theoretical basis for the recently developed phase-field models of fracturing. With regard to the latter we highlight the contributions of
\cite{francfort+marigo98},
\cite{bourdin+francfort+marigo00,bourdin+francfort+marigo08b},
\cite{amor+marigo+maurini09},
\cite{kuhn+mueller10},
\cite{miehe+hofacker+welschinger10}
and
\cite{pham+marigo+maurini11a}.
The phase-field approach to fracture is extremely versatile and has been extended to the modeling of ductile fracture
(\cite{miehe+hofacker+schaenzel+aldakheel15,miehe+aldakheel+raina16,miehe+kienle+aldakheel+teichtmeister16,miehe+teichtmeister+aldakheel16},
\cite{ambati+gerasimov+lorenzis15},
\cite{borden+huges+landis+anvari+lee16},
\cite{alessi+margio+maurini+vidoli18},
\cite{steinke2020modelling},
\cite{yin2020ductile}
),
anisotropic fracturing
(\cite{li2015phase},
\cite{teichtmeisteretal17},
\cite{storm2020concept}
)
and fatigue fracture
(\cite{alessi+vidoli+delorenzis17},
\cite{lo2019phase},
\cite{schreiber2020phase},
\cite{carrara2020framework}),
just to name a few. 

Next to the above mentioned extensions, the phase-field approach to fracture has also seen pronounced activity in the field of porous media. Here we highlight the contributions of
\cite{mikelic+wheeler+wick15,mikelic+wheeler+wick15b,mikelic2015phase},
\cite{miehe+mauthe+teichtmeister15},
\cite{wilson+landis16},
\cite{wu+delorenzis16},
\cite{mauthe+miehe17}
and
\cite{cajuhi+sanavia+delorenzis18}
that are embedded into Biot's theory of consolidation. In that connection we also mention a recent work on the numerical treatment of poro-elastic problems, in which an emphasis is put on suitable finite-element formulations
\cite{teichtmeister+mauthe+miehe19}.
Combinations of the Theory of Porous Media with the phase-field approach to fracture have been provided by
\cite{ehlers17,ehlers+luo18a,ehlers+luo18b,heider17}.

We note that the above mentioned models all relate to elastic deformations of the material. However, there is experimental evidence that elastic material response alone cannot capture all relevant effects
\cite{johnson+clearly91}.
This serves as a motivation to the develop models that incorporate elastic-plastic deformations. Here we would like to refer to the recent developments of
\cite{pise+bluhm+schroeder19}
and
\cite{aldakheel20}.
In the given works, the plastic response is incorporated by introducing a Drucker--Prager-type yield criterion function
\cite{drucker+prager52}.
We note that the works
\cite{pise+bluhm+schroeder19},
\cite{aldakheel20}
are not variational and consider a yield-criterion function in terms of the effective stresses. The present work provides an alternative formulation of elastic-plastic hydraulic fracturing in a rigorously variational setting motivated by the ideas of
\cite{armero99}.
In the latter work, an additive split of the fluid content into an elastic and a plastic part is suggested, which leads to a constitutive fluid pressure as a function of only elastic quantities. This concept is perfectly suitable for a variational framework that combines the elastic-plastic deformation of the porous medium with the fluid transport on the one hand and the fracture initiation and evolution on the other.

The present work is structured as follows. In Section~\ref{Section2} the unknown fields and associated kinematic relations are introduced. Thereafter, in \ref{Section3}, we provide the variational framework and the constitutive functions. Here, we make use of a Darcy--Biot-type formulation for the fluid transport in the porous medium and a Drucker--Prager-type yield-criterion function for plastic flow. This then leads to the set of Euler equations that describe the behavior of a fracturing porous-elastic-plastic solid. The numerical treatment of the problem is based on a time-discrete incremental variational formulation combined with a local return-mapping scheme and discussed in Section~\ref{Section4}. We showcase the capabilities of the presented numerical framework for hydraulically induced fracturing of porous-elastic-plastic solids by means of some numerical examples to be discussed in Section~\ref{Section5}. In detail we consider a rigid-footing test as well as two fluid-injection tests. A summary and an outlook will be provided in Section~\ref{Section6}.

\section{Independent and primary fields}
\label{Section2}

In the present section we introduce the independent and primary fields of the porous-elastic-plastic fracture model. They account for the elastic-plastic deformation of a body $\B$, for the fluid flux and storage in $\B$ as well as for the fracture initiation and evolution in $\B$. The surface of the body will in the following be denotes as $\partial \B$.

\paragraph{Small-strain kinematics.}

The model presented in this work is formulated in the context of the infinitesimal-strain theory. Hence, we consider the displacement field $\Bu(\Bx,t)$ at the material point $\Bx \in \calB$ and time $t$
\begin{equation}
\Bu  : \left\{
\begin{array}{ll}
\B \times \calT \rightarrow \calR^3 \\
(\Bx,t) \mapsto \Bu(\Bx,t)
\end{array}
\right.
\label{eq:displacement}
\end{equation}
as independent field. Based on that, we introduce the infinitesimal strain tensor as a primary variable. It cam be computed from the displacement gradient by
\begin{equation}
\Bve:=\tfrac{1}{2}(\nabla \Bu + \nabla \Bu^T) .
\label{eq:strain}
\end{equation}

\paragraph{Plastic deformations.}

Since we are interested in modeling ductile response of the material, we additively decompose the strain tensor \eqref{eq:strain} into elastic and plastic contributions
\begin{equation}
\Bve = \Bve^e + \Bve^p ,
\label{eq:strain}
\end{equation}
wherein the plastic strain $\Bve^p$ will be treated as a local internal variable. Further, to describe local isotropic hardening, we introduce a local hardening field $\alpha$ formally by
\begin{equation}
\alpha  : \left\{
\begin{array}{ll}
\B \times \calT \rightarrow \calR \\
(\Bx,t) \mapsto \alpha(\Bx,t).
\end{array}
\right.
\label{eq:hardening}
\end{equation}

\paragraph{Fluid mass and fluid flux.}

The initial density (mass per unit volume) of the fluid-solid mixture is denoted by $m_0$. It can be computed from the densities of the fluid $\rho_f$ and the solid $\rho_s$ through the given porosity of the mixture $\varphi$ as
\begin{equation}
m_0  = \rho_f \varphi + \rho_s ( 1 - \varphi)
\WITH
\varphi = \tfrac{V_\text{pore}}{V_\text{pore} + V_\text{solid}}.
\end{equation}  
In the above equation, the porosity $\varphi$ has been obtained from the given pore and solid volume $V_\text{pore}$ and $V_\text{solid}$, respectively. The sum of solid and fluid volume is often denoted as bulk volume $V_\text{bulk}$.

Following Biot's approach to thermodynamically open systems, the mass balance in its global and local representation reads
\citep{biot84}
\begin{equation}
\displaystyle \frac{\text{d}}{\text{d}t}~{\int_\B} ~m_0 + m ~\text{d}V = -{\int_{\partial \B}}~ \nIh \cdot \Bn ~\text{d}A
\qquad \Leftrightarrow \qquad
  \dot m = - \div \nIh \quad \text{in}~\B.
\label{eq:globalfluidmass}
\end{equation}  
Here, $m$ denotes the change in the bulk's density, which is caused by fluid flux through the body's surface. Thus, the vector $\nIh$ given on the right-hand side of \eqref{eq:globalfluidmass} denotes a fluid-flux vector. Formally, the latter two quantities can be introduced as
\begin{equation}
m  : \left\{
\begin{array}{ll}
\B \times \calT \rightarrow \calR \\
(\Bx,t) \mapsto m(\Bx,t)
\end{array}
\right.
\AND
\nIh  : \left\{
\begin{array}{ll}
\B \times \calT \rightarrow \calR^3 \\
(\Bx,t) \mapsto \nIh(\Bx,t).
\end{array}
\right.
\label{eq:fluidunknown}
\end{equation}
In what follows we denote $m$ as the change of fluid content. Analogous to the strain tensor, we decompose $m$ into an elastic and a plastic part
\begin{equation}
m = m^e + m^p.
\label{eq:masscontent}
\end{equation}
Here, the plastic contribution $m^p$ describes a change in bulk density that is caused by plastic deformations. It is associated with fluid that is remanently squeezed out or soaked in due to plastic deformations.
 
\begin{figure}
\centering
{
\scriptsize
\psfrag{x}      [c][c]{$\Bx \in \B$}%
\psfrag{u}      [c][c]{$\Bu$}%
\psfrag{bu}     [c][c]{$\Bu=\bar \Bu$}%
\psfrag{n}      [c][c]{$\Bn$}%
\psfrag{t}      [c][c]{$\Bsigma \Bn =\bar \Bt$}%
\psfrag{field1} [c][c]{displacement field}%
\psfrag{m}      [c][c]{$m^p$}%
\psfrag{field2} [c][c]{plastic fluid content}%

\psfrag{h}      [c][c]{$\nIh$}%
\psfrag{bh}     [c][c]{$\nIh\cdot\Bn=\bar h$}%
\psfrag{mu}     [c][c]{$\mu=\bar \mu$}%
\psfrag{field3} [c][c]{fluid flux}%
\psfrag{d}      [c][c]{$d$}%
\psfrag{bd}     [c][c]{$d=\bar d$}%
\psfrag{gd}     [c][c]{$c {\nabla d}\cdot \Bn =0$}%
\psfrag{field4} [c][c]{fracture phase-field}%
\psfrag{e}      [c][c]{$\Bve^p$}%
\psfrag{field5} [c][c]{plastic strain}%
\psfrag{a}      [c][c]{$\alpha$}%
\psfrag{field6} [c][c]{hardening}%
\psfrag{mp}     [c][c]{$m$}%
\psfrag{field7} [c][c]{fluid content}%
\includegraphics[width=0.85\textwidth]{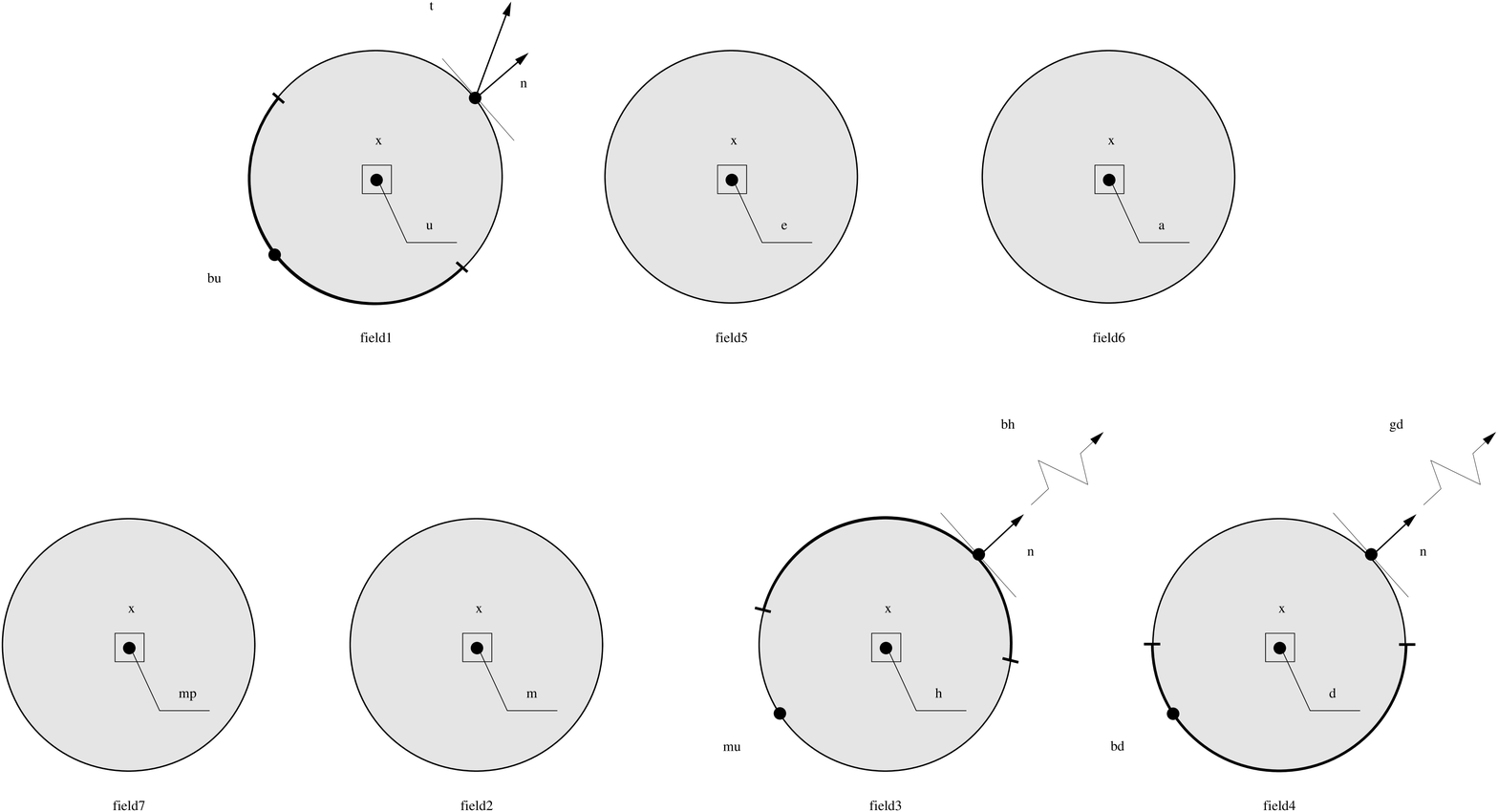}
}
\vspace{1mm}%
\captionsetup{width=0.9\textwidth}
\caption{%
The unknown fields for porous-elastic-plastic solids at fracture. The boundary $\partial \B$ is decomposed into Dirichlet and Neumann parts for the \textit{displacement} $\partial \B_\Bu \cup \partial \B_\Bt$, the \textit{fluid flux} $\partial \B_\nIh \cup \partial \B_\mu$ and the \textit{fracture phase-field} $\partial \B_d \cup \partial \B_\frakk$. For the fracture phase-field, zero Neuman boundary conditions are assumed. Here $c$ is a constant depending on the model formulation.
}%
\label{fig:fields}%
\end{figure}%

\paragraph{Fracture phase field.}

As mentioned above, we will model cracks and their evolution based on the phase-field approach to fracture. The fracture phase field $d$ is thus formally introduced as
\begin{equation}
d  : \left\{
\begin{array}{ll}
\B \times \calT \rightarrow [0,1] \\
(\Bx,t) \mapsto d(\Bx,t).
\end{array}
\right.
\label{eq:damageunknown}
\end{equation}
It denotes with $d=0$ an intact state and with $d=1$ a broken state of the material. The phase field is used to approximate a sharp crack interface $\Gamma$ in a diffuse manner. This results in the definition of a regularized crack surface $\Gamma_l$ in terms of a crack-surface density $\gamma$ and a corresponding length-scale parameter $l$ given by
\begin{equation}
\Gamma \approx \Gamma_l(d) := \int_\B \gamma(d,\nabla d) ~\text{d}V
\WITH
\gamma(d,\nabla d) := \tfrac{1}{2l}d^2+\tfrac{l}{2}|\nabla d|^2,
\label{eq:cracksurface}
\end{equation}
Note that the sharp crack surface is recovered for vanishing length-scale parameter ($l \rightarrow 0 \Rightarrow  \Gamma_l \rightarrow \Gamma$). Here we follow the notation of \cite{miehe+welschinger+hofacker10a}.

\section{Variational formulation of fracturing porous-elastic-plastic solids}
\label{Section3}

Based on the previous section, we are able to introduce the primary fields for the description of fracturing porous-elastic-plastic solids as
\begin{equation}
  \BfrakU:=\{\Bu, \nIh ,d,\Bve^p,\alpha, m^p\} .
\end{equation}
In order to formulate a rate-type variational principle we define the rate of the primary fields as
\begin{equation}
 \dot\BfrakU:=\{\dot \Bu, \nIh, \dot d,\dot \Bve^p, \dot \alpha, \dot m^p\}.
\end{equation}
Furthermore, we identify the constitutive state of the model and its evolution as
\begin{equation}
  \BfrakC:=\{\Bve,m,\nIh,d,\nabla d,\Bve^p,\alpha,m^p\}
  \AND
  \dot\BfrakC:=\{\dot\Bve,\nIh,\div \nIh,\dot d,\dot{\nabla d},\dot\Bve^p,\dot\alpha,\dot m^p\} ,
\end{equation}
respectively.

\subsection{Formulation of the rate-type potential}

The general form of the rate-type potential is given by
\begin{equation}
\Pi(\dot\BfrakU; \BfrakU) := \frac{\text{d}}{\text{d}t}E(\dot\BfrakC;\BfrakC) + D(\dot\BfrakC) - P_\text{ext}(\dot\BfrakU) ,
\end{equation}
where $\tfrac{\text{d}}{\text{d}t}E(\dot\BfrakC;\BfrakC)$ is the rate of the energy, $D(\dot\BfrakC)$ is the dissipation potential and $P_\text{ext}(\dot\BfrakU)$ is the potential of the external loading. The rate of the energy is described in terms of the energy density $\psi(\BfrakC)$
\begin{equation}
  \frac{\text{d}}{\text{d}t} E(\dot\BfrakC;\BfrakC)
= \frac{\text{d}}{\text{d}t} \int_\B \psi(\BfrakC)~\text{d}V ,
\label{eq:rateenergy1}
\end{equation}
which, by application of the chain rule yields
\begin{equation}
\begin{split}
  \frac{\text{d}}{\text{d}t}E(\dot\BfrakC;\BfrakC)=\int_\B
   (\partial_\Bve\psi : \dot \Bve
   -\partial_m\psi \div\nIh
   +\partial_d\psi ~ \dot d
   +\partial_{\Bve^p}\psi : \dot \Bve^p
  +\partial_{\alpha}\psi ~ \dot \alpha
  +\partial_{m^p}\psi~ \dot m^p )~\text{d}V .
\end{split}
\label{eq:rateenergy}
\end{equation}
In the latter equation, we made use of the fluid mass balance \eqref{eq:globalfluidmass}$_2$. Similarly to the rate of energy the dissipation potential can be expressed in terms of a dissipation potential density  $\phi(\dot\BfrakC)$
\begin{equation}
  D(\dot\BfrakC) =\int_\B
  \phi(\dot\BfrakC)
  ~\text{d}V .
  \label{eq:dissipation}
\end{equation}
Combing the right-hand sides of equations \eqref{eq:rateenergy} and \eqref{eq:dissipation} yields the internal rate-potential density per unit volume
\begin{equation}
\begin{split}
  \pi(\dot\BfrakC;\BfrakC) = ~\partial_\Bve\psi : \dot \Bve
   -\partial_m\psi\div\nIh
  +\partial_d\psi ~ \dot d 
  +\partial_{\Bve^p}\psi : \dot \Bve^p
  +\partial_{m^p}\psi ~\dot m^p
  +\partial_{\alpha}\psi ~ \dot \alpha
  + \phi(\dot\BfrakC) .
  \end{split}
\end{equation}
so that 
\begin{equation}
   \Pi(\dot\BfrakU; \BfrakU)
:= \int_\B \pi(\dot\BfrakC;\BfrakC)~\text{d}V
-  P_\text{ext}(\dot\BfrakU) .
\end{equation}
The particular forms of the energy density $\psi(\BfrakC)$, the dissipation-potential density $\phi(\dot\BfrakC)$ and the potential of the external loading $P_\text{ext}(\dot\BfrakU)$ will be discussed in the following sections.

\subsubsection{Constitutive energy density} 

The energy density has two contributions, one from the solid $\psi_\text{solid}$ and one from the fluid $\psi_\text{fluid}$
\begin{equation}
\psi(\BfrakC) = \psi_\text{solid}(\Bve-\Bve^p,\alpha,d) + \psi_\text{fluid}(\Bve-\Bve^p,m-m^p) .
\end{equation}

\paragraph{Energy density of the solid phase.}

The energy density of the solid phase $\psi_\text{solid}$ can be decomposed into an effective elastic and a plastic part
\begin{equation}
  \psi_\text{bulk}(\Bve^e,\alpha,d)=\psi_\text{eff}(\Bve^e,d)+ \psi_\text{plast}(\alpha,d).
\end{equation}
Both parts depend on the fracture phase field $d$ by means of a degradation function $g(d)$. The degraded effective elastic energy reads
\begin{equation}
  \psi_\text{eff}(\Bve^e,d) = [g(d)+k]\psi^{0+}_\text{eff}(\Bve^e) + \psi^{0-}_\text{eff}(\Bve^e) ,
\end{equation}
where the superscript ``0'' indicates the energy density of the undamaged solid matrix. In the latter equation, we have decomposed the undamaged energy into tensile and compressive parts (indicated by the superscripts ``+'' and ``-'', respectively), from which only the tensile part is assumed to contribute to fracture propagation. The parameter $k \ll 1$ ensures the well posedness of the problem. In what follows, we assume that $g(d) =(1-d)^2$.

The undamaged effective energies $\psi^{0 \pm}_\text{eff}$ represent the behavior of the elastic matrix and take the simple quadratic forms
\begin{equation}
\psi^{0 \pm}_\text{eff}(\Bve^e)= \tfrac{\lambda}{2} \big \langle \tr{}(\Bve^e) \big\rangle^2_\pm+ G~\Bve_\pm^e\mkern-3mu:\Bve_\pm^e \WITH \Bve^e_\pm = {\textstyle \sum_{i=1,3}} \langle\lambda^e_i \rangle_\pm \Bn_i \otimes \Bn_i
\label{eq:elasticplasticenergy}
\end{equation}
in terms of the elastic strain $\Bve^e$ and the Lam\'e constants $\lambda \ge - \tfrac{2}{3} G$ and $G \ge 0$. The tensile and compressive strains $\Bve^e_\pm$ are given in terms of the eigenvalues of the strain tensor $\lambda^e_i$ and the ramp function $\langle x \rangle_\pm = (x \pm |x|)/2$.

The plastic energy density considers isotropic saturation-type hardening and takes the form
\begin{equation}
\psi_\text{plast}(\alpha,d) = [g(d)+k] \psi^0_\text{plast}(\alpha)  \WITH \psi^0_\text{plast}(\alpha) = 
\tfrac{h}{2}\alpha^2+\sigma_y\big[\alpha+\tfrac{1}{\omega}\exp(-\omega\alpha)-\tfrac{1}{\omega}\big] .
\label{eq:elasticplasticenergy}
\end{equation}
Here, $h\ge0$ is the hardening modulus, $\sigma_y$ is the saturated yield shift and $\omega$ is a saturation parameter, see
\cite{kienle+aldakheel+keip19}. The resulting hardening function $\beta(\alpha,d)=\partial_\alpha\psi_\text{plast}(\alpha,d)$ is shown in Figure~\ref{fig:yield}.
\paragraph{Energy density of the fluid phase.} 
Based on Biot's theory of consolidation
\citep{biot41}
the fluid energy density is chosen as
\begin{equation}
\psi_\text{fluid}(\Bve^e,m^e)=\tfrac{M}{2}\big[b \tr(\Bve^e)-
\tfrac{m^e}{\rho_f}\big]^2 ,
\label{eq:fluidenergy}
\end{equation}
where $M$ is Biot's modulus, $b$ is Biot's coefficient and $\rho_f$ is the fluid density. It satisfies the following conditions in terms of the fluid pressure $p$
\begin{equation}
  p : = -\tfrac{1}{b}\partial_{\tr \Bve}\psi_\text{fluid} = \rho_f \partial_m \psi_\text{fluid} .
  \label{eq:fluidpressure}
\end{equation}
We refer to
\cite{miehe+mauthe+teichtmeister15}
for a more detailed discussion on the construction of $\psi_\text{fluid}$. Note that according to \eqref{eq:fluidenergy} and \eqref{eq:fluidpressure} the fluid pressure depends only on the elastic quantities $\Bve^e$ and $m^e$.

\subsubsection{Dissipation potential density}

Similar to the energy density the dissipation potential density can be additively decomposed into individual contributions, here associated with dissipative effects arising from fluid flow, fracture evolution and plastic deformations. The dissipation potential density is formally given by
\begin{equation}
\phi(\dot \BfrakC)=\phi_{\text{fluid}}(\nIh)+\phi_\text{frac}(\dot d,\dot{\nabla d})+\phi_\text{plast}(\dot \Bve^p,\dot \alpha,\dot m^p) .
\end{equation}

\paragraph{Dissipation-potential density for fluid flow.}

We follow
\cite{miehe+mauthe+teichtmeister15}
and employ a dissipation-potential density in terms of the fluid flux $\nIh$ given by
\begin{equation}
\phi_{\text{fluid}}(\nIh) = \half \BK^{-1}:(\nIh\otimes \nIh)
\end{equation}
The permeability tensor $\BK$ at a given state $\{\Bve, d,\nabla d\}$ is defined as
\begin{equation}
  \BK(\Bve, d,\nabla d)
= [ 1 - f(d) ] \BK_0
+ f(d) \BK_{\text{frac}}(\Bve,\nabla d) ,
\end{equation}
where $\BK_0$ is the permeability tensor of the undamaged bulk and $\BK_{\text{frac}}$ is the permeability tensor within a crack. Clearly, the function $f(d)$ acts as an interpolation function between intact and fully damaged states of the material. In what follows, we select $f(d)=d^\epsilon$, where $\epsilon$ is an interpolation parameter.

While the permeability tensor of the undamaged bulk can be formulated in an isotropic manner based on the spatial permeability $K$ as $\BK_0 = \rho_f^2 K \Bone$, the permeability tensor within a crack can be expressed in terms of the fracture-opening function $w(\Bve,\nabla d)$
\begin{equation}
\BK_{\text{frac}}(\Bve,\nabla d)=\rho_f^2\tfrac{w^2(\Bve,\nabla d)}{12\eta_f}\big(\Bone-\Bn\otimes\Bn\big) \WITH \Bn:=\tfrac{\nabla d}{|\nabla d|} ,
\label{eq:crackpermablity}
\end{equation}
which is anisotropic in nature. In the above definition, $\eta_f$ is the fluid's dynamic viscosity. Note that the representation of the permeability in \eqref{eq:crackpermablity} is derived based on the lubrication theory and relates to Poiseuille-type flow within the fractures. In this context, the fracture-opening function is given as
\begin{equation}
 w(\Bve, \nabla d)= (\Bn\cdot\Bve\cdot\Bn) L_{\perp} ,
\end{equation}
where $L_{\perp}$ is the length of a line element that is perpendicular to the crack. In a finite-element representation this can be identified as the element size $h^e$.

\begin{figure}
\centering
{
\small
\psfrag{l} [c][c] {$L_\perp$}
\psfrag{w} [c][c] {$w$}
\psfrag{a} [c][c] {a)}
\psfrag{b} [c][c] {b)}
\includegraphics[width=0.9\textwidth]{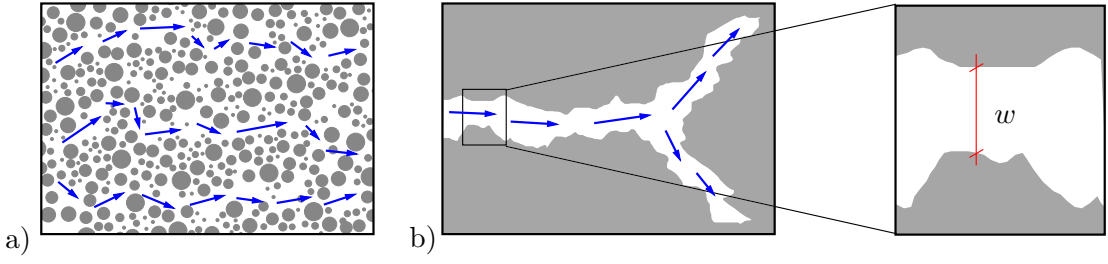}
}
\vspace{1mm}%
\captionsetup{width=0.9\textwidth}
\caption{%
Schematic representation of a) fluid flow in porous medium according to
Dracy's law and b) within developing fractures according to
Poiseuille--type law with fracture opening $w$.
}%
\label{fig:flow}%
\end{figure}%

\paragraph{Dissipation-potential density for fracture evolution.}

The dissipation-potential density associated with fracturing accounts for the dissipation of an evolving crack surface. It preserves thermodynamical consistency by ensuring a local irreversibility condition of the fracture phase field $\dot d \ge 0$, which can be achieved by introducing the indicator function $I(\dot d)$.

For general Griffith-type fracturing the dissipation-potential density reads
\begin{equation}
  \phi^{g_c}_\text{frac}(\dot d,\dot{\nabla d}) := \frac{\text{d}}{\text{d}t}\Big\{g_c\gamma(d,\nabla d)\Big\} + I(\dot d) ,
  \label{eq:fracdiss1}
\end{equation}
where $g_c$ is Griffith's critical energy-release rate. The indicator function $I(\dot d)$ is given as
\begin{equation}
I(\dot d) := \left\{
\begin{array}{ll}
0  & \text{for~} \dot d \ge 0,\\
\infty &\text{otherwise}.
\end{array}
\right.
\label{eq:indicatorfunction}
\end{equation}
Note that the fracture evolution described by \eqref{eq:fracdiss1} does not include any threshold value for the fracture evolution. In that case, damage occurs already at very small load levels. Hence, to clearly separate the material response during loading in elastic and plastic behavior as well as subsequent fracturing, the following dissipation-potential density will be used
\begin{equation}
  \phi_\text{frac}(\dot d,\dot{\nabla d}) := \frac{\text{d}}{\text{d}t}\Big\{[1-g(d)]\psi_c+2\psi_cl\gamma(d,\nabla d)\Big\} + I(\dot d) ,
\end{equation}
where $\psi_c$ is a threshold value.

\paragraph{Dissipation-potential density for plastic response.} 

The dissipation-potential density accounting for plastic behavior can be derived based on the principle of maximum dissipation. For that, the thermodynamical driving forces for the plastic strain, the hardening and the change of plastic fluid content need to be specified. Using the second law of thermodynamics (Clausius--Planck inequality) yields the driving forces as thermodynamical duals of $\Bve^p$, $\alpha$ and $m^p$ as
\begin{equation}
  \begin{split}
  \Bsigma&:= -\partial_{\Bve^p}\psi
           = -\partial_{\Bve^p}\psi_\text{eff} - \partial_{\Bve^p}\psi_\text{fluid}
           = \Bsigma_\text{eff} - b \rho_f \mu \Bone  \\
  \beta  &:=-\partial_{\alpha}\psi=-\partial_{\alpha}\psi_\text{plast} \\
  \mu    &:= -\partial_{m^p}\psi= -\partial_{m^p}\psi_\text{fluid} .
\end{split}
\label{eq:duals}
\end{equation}
We denote the set of diving forces as $\BfrakF=\{\Bsigma,\beta,\mu\}$, wherein $\Bsigma$ is the Cauchy stress, $\beta$ is the hardening function and $\mu$ is the fluid potential. Based on that, we construct the dissipation-potential density $\phi_\text{plast}$, which is formulated in terms of the constrained optimization problem
\begin{equation}
  \phi_\text{plast}(\dot \Bve^p,\dot \alpha,\dot m^p):=
  \underset{\BfrakF\in \nIE}{\sup}
  [{\Bsigma:\dot{\Bve^p} + \beta\dot \alpha+ \mu \dot m^p}]~.
\label{eq:constrainedoptimization}
\end{equation}
In the above definition, the plastic driving forces are constrained to lie within the elastic domain $\nIE:=\{\BfrakF| f^p(\BfrakF)\le 0\}$, which is characterized by the yield function $f^p$. The latter will be specified at a later stage. By introducing the Lagrange multiplier $\lambda^p$ the constrained optimization in \eqref{eq:constrainedoptimization} can be rewritten as
\begin{equation}
  \phi_\text{plast}(\dot \Bve^p,\dot \alpha,\dot m^p):=
  \underset{\BfrakF}{\sup}
  \underset{\lambda^p \ge 0}{\phantom{\text{p}}\inf \phantom{\text{p}}}
  \big[{\Bsigma:\dot{\Bve^p} + \beta\dot \alpha+ \mu \dot m^p - \lambda^p f^p(\BfrakF)}\big],
  \label{eq:constrainedoptimizationlagrange}
\end{equation}
The latter representations are related to the rate-independent elastic-plastic material response, which leads to non-smooth evolution of plastic deformations. The non-smoothness can be relaxed by introducing viscous regularization, which yields the modified dissipation-potential density
\begin{equation}
  \phi_\text{plast}(\dot \Bve^p,\dot \alpha,\dot m^p):=
  \underset{\BfrakF}{\sup}
  ~\big[\underbrace{\Bsigma:\dot{\Bve^p} + \beta\dot \alpha+ \mu \dot m^p -\tfrac{1}{2\eta_p}\langle f^p(\BfrakF)\rangle_+^2}_{\frakD_\text{plast}}\big] .
  \label{eq:constrainedoptimizationviscous}
\end{equation}
From a mathematical point of view this density is obtained from an optimization procedure with side condition, the latter of which is enforced by a quadratic penalty term. The penalty parameter is given by the plastic viscosity $\eta_p$.

\begin{figure}
  \begin{center}
    \small

    \psfrag{a}   [c][c]{a)}
    \psfrag{bb}  [c][c]{b)}
    \psfrag{phi} [c][c]{$f^p=0$}
    \psfrag{elastic domain} [c][c]{\footnotesize \, \, \, \, elastic domain}
    \psfrag{m}   [rc][c]{$M_\phi$}
    \psfrag{1}   [tc][c]{$1$}
    \psfrag{p1}  [lc][rc]{$~s_{\text{max}}^{\phantom{*}}$}
    \psfrag{p2}  [lc][c]{$s_{\text{max}}^*$}
    \psfrag{b}   [lb][c]{$\beta(\alpha)$}
    \psfrag{x}   [c][b]{$\tfrac{1}{3}\tr\Bsigma_\text{eff}$}
    \psfrag{y}   [lc][rc]{$\sqrt{\tfrac{2}{3}}|\dev [\Bsigma_\text{eff}]|$}

    \psfrag{h}   [c][c]{$h$}
    \psfrag{t}   [rc][c]{$\sigma_y$}
    \psfrag{xt}  [ct][c]{$\alpha^*$}
    \psfrag{w}   [c][c]{$\omega$}
    \psfrag{xx}  [ct][c]{$\alpha\phantom{^*}$}
    \psfrag{yy}  [rt][c]{$\beta(\alpha)$}
    \psfrag{psi} [c][c]{}
    \includegraphics[width = 0.95\textwidth]{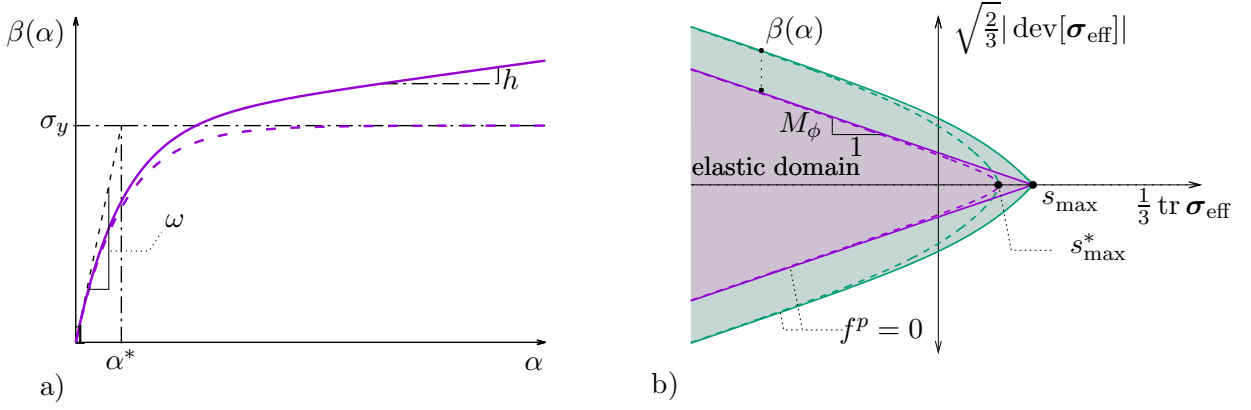}
    \captionsetup{width=0.9\textwidth}
    \caption{Hardening function $\beta(\alpha;d=0)$ with linear hardening $h\ne0$ (solid line) and without linear hardening $h=0$ (dashed line), where $\alpha^*=(\tfrac{h}{\sigma_y}+\omega)^{-1}$ in a). In b) the yield function in two-dimensional hydrostatic-deviatoric plane  with and without regularization (dashed/solid lines) where  $s^*_\text{max} = s_\text{max} - \sqrt{\tfrac{2}{3}} q_1$. In green with hardening and no hardening in purple.
    }
    \label{fig:yield}
  \end{center}
\end{figure}

In the present work it is assumed that the yield function $f^p_\text{solid}$ describes the plastic response of the drained solid matrix. It is thus expressed in terms of the effective stresses $\Bsigma_\text{eff}$ and the hardening function $\beta$. It should be equal to the yield function $f^p$ of the undrained bulk, which can be expressed in terms of the total stress $\Bsigma$, the hardening function $\beta$ and the fluid potential $\mu$
\begin{equation}
  f^p_\text{solid}(\Bsigma_\text{eff},\beta)
= f^p(\Bsigma,\beta,\mu).
\label{eq:yieldsk}
\end{equation}
A similar assumption has been made in
\cite{armero99}.
The starting point for the derivation of an appropriate yield function for the presented model is given by the yield function for frictional materials presented in
\cite{kienle+aldakheel+keip19}, see also
\cite{vermeer+deborst84}
and
\cite{lambrecht+miehe01}.
Making use of equations \eqref{eq:duals}$_1$ and \eqref{eq:yieldsk} yields
\begin{equation}
f^p(\Bsigma,\beta,\mu) = \sqrt{\tfrac{3}{2}}\sqrt{|\dev [\Bsigma]|^2 
+ M_\phi^2 q_1^2} - M_\phi (s_{\text{max}}-\tfrac{1}{3}\tr\Bsigma - b\rho_f\mu) + \beta M_h(\Bsigma,\mu) ,
\label{eq:yield1}
\end{equation}
where $M_\phi$ is related to the friction angle, $q_1$ is a parameter related to the regularization of the tip of the yield surface and $s_\text{max}$ is related to the cohesion of the material. The hardening response is limited to friction hardening by choosing the material function $M_h(\Bsigma,\mu)$ as
\begin{equation}
  M_h(\Bsigma,\mu) = 1 - \sqrt{\tfrac{3}{2}}q_1\exp \big[ \tfrac{1}{3}\tr\Bsigma + b\rho_f\mu - s_{max} \big] .
\label{eq:yield2}
\end{equation}
By inserting equation \eqref{eq:duals}$_1$ in equation \eqref{eq:yield1} the yield function $f^p_\text{solid}(\Bsigma_\text{eff},\beta)$ can be recovered. The yield function in terms of the effective stress $\Bsigma_\text{eff}$ is visualized in Figure~\ref{fig:yield}.

\subsubsection{Potential of external loading}

The external loading in form of mechanical tractions and fluid potential is formulated as
\begin{equation}
P_\text{ext}(\dot \BfrakU)= \int_{\partial \B_\Bt}\bar \Bt \cdot\dot \Bu~\text{d}A - \int_{\partial \B_\mu}\bar \mu \nIh \cdot \Bn ~\text{d}A ,
\end{equation}
where $\bar\Bt$ is the mechanical traction vector applied on the traction boundary $\partial \B_\Bt$ of the domain $\B$. The fluid contribution is due to the fluid transport over the boundary $\partial \B_\mu$ of the domain $\B$, where the fluid potential $\bar \mu$ is applied.

\subsection{Minimization principle and mixed variational principle}

Based on the above introduced functions we can introduce a rate-type minimization principle that governs the boundary-value problems of porous-elastic-plastic solids at fracture
\begin{equation}
\boxed{
\dot \BfrakU^* = 
\arg~\big\{\underset{\dot \BfrakU \in \frakW}{\inf}~\Pi(\dot \BfrakU; \BfrakU) \big\} }
\label{eq:minimization}
\end{equation}
where $\frakW=\{\calW_{\dot \Bu}, \calW_\nIh, \calW_{\dot d}, \calW_{\dot \Bve^p}, \calW_{\dot \alpha}, \calW_{\dot m^p}\}$ is the set of admissible spaces corresponding to the set of the rate of unknowns $\dot \BfrakU$. The admissible spaces are given as
\begin{align}
  \calW_{\dot \Bu}    := \{\dot \Bu \in H^1(\B)|\dot \Bu = \dot{\bar \Bu} ~\text{on}~ \partial \B_\Bu\}, \quad
 &\calW_\nIh          := \{\nIh\in H(\div,\B)|\nIh \cdot \Bn = \bar h ~\text{on}~ \partial \B_\nIh\}\nonumber \\
  \calW_{\dot d}      := \{\dot d\in H^1(\B)\},\quad
 &\calW_{\dot \Bve^p} := \{\Bve^p\in L^2\}  \\
  \calW_{\dot \alpha} := \{\dot \alpha \in L^2\},\quad
 &\calW_{\dot m^p}    := \{\dot m^p \in L^2\} . \nonumber
\end{align}

Combining the global minimization principle \eqref{eq:minimization} with the local maximization principle in \eqref{eq:constrainedoptimizationviscous} yields the mixed variational principle
\begin{equation}
\boxed{
\{\dot \BfrakU^*,\BfrakF^*\} = 
\arg~\bigg\{\underset{\dot \BfrakU \in \frakW}{\inf}~\underset{\BfrakF \in L^2}{\sup}~\int_\B\pi^\star (\dot \BfrakC,\BfrakF;\BfrakC)~\text{d}V - P_\text{ext} (\dot \BfrakU) \bigg\}
}
\label{eq:mixedvar}
\end{equation}
where $\pi^\star(\dot \BfrakC, \BfrakF;\BfrakC)$ is the mixed potential density. It reads
\begin{equation}
  \pi^\star(\dot \BfrakC, \BfrakF;\BfrakC) =
  \tfrac{\text{d}}{\text{d}t}\psi(\BfrakC)
  +\frakD_\text{plast}(\dot \Bve^p,\dot \alpha,\dot m^p,\BfrakF) +\phi_{\text{fluid}}(\nIh)+\phi_\text{frac}(\dot d,\dot{\nabla d}).
    \label{eq:mixedpot}
\end{equation}
Performing the variation of \eqref{eq:mixedpot} at a fixed state $\BfrakC$, we obtain the Euler equations of the mixed variational principle \eqref{eq:mixedvar} for the global unknowns as
\begin{equation}
\boxed{
\begin{aligned}
  \div[\partial_\Bve\psi] &= \phantom{\ni} \Bzero  \phantom{0}              ~~\text{in}\quad \B \\
  \nabla[\partial_m\psi]+\partial_\nIh \phi &= \phantom{\ni} \Bzero \phantom{0} ~~\text{in}\quad \B \\
  \partial_d\psi+\partial_{\dot d}\phi - \div[\partial_{\dot{\nabla d}}\phi]&\ni \phantom{=} 0\phantom{\Bzero}  ~~ \text{in}\quad \B\\
  -\partial_m\psi+\bar\mu &=\phantom{\ni} 0  \phantom{\Bzero}                     ~~\text{on}\quad \partial\B_\mu \\
  \partial_\Bve\psi \cdot \Bn - \bar \Bt &= \phantom{\ni}\Bzero \phantom{0}         ~~\text{on}\quad \partial\B_\Bt \\
  \partial_{\dot{\nabla d}} \phi \cdot \Bn &= \phantom{\ni}0    \phantom{\Bzero}        ~~\text{on}\quad \partial\B_\frakk
\end{aligned}
}
\label{eq:globaleuler}
\end{equation}
and the local unknowns with the corresponding thermodynamic duals as
\begin{equation}
\boxed{
\begin{aligned}
  \partial_{\Bve^p}\psi+\Bsigma &= \Bzero   \phantom{0}      \quad\text{in}\quad \B\\
  \partial_\alpha\psi+\beta &= 0   \phantom{\Bzero}                        \quad\text{in}\quad \B\\
  \partial_{m^p}\psi +\mu &= 0    \phantom{\Bzero}                        \quad\text{in}\quad \B\\
  \dot \Bve^p - \lambda^p_v \, \partial_\Bsigma f^p &= \Bzero  \phantom{0}         \quad\text{in}\quad \B\\
  \dot \alpha - \lambda^p_v \, \partial_\beta f^p &= 0 \phantom{\Bzero}       \quad\text{in}\quad \B\\
  \dot m^p -    \lambda^p_v \, \partial_\mu f^p &= 0  \phantom{\Bzero}    \quad\text{in}\quad \B
\end{aligned}
}
\label{eq:localeuler}
\end{equation}
The equations in \eqref{eq:globaleuler} represent the global balance laws and the corresponding Neumann boundary conditions. In \eqref{eq:localeuler} the definition of the thermodynamic duals of the local fields as well as their evolution equations are given. In the latter, we introduced the visco-plastic multiplier $\lambda^p_v := \tfrac{1}{\eta_p}\langle f^p \rangle_+$.

\paragraph{Condensation of local variables.}

The set of the rate of the primary fields can be split into a local and global part. The former is given as $\dot \BfrakU_\text{l} = \{\dot \Bve^p, \dot \alpha, \dot m^p\}$ and the latter arises as $\dot \BfrakU_\text{g}=\dot \BfrakU \setminus \dot \BfrakU_\text{l} = \{\dot \Bu, \dot d , \nIh\}$. The set of the rate of the local fields $\dot \BfrakU_\text{l}$ is governed by the mixed variational principle
\begin{equation}
  \boxed{
  \{\dot \BfrakU_\text{l}^*,\BfrakF^* \}=
  \arg \big\{ \underset{\dot \BfrakU_\text{l}}{\phantom p \inf \phantom p}~\underset{\BfrakF}{\sup}
  ~ \pi^\star(\dot \BfrakC, \BfrakF;\BfrakC)\big\} }
  \label{eq:localmixprinciple}
\end{equation}
In order to obtain a solution for $\dot \BfrakU_\text{g}$ the reduced potential density is introduced as
\begin{equation}
  \pi_\text{red}^\star(\dot \BfrakC_\text{red};  \BfrakC_\text{red})=
   \underset{\dot \BfrakU_\text{l}}{\phantom p \inf \phantom p}~\underset{\BfrakF}{\sup}
  ~\pi^\star(\dot \BfrakC, \BfrakF;\BfrakC) ,
\end{equation}
where the reduced constitutive state $\BfrakC_\text{red}=\{\Bve,m,d,\nabla d\}$ and its evolution $\dot \BfrakC_\text{red}=\{\dot \Bve,\nIh, \div \nIh,\dot d,\dot {\nabla d}\}$ are introduced. The set of the rate of the global fields is given by the following minimization principle
\begin{equation}
  \boxed{
  \dot \BfrakU^*_\text{g} =
  \{\dot \Bu^*, \dot d^* , \nIh^*\}
  =\arg\bigg\{
  ~\underset{\dot \Bu \in \calW_{\dot \Bu}}{\inf}
  ~\underset{\dot d \in \calW_{\dot d}}{\inf}
  ~\underset{\nIh \in \calW_\nIh}{\inf}
  \int_\B\pi_\text{red}^\star(\dot \BfrakC_\text{red}; \BfrakC_\text{red})~\text{d}V - P_\text{ext} ( \dot \BfrakU ) \bigg \} }
\end{equation}

\subsection{Modification of the fracture driving force} 

In this section a closer look at equation \eqref{eq:globaleuler}$_3$ is taken. Inserting the definitions of the energy density and the dissipation potential density yields
\begin{equation}
  -2(1-d) [\psi^0_\text{eff}+ \psi^0_\text{plast}-\psi_c]+2\psi_c(d - l^2 \Delta d) + \partial_{\dot d}I \ni 0 .
  \label{eq:fracevol1}
\end{equation}
By introducing the crack driving history field $\calH$ this is modified to
\begin{equation}
  -2(1-d) \calH + 2 \psi_c(d - l^2 \Delta d) = 0
  \WITH
  \calH := \underset{s\in [0,t]}{\max} \langle \psi^0_\text{eff}+ \psi^0_\text{plast}-\psi_c \rangle_+ .
  \label{eq:fracevol2}
\end{equation}
This follows the notation for brittle fracture in \cite{miehe+hofacker+welschinger10} and for ductile fracture in \cite{miehe+hofacker+schaenzel+aldakheel15, miehe+aldakheel+raina16}.

\begin{figure}
  \begin{center}
    \footnotesize

    \psfrag{c}   [c][c]{a)}
    \psfrag{a}   [c][c]{b)}
    \psfrag{b}   [c][c]{c)}
    \psfrag{phi} [l][c]{$\psi^0_\text{eff}$}
    \psfrag{p1}  [l][l]{$\psi^0_\text{plast}=0$}
    \psfrag{p2}  [l][l]{$\psi^0_\text{plast}$}
    \psfrag{p3}  [l][l]{$w_\text{plast}$}
    \psfrag{x}   [rc][c]{$\Bve$}
    \psfrag{y}   [l][l]{$\Bsigma^0_\text{eff}$}

    \includegraphics[width = 0.95\textwidth]{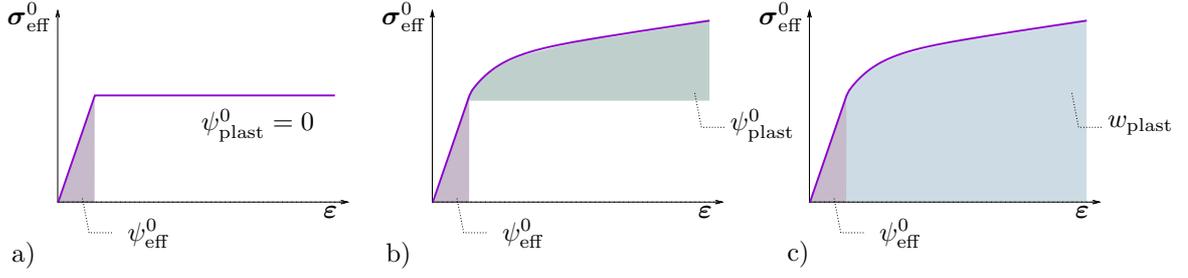}
    \captionsetup{width=0.9\textwidth}
    \caption{Visualization of the different energy contributions to the crack driving history field $\calH$. For the representation of $\calH$ in \eqref{eq:fracevol2}$_2$ only the elastic energy and the energy arising from the hardening, as show in a) and b), will drive the crack. For ideal plasticity the plastic energy vanishes ($\psi^0_\text{plast}=0$) leading to pure brittle fracture, see a). By using \eqref{eq:crackdrive} the crack is driven by the elastic energy and the full plastic work, see c).
    }
    \label{fig:energy}
  \end{center}
\end{figure}

For plasticity models with a yield limit that is independent of the stress state it is possible to formulate a plastic energy density $\widetilde\psi^{0}_\text{plast}(\alpha)$, which contains not only the work of the hardening for the solid matrix but also the work of the ideal plastic deformation of the solid matrix. One such example is given by von-Mises plasticity, for which we can write
\begin{equation}
  \widetilde\psi^{0}_\text{plast}(\alpha)\overset{\text{von Mises}}{=}  w_\text{plast}
\WITH  w_\text{plast}\overset{\text{\eqref{eq:localmixprinciple}}}{=}{\int}\Bsigma^0_\text{eff}:\dot \Bve^p~\text{d}t
\label{eq:vonmisesenergy}
\end{equation}
Note that the undamaged effective stress $\Bsigma^0_\text{eff}$ acting on the solid matrix is used here. The plastic work $w_\text{plast}$ can alternatively be expressed in terms of the undamaged total stress $\Bsigma^0$ and the fluid potential $\mu$ as
\begin{equation}
w_\text{plast}=  {\int}\Bsigma^0_\text{eff}:\dot \Bve^p~\text{d}t=  {\int}\Bsigma^0:\dot \Bve^p+\mu\dot m^p~\text{d}t .
  \label{eq:plasticwork}
\end{equation}
Since the construction of an energy density  $\widetilde\psi^{0}_\text{plast}$ like \eqref{eq:vonmisesenergy}$_1$ is not possible for more complicated plasticity models such as the Drucker--Prager model, the crack driving history field in \eqref{eq:fracevol2}$_2$ is modified to
\begin{equation}
 \calH =\underset{s\in [0,t]}{\max} \langle \psi^0_\text{eff}+w_\text{plast}-\psi_c \rangle_+ .
  \label{eq:crackdrive}
\end{equation}
With the representation of the crack driving history field in \eqref{eq:crackdrive} it is possible to model ductile fracture evolution that is driven by the elastic and the ideal plastic deformation as well as the hardening. With the representation in \eqref{eq:fracevol2}$_2$ and the definition of the plastic energy in \eqref{eq:elasticplasticenergy}$_2$ it is only possible to model a ductile fracture evolution which is driven by the elastic deformation and the hardening.

\subsubsection{Relation between plastic strain and change of fluid content}

Based on the construction of the yield function, see \eqref{eq:yieldsk}, the evolution of the plastic strain \eqref{eq:localeuler}$_4$ can be reformulated as
\begin{equation}
  \dot \Bve^p=\lambda^p_v \, \partial_\Bsigma f^p =
  \lambda^p_v \,\partial_{\Bsigma_\text{eff}}f^p_\text{solid}.
  \label{eq:plasticevol2}
\end{equation}
Furthermore, reformulation of the evolution of the change of plastic fluid content \eqref{eq:localeuler}$_6$ yields
\begin{equation}
\dot m^p = \lambda^p_v \, \partial_\mu f^p = 
\lambda^p_v \, \rho_f b ~\tr (\partial_{\Bsigma_\text{eff}}f^p_\text{solid}).
\label{eq:plasticfulidevol2}
\end{equation}
Combining the above two equations gives a relation between the evolution of the plastic strain and the evolution of the change of the plastic fluid content
\begin{equation}
\dot m^p = \rho_f b \tr \dot \Bve^p .
\end{equation}
We observe that the change of the plastic fluid content depends exclusively on volumetric plastic deformation.

\section{Numerical Treatment}
\label{Section4}

\subsection{Incremental variational formulation}

The incremental version of the rate-type potential introduced in Section~\ref{Section3} is obtained by algorithmic time integration over a given time step ~$\tau=[t_n,t_{n+1})$. For a pure Dirichlet problem ($P_\text{ext}=0$) we arrive at
\begin{equation}
\Pi^\tau(\BfrakU,\BfrakF) = \text{Algo} \Big\{ \int_{t_n}^{t_{n+1}} \Pi(\dot\BfrakU, \BfrakF; \BfrakC) ~\text{d}t \Big\}= \int_\B\pi^{\star\tau}(\BfrakC,\BfrakF;\BfrakC_n) ~\text{d}V ,
\end{equation}
where $\pi^{\star\tau}(\BfrakC,\BfrakF;\BfrakC_n)$ is the incremental potential density. It is given in terms of the energy density $\psi$, the incremental fluid and fracture dissipation-potential density $\phi^\tau_{\text{fluid}}$ and $\phi^\tau_\text{frac}$, respectively, as well as the incremental dissipation density related to the visco-plastic behavior $\frakD^{\tau}_{\text{plast}}$ as
\begin{equation}
  \pi^{\star\tau}(\BfrakC,\BfrakF;\BfrakC_n) =
 \psi(\BfrakC)
 +\frakD^{\tau}_\text{plast}(\Bve^p,\alpha, m^p,\BfrakF;\Bve_n^p,\alpha_n, m_n^p) + \phi^\tau_{\text{fluid}}(\nIh)+ \phi^\tau_\text{frac}(d,\nabla d,d_n) .
\end{equation}
The individual dissipative contributions read
\begin{equation}
  \begin{split}
    \phi^\tau_{\text{fluid}} &= \tau \phi_{\text{fluid}}(\nIh)\\
    \phi^\tau_\text{frac}~ &= [1-g(d)] \psi_c +2 \psi_c l \gamma(d,\nabla d)
    + I^\tau(d,d_n) \\
    \frakD^{\tau}_\text{plast}&= \Bsigma:(\Bve^p-\Bve^p_{n}) + \beta (\alpha-\alpha_n) + \mu (m^p-m^p_n) + \tfrac{\tau}{2\eta_p}\langle f^p(\BfrakF)\rangle_+^2.
  \end{split}
\end{equation}
Furthermore the fluid mass balance \eqref{eq:globalfluidmass}$_2$ is satisfied by the implicit update
\begin{equation}
  m = m_n - \tau \div \nIh .
\end{equation}

\paragraph{Condensation of local variables.}

Similar as in the continuous problem, the set of primary fields can be decomposed into a local and global part. Again, the local fields are identified as $\BfrakU_\text{l} = \{\Bve^p, \alpha, m^p\}$ and the global fields as $\BfrakU_\text{g}=\BfrakU \setminus \BfrakU_\text{l} = \{\Bu, d , \nIh\}$. The local fields are governed by the mixed variational principle
\begin{equation}
  \boxed{
  \BfrakU^*_\text{l}=
  \arg \{ \underset{\BfrakU_\text{l}}{\phantom p \inf \phantom p}~\underset{\BfrakF}{\sup}
  ~\pi^{\star\tau}(\BfrakC, \BfrakF;\BfrakC_n)\}}
  \label{eq:mixvariationltime}
\end{equation}
Using the representation \eqref{eq:constrainedoptimizationviscous} the mixed variational principle \eqref{eq:mixvariationltime} leads to the following condition
\begin{equation}
  \partial_{\BfrakU_\text{l},\BfrakF}\pi^\star=
  \left[
\begin{array}{c}
\partial_{\Bve^p}\psi+\Bsigma \\
\partial_{\alpha}\psi+\beta \\
\partial_{m^p}\psi+\mu \\
\Bve^p-\Bve^p_{n}- \gamma^p_v \ \partial_\Bsigma f^{p}\\
\alpha-\alpha_n- \gamma^p_v \ \partial_\beta f^{p}\\
m^p-m^p_n - \gamma^p_v \ \partial_\mu f^p\\
\end{array}
\right] = \Bzero .
\label{eq:localstress}
\end{equation}
Here, $\gamma^p_v:=\tau \lambda^{p}_v=\tfrac{\tau}{\eta_p}\langle f^p \rangle_+$ is the incremental visco-plastic multiplier. The local system of equations in \eqref{eq:localstress} is solved via a general return mapping scheme summarized in Box~1.

\hspace*{-6mm}\fbox{
\hspace{-7.5mm}
\parbox{160mm}{   
\vspace{-4mm}    
\begin{itemize}
\item[{\bf 0.}] Get trial values
  $[\Bve^{e,\text{tr}};~m^{e,\text{tr}};~\alpha^\text{tr}]^T =  [(\Bve-\Bve_n^p);~(m-m_n^p);~\alpha_n]^T$
\vspace{-2mm}
\item[{\bf 1.}]Set initial values
  $[\Bve^e;~m^e;~\alpha]^T = [\Bve^{e,\text{tr}};~m^{e,\text{tr}};~\alpha^\text{tr}]$
and $\gamma^p_v = 0$.
\vspace{-2mm}
\item[{\bf 2.}]Compute derivatives of energy density and yield function\\
  \vspace{-2mm}
\begin{equation*}
\fbox{\ 
\parbox{11cm}{
\hspace{-5mm}
$
\begin{array}{l@{\hspace{5mm}}r@{\ =\ }lr@{\ =\ }l}
\psi(\Bve^e,m^e,\alpha;d)&
            \Bs &
                  \left[\begin{array}{l}
                           \psi_{,\Bve^e} \\
                           \psi_{,m^e} \\
                           \psi_{,\alpha}
                          \end{array} \right]
     &
            \IE &
                  \left[\begin{array}{lll}
                    \psi_{,\Bve^e\Bve^e} & \psi_{,\Bve^e m^e} &  \psi_{,\Bve^e \alpha}    \hspace{-1mm} \\
                    \psi_{,m^e \Bve^e} & \psi_{,m^e m^e} & \psi_{,m^e\alpha}             \hspace{-1mm}\\
                    \psi_{,\alpha \Bve^e} & \psi_{,\alpha m^e} & \psi_{,\alpha\alpha}\hspace{-1mm}
                  \end{array}
\right]
\\[8mm]
f^p{}(\Bsigma, \mu,\beta)&
            \bn &
                  \left[\begin{array}{c}
                    f^p{}_{,\Bsigma} \\
                    f^p{}_{,\mu} \\
                    f^p{}_{,\beta}
                  \end{array}  \right]
&
            \IF &
            \left[\begin{array}{lll}
              f^p{}_{,\Bsigma \Bsigma} &   f^p{}_{,\Bsigma \mu}  &   f^p{}_{,\Bsigma \beta}\hspace{-1mm}  \\
              f^p{}_{,\mu \Bsigma}     & f^p{}_{,\mu \mu}       &   f^p{}_{,\mu \beta}   \hspace{-1mm}  \\
              f^p{}_{,\beta \Bsigma}   & f^p{}_{,\beta \mu}      & f^p{}_{,\beta \beta}   \hspace{-1mm} 
            \end{array}
\right]
  
 \\[5mm] 
\end{array}$
}
\ }
\end{equation*}
\vspace{-6mm}
\item[{\bf 3.}] Check for yielding. If yielding do a local Newton iteration\\[-9.5mm]

\hspace{-6mm}
\begin{minipage}[t][][c]{\textwidth}
\begin{algorithm}[H]
\linespread{1.2}\selectfont
\SetAlgoVlined
\SetKwBlock{Compute}{compute}{}
\SetKwBlock{Check}{check}{}
\SetKwBlock{Update}{update}{}
 \eIf(\tcp*[h]{elastic step}){$f^\text{p}{} < 0$}
  {set $[\Bsigma;~\psi_{,m} ]^T = [\psi_{,\Bve^e};~\psi_{,m^e}]^T$ and $\IE^\text{ep} =\IE$ \\
  \Return}
  (\tcp*[h]{plastic step})
  {
\Compute(residual vector){
$\br{} :=
\left[( \Bve^p_n - \Bve^p );~(m^p_n-m^p);~(\alpha_n-\alpha) \right]^T
+ \gamma^p_v \, \Bn$}
\Check(if local Newton is converged){ 
$\mathrm{if\ }
[\; \sqrt{ \br^{T}\br
+ \; [\, f^p{} - \frac{\eta_\text{p}}{\tau} \gamma^p_v \, ]^2
} < tol \; ]
\ \hbox{go to {\bf 4.}}
$}
\Compute(incremental plastic parameter){
$\increment \gamma^p_v = { 1 \over C} [\; f^p{} -
\bn^{T} \IX \, \br \; ]
~~ \hbox{with} ~~
C:= \bn^{T}\IX\,\Bn   + \frac{\eta^p}{\tau}
~~ \hbox{and} ~~ \IX := [\; \IE^{-1}
+ \gamma^p_v \IF \; ]^{-1}$}
\Compute(incremental strains, plastic fluid content and hardening variable){
$\left[\increment \Bve^p;~\increment m^p; ~\increment \alpha \right]^T
= - \IE^{-1} \IX \; [\;
    \br + \increment \gamma^p_v \, \bn \; ]
$}
\Update(plastic quantities){
$[\Bve^p;~m^p~;\alpha;~\gamma^p_v]^T ~\Leftarrow ~[\Bve^p;~m^p;~\alpha;~\gamma^p_v]^T ~ + ~ [\increment \Bve^p;~ \increment m^p;~ \increment \alpha;~ \increment \gamma^p_v]^T$
}
\Update(elastic quantities){
$[\Bve^e;~m^e]^T =  [(\Bve-\Bve^p);~(m-m^p)]^T$
} 
go to {\bf 2}.
}
\end{algorithm}
\end{minipage}
\\[-4mm]
\item[{\bf 4.}]For plastic step: Obtain stresses and consistent moduli
\vspace{-3mm}
\ebn
[\Bsigma;~\psi_{,m} ]^T = [\psi_{,\Bve^e};~\psi_{,m^e}]^T
\quad \mbox{and} \quad
\IE{}^\text{ep} = \IX{}\ - \ { 1 \over C} \,
[\, \IX{} \cdot \Bn  \, ] 
 \otimes
[\, \Bn \cdot \IX  \, ]
\een
\end{itemize}
\vspace{-6mm}
} }

\vspace*{0.2cm}
\noindent \textbf{Box 1:} Return mapping and tangent moduli for poro-elasto-plasticity. It is based on the algorithm for elasto-plasticity in \cite{miehe98a}.

\paragraph{Reduced global problem.}

In order to obtain a solution for $\BfrakU_\text{g}$ the reduced potential density is introduced 
\begin{equation}
  \pi_\text{red}^{\star \tau}(\BfrakC_\text{red})=
   \underset{\BfrakU_\text{l}}{\phantom p \inf \phantom p}~\underset{\BfrakF}{\sup}
  ~\pi^{\star\tau}(\BfrakC, \BfrakF; \BfrakC_n) ,
\label{eq:timemoddensity}
\end{equation}
where the reduced constitutive state $\BfrakC_\text{red}=\{\Bve,\nIh, \div \nIh, d,{\nabla d}\}$ is introduced. The global fields are then given by the minimization principle
\begin{equation}
  \boxed{
  \BfrakU^*_\text{g}=\{\Bu^*, d^* , \nIh^*\}=
  \arg\bigg\{
  \underset{\Bu \in \calW_\Bu}{\inf}
  ~\underset{d \in \calW_d}{\inf}
  ~\underset{\nIh \in \calW_\nIh}{\inf}
  \int_\B\pi_\text{red}^{\star\tau}(\BfrakC_\text{red})~\text{d}V \bigg\},
}
  \label{eq:reducedprinciple}
\end{equation}
with the admissible spaces
\begin{align}
  \calW_{\Bu} & := \{\Bu \in H^1(\B)| \Bu = {\bar \Bu} ~\text{on}~ \partial \B_\Bu\}, \nonumber \\
  \calW_\nIh  & := \{\nIh\in H(\div,\B)|\nIh \cdot \Bn = \bar h ~\text{on}~ \partial \B_\nIh\}, \\
  \calW_{d}   & := \{d\in H^1(\B)\} \nonumber .
\end{align}

\subsection{Space-time-discrete finite-element formulation}

Considering a finite-element discretization $\calT^h(\B)$, the discrete state vector $\textbf{d}$ containing the discrete values of $\{\Bu,\nIh,d\}$ and the interpolation of the constitutive state $\BfrakC_\text{red}^h = \textbf{B} \textbf{d}$ the global minimization principle \eqref{eq:reducedprinciple} can be written as
\begin{equation}
  \textbf{d}^*=
  \arg\{
  \underset{\textbf{d}}{\inf}~
  \Pi^{\tau h}(\textbf{d}) \} \WITH  \Pi^{\tau h}(\textbf{d}) = \int_\B\pi_\text{red}^{\star\tau h}(\textbf{B}\textbf{d})~\text{d}V.
  \label{eq:algebaricmin}
\end{equation}
Here, we employ the shape functions from
\cite{raviart+thomas77}
for the interpolation of the fluid flux, see also
\cite{teichtmeister+mauthe+miehe19}.
The nodal displacement is interpolated by the shape functions of the enhanced-assumed-strain formulation, see
\cite{simo+rifai90}.
The interpolation of the nodal phase-field values is done by the standard Q$_1$-type shape functions.

The global algebraic minimization principle \eqref{eq:algebaricmin} leads to the following condition
\begin{equation}
  \textbf{R}:=\Pi_{,\textbf{d}}^{\tau h} = \int_\B\textbf{B}^T \textbf{S}~\text{d}V = \bzero
  \label{eq:globalcondition}
\end{equation}
where the generalized array $\textbf{S}$ is introduced. This array is defined as follows
\begin{equation}
\textbf{S}
:=
  \left[
\begin{array}{c}
\psi_{,\Bvarepsilon^h}\\
-\tau \psi_{,m^h}\\
  \tau\phi_{,\nIh^h}\\[1mm]\cline{1-1}\\[-3mm]
    -2(1-d)\calH+\phi^{\tau \prime}_{,d^h}\\
  \phi^{\tau}_{,\nabla d^h}\\
\end{array}
\right]
\label{eq:globalstress}
\end{equation}
Note that this array is not obtained by straightforward differentiation of \eqref{eq:algebaricmin}$_2$, i.e.\ we have in general $\textbf{S} \ne \partial_{\BfrakC^h_\text{red}}\pi^ {\star\tau h}_\text{red}$. This is due to the consideration of the history field $\cal H$ in \eqref{eq:fracevol2}. Above, we have introduced the $\phi^{\tau \prime} = \phi^{\tau}-I^\tau(d,d_n) $, where $I^\tau(d,d_n)$ is the time discrete version of the indicator function \eqref{eq:indicatorfunction}. The system of equations \eqref{eq:globalcondition} is solved by a Newton--Raphson-type iteration yielding
\begin{equation}
  \textbf{d} \leftarrow \textbf{d} - \textbf{K}^{-1}\textbf{R}  \WITH \textbf{K} := \Pi_{,\textbf{d}\textbf{d}}^{\tau h} = \int_\B \textbf{B}^T\textbf{C} \BB~\text{d}V .
\end{equation}
The generalized tangent array $\textbf{C}$ is given as
\begin{equation}
  \textbf{C}:=\partial_{\BfrakC^h_\text{red}}\textbf{S}=
  \left[
     \begin{array}{ccccc}
      \\[-3.7mm]
\cline{1-3}
\multicolumn{1}{|c}{ \nIE^\text{ep}_{\Bve \Bve}}   & -\tau \nIE^\text{ep}_{\Bve m} & \multicolumn{1}{c|}{0} & \cdot & \cdot \\
\multicolumn{1}{|c}{-\tau \nIE^\text{ep}_{m \Bve}}&  \tau^2\nIE^\text{ep}_{m m}           & \multicolumn{1}{c|}{0} & \cdot & \cdot \\
 \multicolumn{1}{|c}{ 0 }                         &                0           &  \multicolumn{1}{c|}{ \tau \phi_{,\nIh^h \nIh^h}} & 0  & 0 \\
 \cline{1-3} \\[-4mm]
 \cline{4-5} 
 \cdot &\cdot & 0 &\multicolumn{1}{|c}{ 2 \calH + \phi^{\tau \prime}_{,d^h d^h}}& \multicolumn{1}{c|}{ 0} \\
 \cdot & \cdot& 0 &\multicolumn{1}{|c}{ 0}& \multicolumn{1}{c|}{\phi^{\tau}_{,\nabla d^h \nabla d^h}}\\ \cline{4-5} \\[-3mm]
  \end{array}
\right] .
\label{eq:gloabltangent}
\end{equation}
Here "$0$" indicates that the corresponding derivative does not exists. The derivatives at the slots labeled with "$\cdot$" indeed exist but are not needed due to the modification of the fracture driving force \eqref{eq:fracevol2} as consequence of an operator split. The latter leads to a decoupling of the related fields so that the two boxed sub-blocks in \eqref{eq:gloabltangent} can be treated in separate solution steps. The so called one-pass solution strategy is utilized here
\citep{miehe+hofacker+welschinger10,miehe+hofacker+schaenzel+aldakheel15}.
This means that the displacement and flux is updated first and then the fracture phase-field is updated. This might underestimate the speed of the fracture evolution but can be controlled by the choice of an appropriated time step size $\tau$
\citep{kienle+aldakheel+keip19}.

\section{Numerical Examples}
\label{Section5}

In the following we present a sequence of numerical examples that demonstrate the capabilities of the model formulation. The examples start with a test of a porous-elastic-plastic medium that is surcharged with a rigid footing leading to the creation of shear bands. Furthermore, we analyze the effect of different driving forces on porous-elastic-plastic fracture evolution. The latter analysis is extended to the comparison of porous-elastic and porous-elastic-plastic materials response leading to the evolution of Hydraulically induced fractures.

\subsection{Rigid-footing test on porous-elastic-plastic medium} 

\begin{table}
\captionsetup{width=0.9\textwidth}
\caption{Material parameters for rigid footing test.}
\centering
\footnotesize
\renewcommand{\arraystretch}{1.3}
\begin{tabular}{|lr@{\ =\ }lc|lr@{\ =\ }l@{\hspace{1.5mm}}c|}
\hline
Lam\'e parameter      & $\lambda$ & 180.0 & GN/m$^2$ & slope yield function   & $M_\phi$ & 0.6 & -- \\
Lam\'e parameter      & $G$     & 31.0 & GN/m$^2$  & position of peak      & $s_\text{max}$  & 4.0      & MN/m$^2$ \\ 
hardening modulus     & $h$       & 0.035 & MN/m$^2$ &     Biot's modulus & $M$  & 25.0 & GN/m$^2$ \\                             
saturated yield shift & $\sigma_y$& 0.1  & MN/m$^2$&    Biot's coefficient     & $b$   & 0.5 & -- \\                           
saturation parameter  & $\omega$  & 2.0  & --        &  fluid dyn.  viscosity   & $\eta_f$        & $1.0\cdot 10^{-3}$    & Ns/m$^2$ \\  
plastic viscosity     & $\eta_p$  & $5\cdot 10^{-6}$  & s &permeability      & $K$ &$9.8\cdot 10^{-12}$             & m$^3$s/kg \\
perturbation parameter & $q_1$   &  0.04 & MN/m$^2$ & fluid density &$\rho_f$ & 1000.0 & kg/m$^3$ \\
\hline
\end{tabular}
\label{tab:footing}
\end{table}

In the first example we consider a rigid footing test without fracture evolution. Our goal is to analyze the effect of plasticity as well as fluid flux and storage on the system's response. We thus take into account three different material types:
\vspace{-0.3cm}
\begin{itemize}
\setlength\itemsep{-0.2cm}
\item[i)]   drained elastic-plastic material,
\item[ii)]  undrained porous-elastic material,
\item[iii)] undrained porous-elastic-plastic material with different permeabilities.
\end{itemize}
\vspace{-0.3cm}
The drained elastic-plastic material is recovered by setting Biot's modulus and coefficient to zero ($M=0$, $b=0$). In order to model the undrained porous-elastic material, the yield limit is increased to a very high value by setting $s_\text{max}=1\cdot 10^4$~MN/m$^2$. For the undrained porous-elastic-plastic material all contributions of the model are active, hence no artificial choice of any material parameter is necessary. In order to analyze the effect of the permeability on the overall model response, we consider different magnitudes of permeabilities given by an original value $K$ as well as a reduced and an increased permeability ($K/5$ and $5K$, respectively). The chosen material parameters are listed in Table~\ref{tab:footing}.

\begin{figure}
\centering
{
  \small
\psfrag{u}     [c][c]{$\bar \Bu$}%
\psfrag{w1}    [c][c]{$a$}%
\psfrag{w22}   [c][c]{$a/2$}%
\psfrag{l2}    [c][c]{$H$}%
\psfrag{w}     [c][c]{$W$}%
\psfrag{w2}    [c][c]{$W/2$}%

\includegraphics[width=0.95\textwidth]{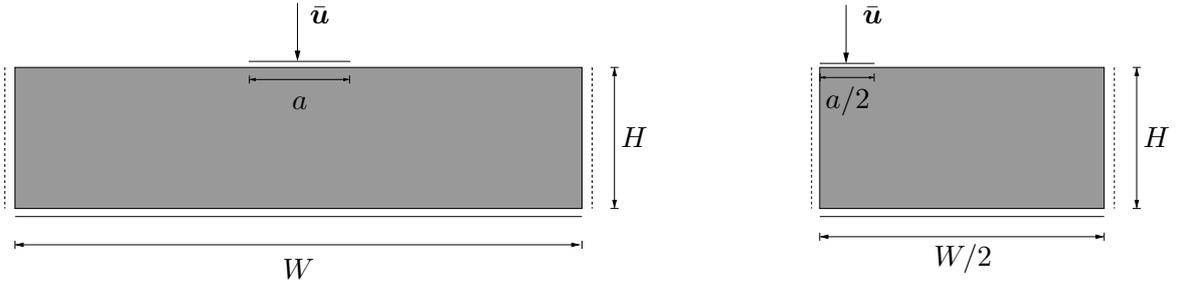}
}
\vspace{1mm}%
\captionsetup{width=0.9\textwidth}
\caption{Rigid footing test: Geometry and boundary conditions. Bottom is mechanically fixed, left and right edge is mechanically fixed in horizontal direction. Bottom, left and right is impermeable.
}
\label{fig:bound1}%
\end{figure}%

The geometry and boundary conditions are shown in Figure~\ref{fig:bound1}. Due to the symmetry of loading and geometry only one half of the specimen is discretized by $2,376$ quadrilateral Raviart--Thomas-type enhanced-assumed-strain elements. The dimensions are $H=4.758$~m, $W=23.088$~m and $a=4.587$~m. The loading increment is $\increment \bar \Bu=5\cdot 10^{-6}$~m. The loading is linearly increased until a total displacement of $\bar \Bu = 0.0023$~m is reached.

\begin{figure}%
\centering
\footnotesize
\psfrag{a}    [c][c] {a)}
\psfrag{b}    [c][c] {b)}
\psfrag{pmax} [l][l] {$9.4\cdot 10^{-3}$}
\psfrag{p}    [l][l] {$\alpha$}
\psfrag{pmin} [l][l] {$0$}
\includegraphics*[width = 0.75\textwidth]{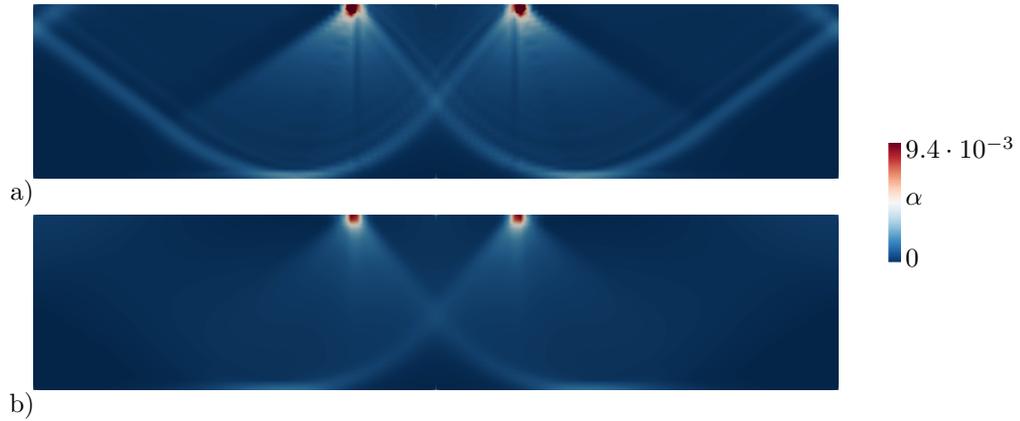}
\captionsetup{width=0.9\textwidth}
\caption{Rigid footing test: Distribution of the hardening variable $\alpha$ in drained elastic-plastic material in a) and undrained porous-elastic-plastic material in b).
}
\label{fig:flux_noflux}
\end{figure}

\begin{figure}
  \newcommand{\fsize}[0]{\Large}
  \centering
  \resizebox{!}{6cm}{
\begingroup
  \makeatletter
  \providecommand\color[2][]{%
    \GenericError{(gnuplot) \space\space\space\@spaces}{%
      Package color not loaded in conjunction with
      terminal option `colourtext'%
    }{See the gnuplot documentation for explanation.%
    }{Either use 'blacktext' in gnuplot or load the package
      color.sty in LaTeX.}%
    \renewcommand\color[2][]{}%
  }%
  \providecommand\includegraphics[2][]{%
    \GenericError{(gnuplot) \space\space\space\@spaces}{%
      Package graphicx or graphics not loaded%
    }{See the gnuplot documentation for explanation.%
    }{The gnuplot epslatex terminal needs graphicx.sty or graphics.sty.}%
    \renewcommand\includegraphics[2][]{}%
  }%
  \providecommand\rotatebox[2]{#2}%
  \@ifundefined{ifGPcolor}{%
    \newif\ifGPcolor
    \GPcolortrue
  }{}%
  \@ifundefined{ifGPblacktext}{%
    \newif\ifGPblacktext
    \GPblacktexttrue
  }{}%
  \let\gplgaddtomacro\g@addto@macro
  \gdef\gplbacktext{}%
  \gdef\gplfronttext{}%
  \makeatother
  \ifGPblacktext
    \def\colorrgb#1{}%
    \def\colorgray#1{}%
  \else
    \ifGPcolor
      \def\colorrgb#1{\color[rgb]{#1}}%
      \def\colorgray#1{\color[gray]{#1}}%
      \expandafter\def\csname LTw\endcsname{\color{white}}%
      \expandafter\def\csname LTb\endcsname{\color{black}}%
      \expandafter\def\csname LTa\endcsname{\color{black}}%
      \expandafter\def\csname LT0\endcsname{\color[rgb]{1,0,0}}%
      \expandafter\def\csname LT1\endcsname{\color[rgb]{0,1,0}}%
      \expandafter\def\csname LT2\endcsname{\color[rgb]{0,0,1}}%
      \expandafter\def\csname LT3\endcsname{\color[rgb]{1,0,1}}%
      \expandafter\def\csname LT4\endcsname{\color[rgb]{0,1,1}}%
      \expandafter\def\csname LT5\endcsname{\color[rgb]{1,1,0}}%
      \expandafter\def\csname LT6\endcsname{\color[rgb]{0,0,0}}%
      \expandafter\def\csname LT7\endcsname{\color[rgb]{1,0.3,0}}%
      \expandafter\def\csname LT8\endcsname{\color[rgb]{0.5,0.5,0.5}}%
    \else
      \def\colorrgb#1{\color{black}}%
      \def\colorgray#1{\color[gray]{#1}}%
      \expandafter\def\csname LTw\endcsname{\color{white}}%
      \expandafter\def\csname LTb\endcsname{\color{black}}%
      \expandafter\def\csname LTa\endcsname{\color{black}}%
      \expandafter\def\csname LT0\endcsname{\color{black}}%
      \expandafter\def\csname LT1\endcsname{\color{black}}%
      \expandafter\def\csname LT2\endcsname{\color{black}}%
      \expandafter\def\csname LT3\endcsname{\color{black}}%
      \expandafter\def\csname LT4\endcsname{\color{black}}%
      \expandafter\def\csname LT5\endcsname{\color{black}}%
      \expandafter\def\csname LT6\endcsname{\color{black}}%
      \expandafter\def\csname LT7\endcsname{\color{black}}%
      \expandafter\def\csname LT8\endcsname{\color{black}}%
    \fi
  \fi
    \setlength{\unitlength}{0.0500bp}%
    \ifx\gptboxheight\undefined%
      \newlength{\gptboxheight}%
      \newlength{\gptboxwidth}%
      \newsavebox{\gptboxtext}%
    \fi%
    \setlength{\fboxrule}{0.5pt}%
    \setlength{\fboxsep}{1pt}%
\begin{picture}(9070.00,7086.00)%
    \gplgaddtomacro\gplbacktext{%
      \csname LTb\endcsname
      \put(682,719){\makebox(0,0)[r]{\strut{}\fsize $0$}}%
      \put(682,1485){\makebox(0,0)[r]{\strut{}\fsize $10$}}%
      \put(682,2252){\makebox(0,0)[r]{\strut{}\fsize $20$}}%
      \put(682,3018){\makebox(0,0)[r]{\strut{}\fsize $30$}}%
      \put(682,3785){\makebox(0,0)[r]{\strut{}\fsize $40$}}%
      \put(682,4551){\makebox(0,0)[r]{\strut{}\fsize $50$}}%
      \put(682,5317){\makebox(0,0)[r]{\strut{}\fsize $60$}}%
      \put(682,6084){\makebox(0,0)[r]{\strut{}\fsize $70$}}%
      \put(682,6850){\makebox(0,0)[r]{\strut{}\fsize $80$}}%
      \put(814,499){\makebox(0,0){\strut{}\fsize $0$}}%
      \put(2522,499){\makebox(0,0){\strut{}\fsize $0.0005$}}%
      \put(4231,499){\makebox(0,0){\strut{}\fsize $0.001$}}%
      \put(5939,499){\makebox(0,0){\strut{}\fsize $0.0015$}}%
      \put(7648,499){\makebox(0,0){\strut{}\fsize $0.002$}}%
    }%
    \gplgaddtomacro\gplfronttext{%
      \csname LTb\endcsname
      \put(198,3784){\rotatebox{-270}{\makebox(0,0){\strut{}\fsize $|F|$/(MN/m)}}}%
      \put(4743,169){\makebox(0,0){\strut{}\fsize $|u|$/m}}%
      \csname LTb\endcsname
      \put(3921,2267){\makebox(0,0)[l]{\strut{}\fsize drained elastic-plastic}}%
      \csname LTb\endcsname
      \put(3921,1937){\makebox(0,0)[l]{\strut{}\fsize undrained poro-elastic}}%
      \csname LTb\endcsname
      \put(3921,1607){\makebox(0,0)[l]{\strut{}\fsize undrained poro-elastic-plastic $K$}}%
      \csname LTb\endcsname
      \put(3921,1277){\makebox(0,0)[l]{\strut{}\fsize undrained poro-elastic-plastic $\tfrac{1}{5}K$}}%
      \csname LTb\endcsname
      \put(3921,947){\makebox(0,0)[l]{\strut{}\fsize undrained poro-elastic-plastic $5K$}}%
    }%
    \gplbacktext
    \put(0,0){\includegraphics{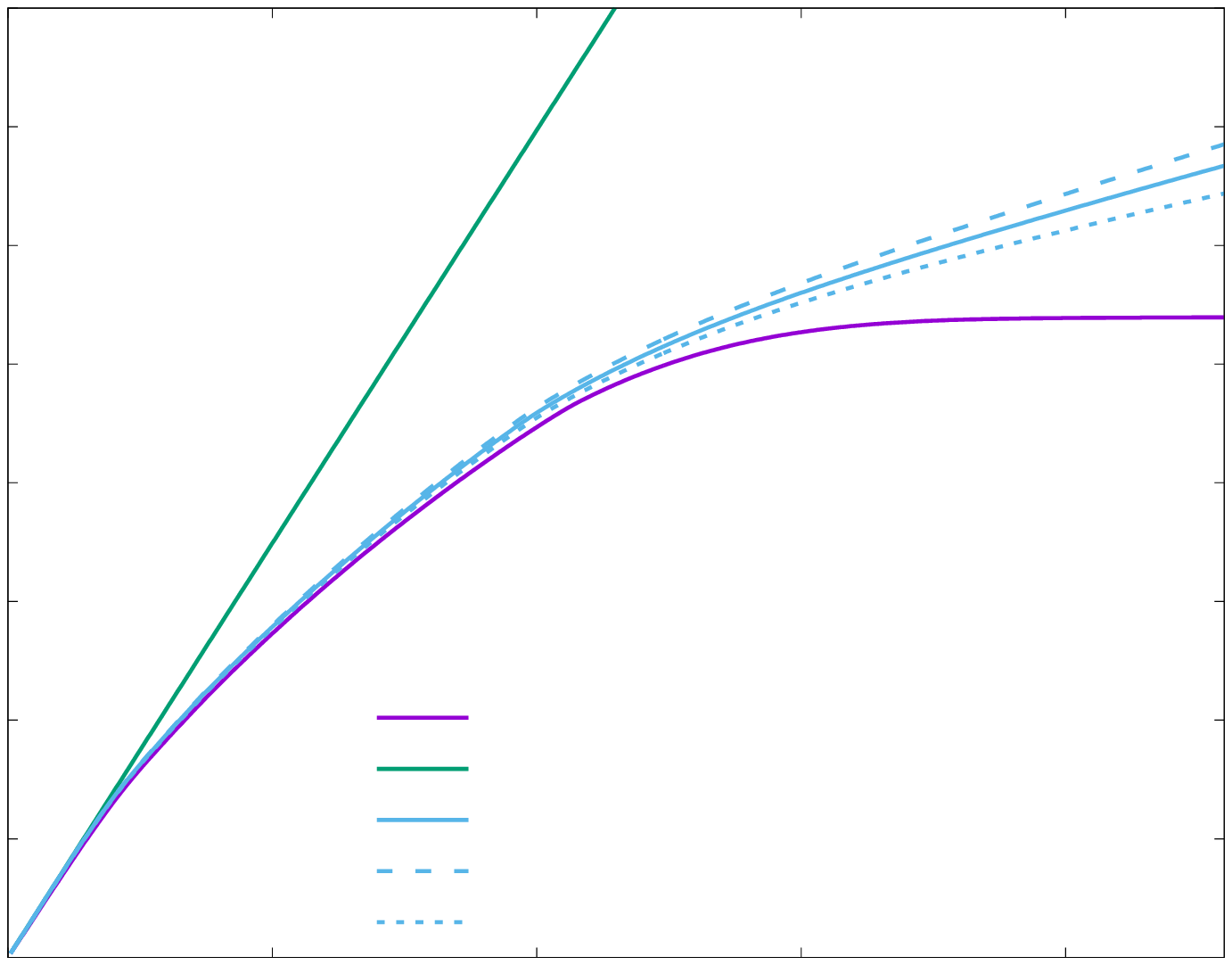}}%
    \gplfronttext
  \end{picture}%
  \endgroup
  }
  \captionsetup{width=0.9\textwidth}
  \caption{Rigid footing: Load displacement curve for the different tested material types. The load and displacement is taken from the area were the displacement is prescribed, see Figure~\ref{fig:bound1}.}
  \label{fig:load_displacement}
\end{figure}

In Figure~\ref{fig:flux_noflux} the distribution of the hardening variable $\alpha$ for the drained elastic-plastic material and undrained porous-elastic-plastic material is shown. It can be seen that the plastic deformation is more pronounced in the case of the drained material. This leads to the conclusion that the fluid within the material leads to an additional hardening effect. The load-displacement curve in Figure~\ref{fig:load_displacement} also shows this behaviour. Furthermore, a lower permeability leads to more pronounced hardening. This can be explained by the fact that the transport of the fluid is hindered and thus requires more work, see Figure~\ref{fig:load_displacement}.

\begin{figure}%
\centering
\footnotesize
\psfrag{a}    [c][c] {a)}
\psfrag{b}    [c][c] {b)}
\psfrag{pmax} [l][l] {$\phantom{-}1.2\cdot 10^{-1}$}
\psfrag{p}    [l][l] {$\phantom{-}m^e / ($kg$/$m$^3)$}
\psfrag{pmin} [l][l] {$-5.6\cdot 10^{-1}$}
\includegraphics*[width = 0.75\textwidth]{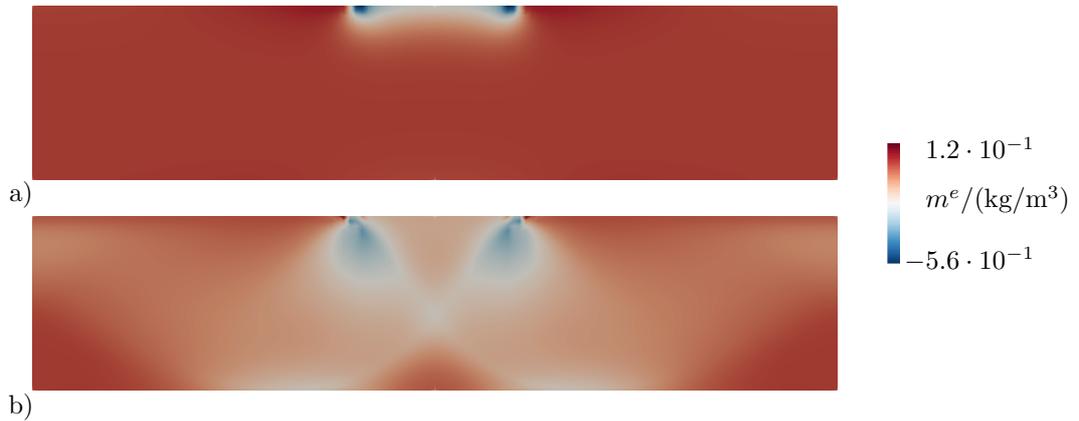}
\captionsetup{width=0.9\textwidth}
\caption{Rigid footing test. Change of the elastic fluid content $m^e$ for undrained porous-elastic material in a) and in b) for undrained porous-elastic-plastic material.
}
\label{fig:me_plast_elast}%
\end{figure}

Next, the distribution of the change of elastic fluid content $m^e$ in the domain for the undrained porous-elastic material is compared with the one of the undrained porous-elastic-plastic material, see Figure~\ref{fig:me_plast_elast}.

Note that for the elastic material the fluid is squeezed out right underneath the area where the loading is applied (see the negative change of the elastic fluid content at the boundary of the applied footing depicted in Figure~\ref{fig:me_plast_elast}a). Opposed to that, in case of the elastic-plastic material, the highest (negative) change of elastic fluid content is occurring in a more diffuse region that also extends to the bulk (see Figure~\ref{fig:me_plast_elast}b). This is precisely the area where most of the plastic deformation is happening, as can be observed in Figure~\ref{fig:flux_noflux} b). This phenomenon can be explained by the fact that the plastic deformation leads to a positive change of the plastic fluid content. Due to the additive decomposition of the change of the fluid content and the fact that no fluid is injected, the change of the elastic fluid content becomes negative in the plastifying areas.

\begin{figure}%
\centering
\footnotesize
\psfrag{a}     [c][c] {a)}
\psfrag{b}     [c][c] {b)}
\psfrag{c}     [c][c] {c)}
\psfrag{pmax1} [l][l] {$\phantom{-}1.2\cdot 10^{-1}$}
\psfrag{p1}    [l][l] {$\phantom{-}m^e / ($kg$/$m$^3)$}
\psfrag{pmin1} [l][l] {$-5.6\cdot 10^{-1}$}
\psfrag{pmax2} [l][l] {$\phantom{-}2.8$}
\psfrag{p2}    [l][l] {$\phantom{-}m^p / ($kg$/$m$^3)$}
\psfrag{pmin2} [l][l] {$\phantom{-}0$}
\psfrag{pmax3} [l][l] {$\phantom{-}2.5$}
\psfrag{p3}    [l][l] {$\phantom{-}m / ($kg$/$m$^3)$}
\psfrag{pmin3} [l][l] {$-1.4\cdot 10^{-1}$}
\includegraphics*[width = 0.75\textwidth]{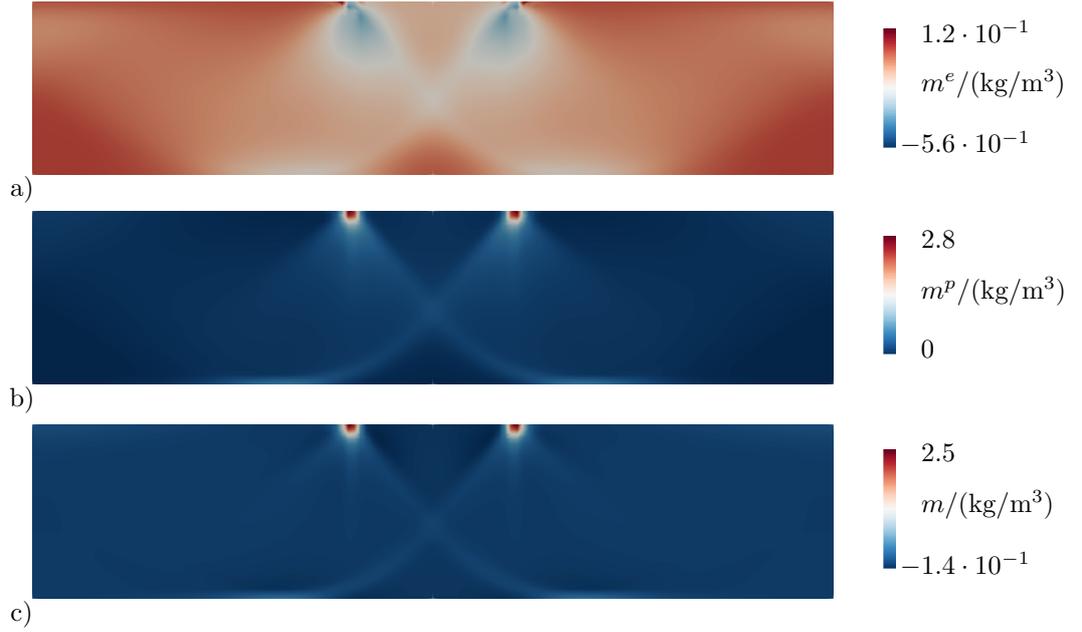}
\captionsetup{width=0.9\textwidth}
\caption{Rigid footing test: Change of fluid content for undrained porous-elastic-plastic material. Change of elastic fluid content $m^e$ in a), change of plastic fluid content $m^p$ in b) and change of total fluid content $m$ in c).
}
\label{fig:m_plast}%
\end{figure}

We depict the distribution of the change of the elastic $m^e$, the plastic $m^p$ and total fluid content $m$ for the undrained porous-elastic-plastic material in Figure~\ref{fig:m_plast}. One can observe that the change of the total fluid content $m$ is strongly dominated by the change of the plastic fluid content $m^p$.

\subsection{Comparison of different fracture driving forces for porous-elastic-plastic fracturing}

\begin{table}
\captionsetup{width=0.9\textwidth}
\caption{Material parameters for hydraulically induced fracture.}
\centering
\footnotesize
\renewcommand{\arraystretch}{1.3}
\begin{tabular}{|lr@{\ =\ }lc|lr@{\ =\ }l@{\hspace{1.5mm}}c|}
\hline
Lam\'e parameter      & $\lambda$ & 180.0 & GN/m$^2$         &  Biot's modulus & $M$  & 25.0 & GN/m$^2$ \\
Lam\'e parameter      & $G$     & 31.0 & GN/m$^2$          &  Biot's coefficient     & $b$   & 0.5 & -- \\     
hardening modulus     & $h$       & 5.0 & MN/m$^2$           &  fluid dyn. viscosity   & $\eta_f$     & $1.0\cdot 10^{-3}$    & Ns/m$^2$ \\
saturated yield shift & $\sigma_y$& 0.1  & MN/m$^2$               &  permeability      & $K$ &$9.8\cdot 10^{-12}$             & m$^3$s/kg \\                          
saturation            & $\omega$  & 2.0  & --                     & fluid density &$\rho_f$ & 1000.0 & kg/m$^3$ \\    
plastic viscosity     & $\eta_p$  & $5 \cdot 10^{-6}$ & s        &  crit. fracture energy & $\psi_c$   &  $5.0\cdot 10^{-8}$ & MN/m$^2$  \\
perturbation param.   & $q_1$     & $2 \cdot 10^{-5}$ & MN/m$^2$ & length scale &$l$                     & 0.5 & m \\  
slope yield function  & $M_\phi$ & 1.8 & --                       & residual stiffness &$k$                     & $1\cdot10^{-5}$ & -- \\
position of peak      & $s_\text{max}$  &   $2\cdot 10^{-3}$      & MN/m$^2$ & interpolation param. &$\epsilon$                     & 50 & -- \\
\hline
\end{tabular}
\label{tab:05crack}
\end{table}

\begin{figure}
\centering
{
\small
\psfrag{a}   [c][c]{$a$}%
\psfrag{aa}  [c][c]{$a/2$}%
\psfrag{m}   [c][c]{$\dot {\bar m}$}%
\psfrag{l1}  [c][c]{$L/2$}%
\psfrag{l2}  [c][c]{$L$}%
\psfrag{l11} [c][c]{$L$}%
\psfrag{x}   [l][l]{\footnotesize $x$}%
\psfrag{y}   [l][l]{\footnotesize $y$}%
\psfrag{h}   [c][r]{$\bar h$}%

\includegraphics[width=0.7\textwidth]{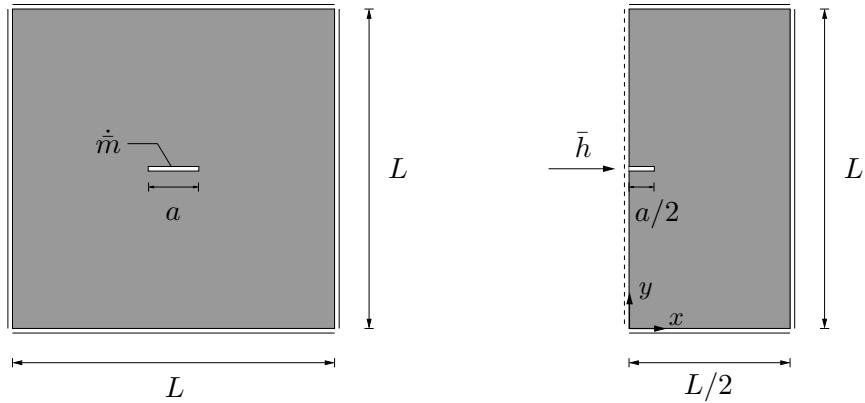}
}
\vspace{1mm}%
\captionsetup{width=0.9\textwidth}
\caption{Comparison of different fracture driving forces: Geometry and boundary conditions. Due to the symmetry of both only one half of the specimen is discretized. Here the prescribed flux is $\bar h=\dot{\bar m} = 0.01$~kg/s. All edges are mechanically fixed and permeable.
}
\label{fig:bound2}%
\end{figure}

In the present example, we investigate the influence of the presented fracture driving forces defined in \eqref{eq:fracevol2} and \eqref{eq:crackdrive}. For that purpose, a squared domain with the dimensions of $80$~m $\times$ $80$~m and a notch of the length $a=8$~m in its center is considered, see Figure~\ref{fig:bound2}. The fracture evolution is triggered by fluid injection into the notch.

Due to the symmetry of loading and geometry only one half of the domain is discretized with $12,060$ quadrilateral Raviart--Thomas-type enhanced-assumed-strain elements. The elements in the area surrounding the anticipated crack are refined yielding an element size of $h^e=0.25$~m in that region ($[0~$m$,40~$m$]$ $\times$ {$[31.875~$m$,48.125~$m$]$} around the notch).

The fluid injection is modeled by a prescribed fluid flux of $\bar h = 0.01~\tfrac{\text{kg}^3}{\text{s}}$. The time step is set to $\tau=1\cdot 10^{-3}$~s and the material parameters are listed in Table~\ref{tab:05crack}.

\begin{figure}%
\centering
\footnotesize
\psfrag{b}    [c][c] {a)}
\psfrag{c}    [c][c] {b)}
\psfrag{aa}   [c][c] {$\calH(\psi^0_\text{eff},\psi^0_\text{plast})$~~~~}
\psfrag{bb}   [c][c] {$\calH(\psi^0_\text{eff},w_\text{plast})$~~~~}
\psfrag{pmax} [l][l] {$1$}
\psfrag{p}    [l][l] {$d$}
\psfrag{pmin} [l][l] {$0$}
\includegraphics*[width = 0.75\textwidth]{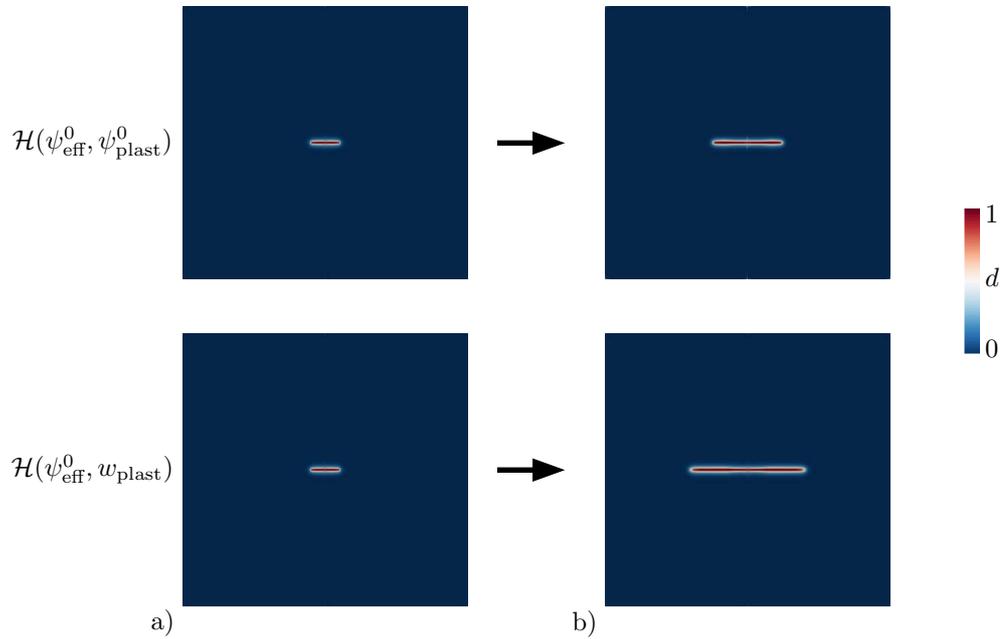}
\captionsetup{width=0.9\textwidth}
\caption{Comparison of different fracture driving forces: Distribution of the fracture phase-field for a porous-elastic-plastic with $\calH=f(\psi^0_\text{eff},\psi^0_\text{plast})$ and a porous-elastic-plastic material with $\calH=f(\psi^0_\text{eff},w_\text{plast})$ at two time steps: a) $t=0$~s and b) $t=90$~s.
}
\label{fig:d_plast}%
\end{figure}

The test was performed for the following two undrained settings with different choices of fracture driving forces:
\vspace{-0.3cm}
\begin{itemize}
  \setlength\itemsep{-0.2cm}
 \item[i)] porous-elastic-plastic material with $\calH = f(\psi^0_\text{eff},\psi^0_\text{plast})$ according to \eqref{eq:fracevol2}
 \item[ii)] porous-elastic-plastic material with $\calH = f(\psi^0_\text{eff},w_\text{plast})$ according to \eqref{eq:crackdrive}
 \end{itemize}
 \vspace{-0.3cm}
We now compare the hydraulically induced fracture lengths for the two porous-elastic-plastic settings. As can be observed in Figure~\ref{fig:d_plast}, both driving forces lead to the evolution of cracks. The fracture evolution in consideration of the fracture driving force $\calH = f(\psi^0_\text{eff},\psi^0_\text{plast})$ is however less prominent. Note that in that setting, only the elastic and hardening energies contribute to the fracture driving force. Thus, we would not obtain ductile fracture evolution in case of ideal plasticity with ($h=0$, $\sigma_y = 0$; not investigated here). In particular the latter observation justifies the presented modification of the fracture driving force in \eqref{eq:crackdrive}.

\subsection{Detailed analysis of hydraulically induced porous-elastic-plastic fracture} 

Finally, we investigate the ductile fracture evolution driven by an injected fluid in detail. The setup of the geometry and boundary conditions as well as the material parameters are taken from the previous example (please refer to Figure~\ref{fig:bound2} and Table~\ref{tab:05crack}).

The test was performed for two kinds of undrained materials:
\vspace{-0.3cm}
\begin{itemize}
  \setlength\itemsep{-0.2cm}
 \item[i)] porous-elastic material
 \item[ii)] porous-elastic-plastic material with $\calH = f(\psi^0_\text{eff},w_\text{plast})$.
\end{itemize}
\vspace{-0.3cm}

\begin{figure}%
\centering
\footnotesize
\psfrag{a}    [c][c] {a)}
\psfrag{c}    [c][c] {b)}
\psfrag{pmax} [l][l] {$1$}
\psfrag{p}    [l][l] {$d$}
\psfrag{pmin} [l][l] {$0$}
\includegraphics*[width = 0.7\textwidth]{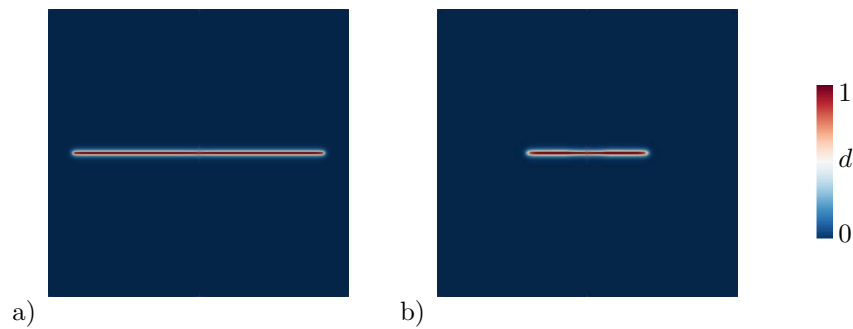}
\captionsetup{width=0.9\textwidth}
\caption{Hydraulically induced ductile fracture: Distribution of the fracture phase-field in a) porous-elastic material  in b) porous-elastic-plastic material at $t=90$~s.
}
\label{fig:d_elast_plast}%
\end{figure}%

\begin{figure}[h]
  \newcommand{\fsize}[0]{\Large}
  \centering
  \resizebox{!}{6cm}{%
\begingroup
  \makeatletter
  \providecommand\color[2][]{%
    \GenericError{(gnuplot) \space\space\space\@spaces}{%
      Package color not loaded in conjunction with
      terminal option `colourtext'%
    }{See the gnuplot documentation for explanation.%
    }{Either use 'blacktext' in gnuplot or load the package
      color.sty in LaTeX.}%
    \renewcommand\color[2][]{}%
  }%
  \providecommand\includegraphics[2][]{%
    \GenericError{(gnuplot) \space\space\space\@spaces}{%
      Package graphicx or graphics not loaded%
    }{See the gnuplot documentation for explanation.%
    }{The gnuplot epslatex terminal needs graphicx.sty or graphics.sty.}%
    \renewcommand\includegraphics[2][]{}%
  }%
  \providecommand\rotatebox[2]{#2}%
  \@ifundefined{ifGPcolor}{%
    \newif\ifGPcolor
    \GPcolortrue
  }{}%
  \@ifundefined{ifGPblacktext}{%
    \newif\ifGPblacktext
    \GPblacktexttrue
  }{}%
  \let\gplgaddtomacro\g@addto@macro
  \gdef\gplbacktext{}%
  \gdef\gplfronttext{}%
  \makeatother
  \ifGPblacktext
    \def\colorrgb#1{}%
    \def\colorgray#1{}%
  \else
    \ifGPcolor
      \def\colorrgb#1{\color[rgb]{#1}}%
      \def\colorgray#1{\color[gray]{#1}}%
      \expandafter\def\csname LTw\endcsname{\color{white}}%
      \expandafter\def\csname LTb\endcsname{\color{black}}%
      \expandafter\def\csname LTa\endcsname{\color{black}}%
      \expandafter\def\csname LT0\endcsname{\color[rgb]{1,0,0}}%
      \expandafter\def\csname LT1\endcsname{\color[rgb]{0,1,0}}%
      \expandafter\def\csname LT2\endcsname{\color[rgb]{0,0,1}}%
      \expandafter\def\csname LT3\endcsname{\color[rgb]{1,0,1}}%
      \expandafter\def\csname LT4\endcsname{\color[rgb]{0,1,1}}%
      \expandafter\def\csname LT5\endcsname{\color[rgb]{1,1,0}}%
      \expandafter\def\csname LT6\endcsname{\color[rgb]{0,0,0}}%
      \expandafter\def\csname LT7\endcsname{\color[rgb]{1,0.3,0}}%
      \expandafter\def\csname LT8\endcsname{\color[rgb]{0.5,0.5,0.5}}%
    \else
      \def\colorrgb#1{\color{black}}%
      \def\colorgray#1{\color[gray]{#1}}%
      \expandafter\def\csname LTw\endcsname{\color{white}}%
      \expandafter\def\csname LTb\endcsname{\color{black}}%
      \expandafter\def\csname LTa\endcsname{\color{black}}%
      \expandafter\def\csname LT0\endcsname{\color{black}}%
      \expandafter\def\csname LT1\endcsname{\color{black}}%
      \expandafter\def\csname LT2\endcsname{\color{black}}%
      \expandafter\def\csname LT3\endcsname{\color{black}}%
      \expandafter\def\csname LT4\endcsname{\color{black}}%
      \expandafter\def\csname LT5\endcsname{\color{black}}%
      \expandafter\def\csname LT6\endcsname{\color{black}}%
      \expandafter\def\csname LT7\endcsname{\color{black}}%
      \expandafter\def\csname LT8\endcsname{\color{black}}%
    \fi
  \fi
    \setlength{\unitlength}{0.0500bp}%
    \ifx\gptboxheight\undefined%
      \newlength{\gptboxheight}%
      \newlength{\gptboxwidth}%
      \newsavebox{\gptboxtext}%
    \fi%
    \setlength{\fboxrule}{0.5pt}%
    \setlength{\fboxsep}{1pt}%
\begin{picture}(9070.00,7086.00)%
    \gplgaddtomacro\gplbacktext{%
      \csname LTb\endcsname
      \put(946,822){\makebox(0,0)[r]{\strut{}\fsize \fsize \fsize \fsize \fsize $0$}}%
      \put(946,1668){\makebox(0,0)[r]{\strut{}\fsize \fsize \fsize \fsize \fsize $0.02$}}%
      \put(946,2515){\makebox(0,0)[r]{\strut{}\fsize \fsize \fsize \fsize \fsize $0.04$}}%
      \put(946,3361){\makebox(0,0)[r]{\strut{}\fsize \fsize \fsize \fsize \fsize $0.06$}}%
      \put(946,4208){\makebox(0,0)[r]{\strut{}\fsize \fsize \fsize \fsize \fsize $0.08$}}%
      \put(946,5054){\makebox(0,0)[r]{\strut{}\fsize \fsize \fsize \fsize \fsize $0.1$}}%
      \put(946,5901){\makebox(0,0)[r]{\strut{}\fsize \fsize \fsize \fsize \fsize $0.12$}}%
      \put(946,6747){\makebox(0,0)[r]{\strut{}\fsize \fsize \fsize \fsize \fsize $0.14$}}%
      \put(1078,602){\makebox(0,0){\strut{}\fsize \fsize \fsize \fsize \fsize $0$}}%
      \put(2344,602){\makebox(0,0){\strut{}\fsize \fsize \fsize \fsize \fsize $0.00015$}}%
      \put(3610,602){\makebox(0,0){\strut{}\fsize \fsize \fsize \fsize \fsize $0.0003$}}%
      \put(4876,602){\makebox(0,0){\strut{}\fsize \fsize \fsize \fsize \fsize $0.00045$}}%
      \put(6141,602){\makebox(0,0){\strut{}\fsize \fsize \fsize \fsize \fsize $0.0006$}}%
      \put(7407,602){\makebox(0,0){\strut{}\fsize \fsize \fsize \fsize \fsize $0.00075$}}%
      \put(8673,602){\makebox(0,0){\strut{}\fsize \fsize \fsize \fsize \fsize $0.0009$}}%
    }%
    \gplgaddtomacro\gplfronttext{%
      \csname LTb\endcsname
      \put(0,3784){\rotatebox{-270}{\makebox(0,0){\strut{}\fsize $p$/Mpa}}}%
      \put(4875,0){\makebox(0,0){\strut{}\fsize  $V$/m$^3$}}%
      \csname LTb\endcsname
      \put(5901,1380){\makebox(0,0)[l]{\strut{}\fsize poro-elastic}}%
      \csname LTb\endcsname
      \put(5901,1050){\makebox(0,0)[l]{\strut{}\fsize poro-elastic-plastic}}%
    }%
    \gplbacktext
    \put(0,0){\includegraphics{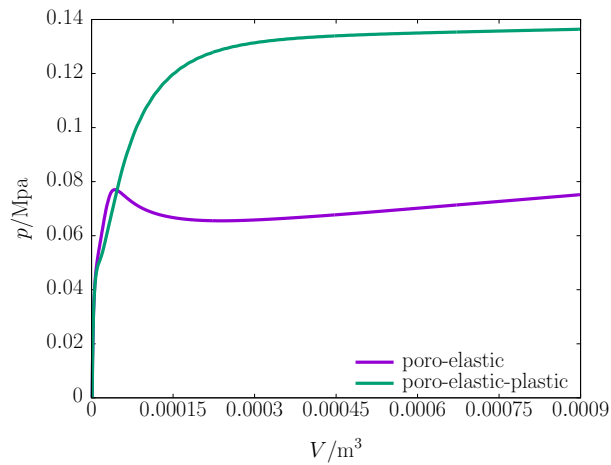}}%
    \gplfronttext
  \end{picture}%
  \endgroup
  }
  \captionsetup{width=0.9\textwidth}
  \caption{Hydraulically induced ductile fracture: Fluid pressure within the fracture at $x=1.0$~m, $y=39.875$~m over inject volume for porous-elastic and porous-elastic-plastic material.}
  \label{fig:fluidpressure}
\end{figure}

As can be seen in Figure~\ref{fig:d_elast_plast}, the length of the finally induced crack for the porous-elastic material is much more pronounced that in case of the porous-elastic-plastic material. This goes along with the observation of a higher fluid pressure inside the crack, see Figure~\ref{fig:fluidpressure}. We conclude that in case of an elastic-plastic material more fluid needs to be injected into the crack to drive fracturing. For the elastic material we can observe a characteristic drop of the pressure within the fracture at the onset of fracture propagation (injected fluid volume $V\approx 0.00004$~m$^3$). In the elastic-plastic material this drop cannot be observed, see again Figure~\ref{fig:fluidpressure}.

\begin{figure}%
\centering
\footnotesize
\psfrag{a}    [c][c] {a)}
\psfrag{c}    [c][c] {b)}
\psfrag{pmax} [l][l] {$\phantom{-}4.2\cdot 10^{-2}$}
\psfrag{p}    [l][l] {$\phantom{-}m^e/($kg/m$^3)$}
\psfrag{pmin} [l][l] {$-3.2\cdot 10^{-4}$}
\includegraphics*[width = 0.7\textwidth]{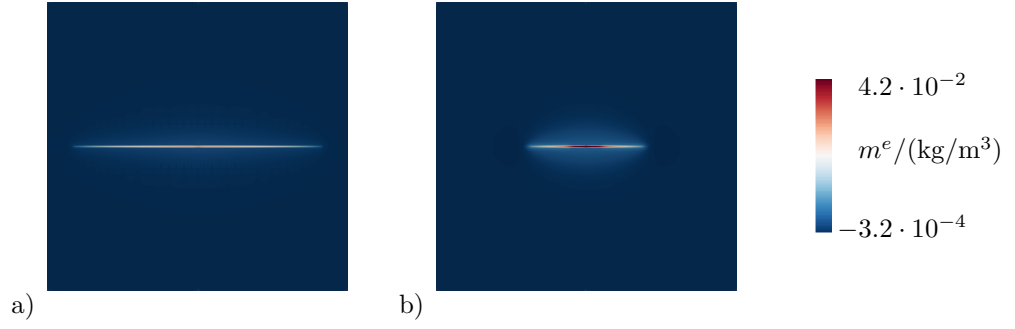}
\captionsetup{width=0.9\textwidth}
\caption{Hydraulically induced ductile fracture: Change of elastic
  fluid content $m^e$ for porous-elastic material in a) and porous-elastic-plastic material in b).
}
\label{fig:m_all}%
\end{figure}

\begin{figure}%
\centering
%
\scriptsize
\psfrag{a}     [c][c] {a)}
\psfrag{b}     [c][c] {b)}
\psfrag{c}     [c][c] {c)}
\psfrag{pmax1} [r][r] {$\phantom{-}4.2\cdot 10^{-2}$}
\psfrag{p1}    [r][r] {$\phantom{-}m,m^e/($kg/m$^3)$}
\psfrag{pmin1} [r][r] {$-3.2\cdot 10^{-4}$}
\psfrag{pmax}  [l][l] {$1.3\cdot 10^{-3}$}
\psfrag{p}     [l][l] {$m^p/($kg/m$^3)$}
\psfrag{pmin}  [l][l] {$0.0$}
\includegraphics*[width = 1.0\textwidth]{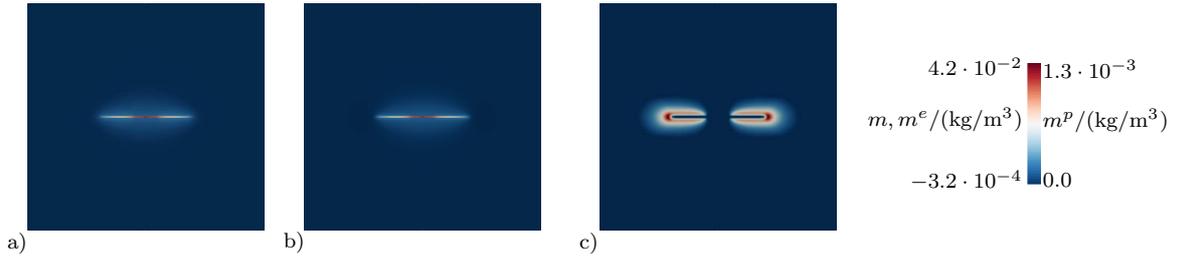}
\captionsetup{width=0.9\textwidth}
\caption{Hydraulically induced ductile fracture: Change of total fluid
  content $m$ in a), change of elastic fluid content $m^e$ in b) and change of plastic fluid content $m^p$ in c) for porous-elastic-plastic material.
}
\label{fig:mp}%
\end{figure}

In Figure~\ref{fig:m_all} the distribution of the change of the elastic fluid content is shown for the final equilibrium state. The individual lengths of the cracks are clearly visible. Due to the short crack length in the elastic-plastic case, the injected fluid is distributed over a smaller region. This then gives rise to a higher change of the elastic fluid content, in particular close to the fracture center. Note carefully that the change of the elastic fluid content in front of the fracture tips is negative for the elastic-plastic material. This phenomenon is investigated in a more detailed way in Figure~\ref{fig:mp}, where the contributions of the change of the fluid content in the elastic-plastic material are shown.

By taking a look at Figure~\ref{fig:mp}c, it can be seen that most of the volumetric plastic deformation occurs at the fracture tips ($\dot m^p = \rho_f b \tr \dot \Bve^p$). This leads to a positive change of the plastic fluid content. Since the permeability in the bulk is comparably low, very little amount of fluid diffuses from the fracture into the bulk. In other words, the change of the fluid content at the fracture tips is almost zero. Due to the definition $m=m^e+m^p$ the positive change of the plastic fluid content leads to a negative change of elastic fluid content. Hence the fluid in the fully saturated medium, which is initially stored elastically, is now stored plastically due to the plastic deformation of the solid matrix.

\begin{figure}
  \newcommand{\fsize}[0]{\Large}
  \begin{minipage}{0.34\textwidth}
  \resizebox{!}{4cm}{
\begingroup
  \makeatletter
  \providecommand\color[2][]{%
    \GenericError{(gnuplot) \space\space\space\@spaces}{%
      Package color not loaded in conjunction with
      terminal option `colourtext'%
    }{See the gnuplot documentation for explanation.%
    }{Either use 'blacktext' in gnuplot or load the package
      color.sty in LaTeX.}%
    \renewcommand\color[2][]{}%
  }%
  \providecommand\includegraphics[2][]{%
    \GenericError{(gnuplot) \space\space\space\@spaces}{%
      Package graphicx or graphics not loaded%
    }{See the gnuplot documentation for explanation.%
    }{The gnuplot epslatex terminal needs graphicx.sty or graphics.sty.}%
    \renewcommand\includegraphics[2][]{}%
  }%
  \providecommand\rotatebox[2]{#2}%
  \@ifundefined{ifGPcolor}{%
    \newif\ifGPcolor
    \GPcolortrue
  }{}%
  \@ifundefined{ifGPblacktext}{%
    \newif\ifGPblacktext
    \GPblacktexttrue
  }{}%
  \let\gplgaddtomacro\g@addto@macro
  \gdef\gplbacktext{}%
  \gdef\gplfronttext{}%
  \makeatother
  \ifGPblacktext
    \def\colorrgb#1{}%
    \def\colorgray#1{}%
  \else
    \ifGPcolor
      \def\colorrgb#1{\color[rgb]{#1}}%
      \def\colorgray#1{\color[gray]{#1}}%
      \expandafter\def\csname LTw\endcsname{\color{white}}%
      \expandafter\def\csname LTb\endcsname{\color{black}}%
      \expandafter\def\csname LTa\endcsname{\color{black}}%
      \expandafter\def\csname LT0\endcsname{\color[rgb]{1,0,0}}%
      \expandafter\def\csname LT1\endcsname{\color[rgb]{0,1,0}}%
      \expandafter\def\csname LT2\endcsname{\color[rgb]{0,0,1}}%
      \expandafter\def\csname LT3\endcsname{\color[rgb]{1,0,1}}%
      \expandafter\def\csname LT4\endcsname{\color[rgb]{0,1,1}}%
      \expandafter\def\csname LT5\endcsname{\color[rgb]{1,1,0}}%
      \expandafter\def\csname LT6\endcsname{\color[rgb]{0,0,0}}%
      \expandafter\def\csname LT7\endcsname{\color[rgb]{1,0.3,0}}%
      \expandafter\def\csname LT8\endcsname{\color[rgb]{0.5,0.5,0.5}}%
    \else
      \def\colorrgb#1{\color{black}}%
      \def\colorgray#1{\color[gray]{#1}}%
      \expandafter\def\csname LTw\endcsname{\color{white}}%
      \expandafter\def\csname LTb\endcsname{\color{black}}%
      \expandafter\def\csname LTa\endcsname{\color{black}}%
      \expandafter\def\csname LT0\endcsname{\color{black}}%
      \expandafter\def\csname LT1\endcsname{\color{black}}%
      \expandafter\def\csname LT2\endcsname{\color{black}}%
      \expandafter\def\csname LT3\endcsname{\color{black}}%
      \expandafter\def\csname LT4\endcsname{\color{black}}%
      \expandafter\def\csname LT5\endcsname{\color{black}}%
      \expandafter\def\csname LT6\endcsname{\color{black}}%
      \expandafter\def\csname LT7\endcsname{\color{black}}%
      \expandafter\def\csname LT8\endcsname{\color{black}}%
    \fi
  \fi
    \setlength{\unitlength}{0.0500bp}%
    \ifx\gptboxheight\undefined%
      \newlength{\gptboxheight}%
      \newlength{\gptboxwidth}%
      \newsavebox{\gptboxtext}%
    \fi%
    \setlength{\fboxrule}{0.5pt}%
    \setlength{\fboxsep}{1pt}%
\begin{picture}(9070.00,7086.00)%
    \gplgaddtomacro\gplbacktext{%
      \csname LTb\endcsname
      \put(946,1050){\makebox(0,0)[r]{\strut{}\fsize $0$}}%
      \put(946,2189){\makebox(0,0)[r]{\strut{}\fsize $0.01$}}%
      \put(946,3329){\makebox(0,0)[r]{\strut{}\fsize $0.02$}}%
      \put(946,4468){\makebox(0,0)[r]{\strut{}\fsize $0.03$}}%
      \put(946,5608){\makebox(0,0)[r]{\strut{}\fsize $0.04$}}%
      \put(946,6747){\makebox(0,0)[r]{\strut{}\fsize $0.05$}}%
      \put(1078,602){\makebox(0,0){\strut{}\fsize $0$}}%
      \put(2027,602){\makebox(0,0){\strut{}\fsize $10$}}%
      \put(2977,602){\makebox(0,0){\strut{}\fsize $20$}}%
      \put(3926,602){\makebox(0,0){\strut{}\fsize $30$}}%
      \put(4876,602){\makebox(0,0){\strut{}\fsize $40$}}%
      \put(5825,602){\makebox(0,0){\strut{}\fsize $50$}}%
      \put(6774,602){\makebox(0,0){\strut{}\fsize $60$}}%
      \put(7724,602){\makebox(0,0){\strut{}\fsize $70$}}%
      \put(8673,602){\makebox(0,0){\strut{}\fsize $80$}}%
    }%
    \gplgaddtomacro\gplfronttext{%
      \csname LTb\endcsname
      \put(0,3784){\rotatebox{-270}{\makebox(0,0){\strut{}\fsize $m$/(kg/m$^3$)}}}%
      \put(4875,272){\makebox(0,0){\strut{}\fsize $x$/m}}%
      \csname LTb\endcsname
      \put(2065,6519){\makebox(0,0)[l]{\strut{}\fsize porous-elastic-plastic $m$}}%
      \csname LTb\endcsname
      \put(2065,6189){\makebox(0,0)[l]{\strut{}\fsize porous-elastic-plastic $m^p$}}%
      \csname LTb\endcsname
      \put(2065,5859){\makebox(0,0)[l]{\strut{}\fsize porous-elastic $m$}}%
    }%
    \gplbacktext
    \put(0,0){\includegraphics{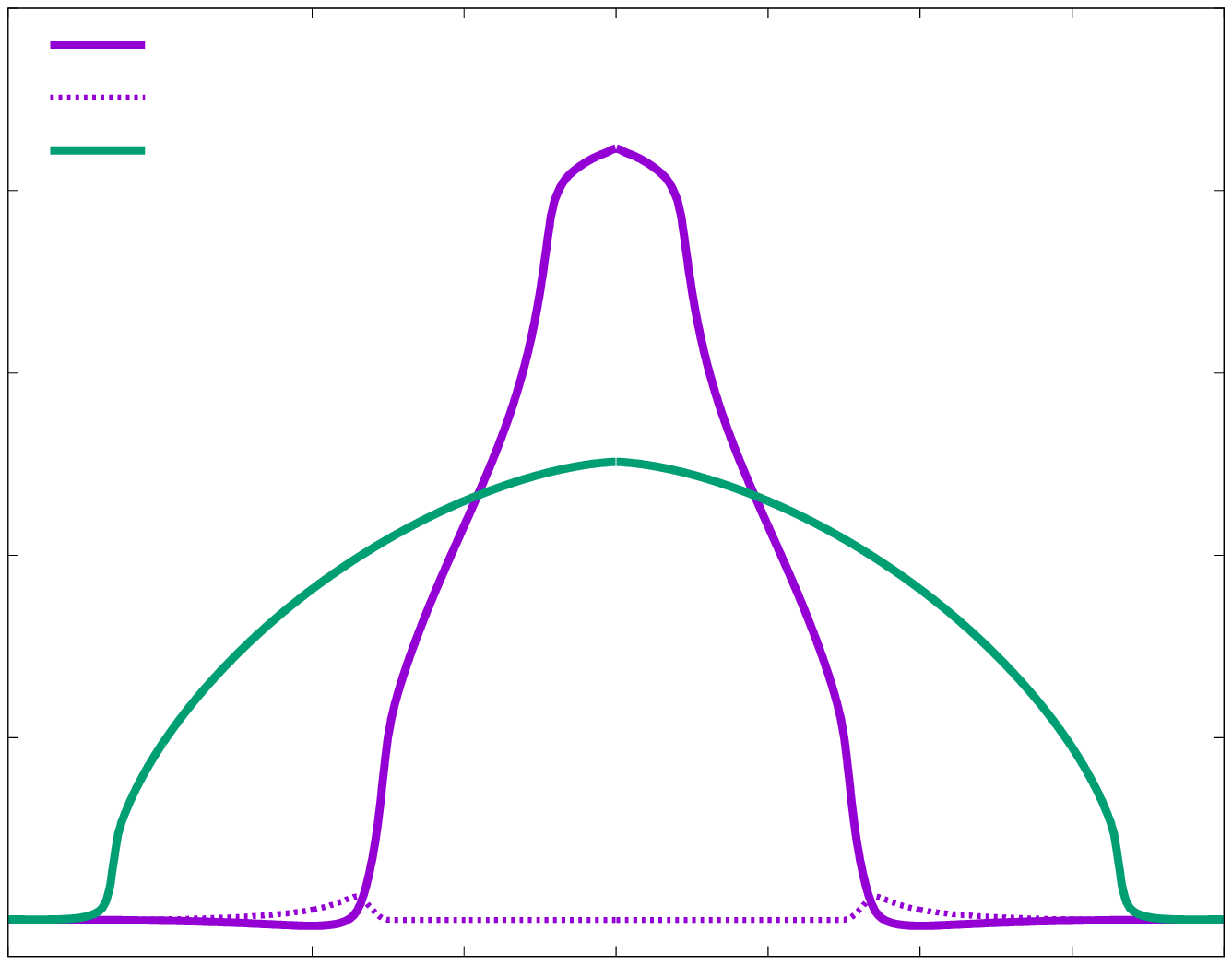}}%
    \gplfronttext
  \end{picture}%
\endgroup
}
  \end{minipage} \hspace{-5mm}
  \begin{minipage}{0.34\textwidth}
  \resizebox{!}{4cm}{
\begingroup
  \makeatletter
  \providecommand\color[2][]{%
    \GenericError{(gnuplot) \space\space\space\@spaces}{%
      Package color not loaded in conjunction with
      terminal option `colourtext'%
    }{See the gnuplot documentation for explanation.%
    }{Either use 'blacktext' in gnuplot or load the package
      color.sty in LaTeX.}%
    \renewcommand\color[2][]{}%
  }%
  \providecommand\includegraphics[2][]{%
    \GenericError{(gnuplot) \space\space\space\@spaces}{%
      Package graphicx or graphics not loaded%
    }{See the gnuplot documentation for explanation.%
    }{The gnuplot epslatex terminal needs graphicx.sty or graphics.sty.}%
    \renewcommand\includegraphics[2][]{}%
  }%
  \providecommand\rotatebox[2]{#2}%
  \@ifundefined{ifGPcolor}{%
    \newif\ifGPcolor
    \GPcolortrue
  }{}%
  \@ifundefined{ifGPblacktext}{%
    \newif\ifGPblacktext
    \GPblacktexttrue
  }{}%
  \let\gplgaddtomacro\g@addto@macro
  \gdef\gplbacktext{}%
  \gdef\gplfronttext{}%
  \makeatother
  \ifGPblacktext
    \def\colorrgb#1{}%
    \def\colorgray#1{}%
  \else
    \ifGPcolor
      \def\colorrgb#1{\color[rgb]{#1}}%
      \def\colorgray#1{\color[gray]{#1}}%
      \expandafter\def\csname LTw\endcsname{\color{white}}%
      \expandafter\def\csname LTb\endcsname{\color{black}}%
      \expandafter\def\csname LTa\endcsname{\color{black}}%
      \expandafter\def\csname LT0\endcsname{\color[rgb]{1,0,0}}%
      \expandafter\def\csname LT1\endcsname{\color[rgb]{0,1,0}}%
      \expandafter\def\csname LT2\endcsname{\color[rgb]{0,0,1}}%
      \expandafter\def\csname LT3\endcsname{\color[rgb]{1,0,1}}%
      \expandafter\def\csname LT4\endcsname{\color[rgb]{0,1,1}}%
      \expandafter\def\csname LT5\endcsname{\color[rgb]{1,1,0}}%
      \expandafter\def\csname LT6\endcsname{\color[rgb]{0,0,0}}%
      \expandafter\def\csname LT7\endcsname{\color[rgb]{1,0.3,0}}%
      \expandafter\def\csname LT8\endcsname{\color[rgb]{0.5,0.5,0.5}}%
    \else
      \def\colorrgb#1{\color{black}}%
      \def\colorgray#1{\color[gray]{#1}}%
      \expandafter\def\csname LTw\endcsname{\color{white}}%
      \expandafter\def\csname LTb\endcsname{\color{black}}%
      \expandafter\def\csname LTa\endcsname{\color{black}}%
      \expandafter\def\csname LT0\endcsname{\color{black}}%
      \expandafter\def\csname LT1\endcsname{\color{black}}%
      \expandafter\def\csname LT2\endcsname{\color{black}}%
      \expandafter\def\csname LT3\endcsname{\color{black}}%
      \expandafter\def\csname LT4\endcsname{\color{black}}%
      \expandafter\def\csname LT5\endcsname{\color{black}}%
      \expandafter\def\csname LT6\endcsname{\color{black}}%
      \expandafter\def\csname LT7\endcsname{\color{black}}%
      \expandafter\def\csname LT8\endcsname{\color{black}}%
    \fi
  \fi
    \setlength{\unitlength}{0.0500bp}%
    \ifx\gptboxheight\undefined%
      \newlength{\gptboxheight}%
      \newlength{\gptboxwidth}%
      \newsavebox{\gptboxtext}%
    \fi%
    \setlength{\fboxrule}{0.5pt}%
    \setlength{\fboxsep}{1pt}%
\begin{picture}(9070.00,7086.00)%
    \gplgaddtomacro\gplbacktext{%
      \csname LTb\endcsname
      \put(1078,874){\makebox(0,0)[r]{\strut{}\fsize $-0.02$}}%
      \put(1078,1602){\makebox(0,0)[r]{\strut{}\fsize $0$}}%
      \put(1078,2329){\makebox(0,0)[r]{\strut{}\fsize $0.02$}}%
      \put(1078,3057){\makebox(0,0)[r]{\strut{}\fsize $0.04$}}%
      \put(1078,3785){\makebox(0,0)[r]{\strut{}\fsize $0.06$}}%
      \put(1078,4512){\makebox(0,0)[r]{\strut{}\fsize $0.08$}}%
      \put(1078,5240){\makebox(0,0)[r]{\strut{}\fsize $0.1$}}%
      \put(1078,5967){\makebox(0,0)[r]{\strut{}\fsize $0.12$}}%
      \put(1078,6695){\makebox(0,0)[r]{\strut{}\fsize $0.14$}}%
      \put(1210,654){\makebox(0,0){\strut{}\fsize $0$}}%
      \put(2143,654){\makebox(0,0){\strut{}\fsize $10$}}%
      \put(3076,654){\makebox(0,0){\strut{}\fsize $20$}}%
      \put(4009,654){\makebox(0,0){\strut{}\fsize $30$}}%
      \put(4942,654){\makebox(0,0){\strut{}\fsize $40$}}%
      \put(5874,654){\makebox(0,0){\strut{}\fsize $50$}}%
      \put(6807,654){\makebox(0,0){\strut{}\fsize $60$}}%
      \put(7740,654){\makebox(0,0){\strut{}\fsize $70$}}%
      \put(8673,654){\makebox(0,0){\strut{}\fsize $80$}}%
    }%
    \gplgaddtomacro\gplfronttext{%
      \csname LTb\endcsname
      \put(198,3784){\rotatebox{-270}{\makebox(0,0){\strut{}\fsize $p$/MPa}}}%
      \put(4941,324){\makebox(0,0){\strut{}\fsize $x$/m}}%
      \csname LTb\endcsname
      \put(2197,6467){\makebox(0,0)[l]{\strut{}\fsize porous-elastic-plastic}}%
      \csname LTb\endcsname
      \put(2197,6137){\makebox(0,0)[l]{\strut{}\fsize porous-elastic}}%
    }%
    \gplbacktext
    \put(0,0){\includegraphics{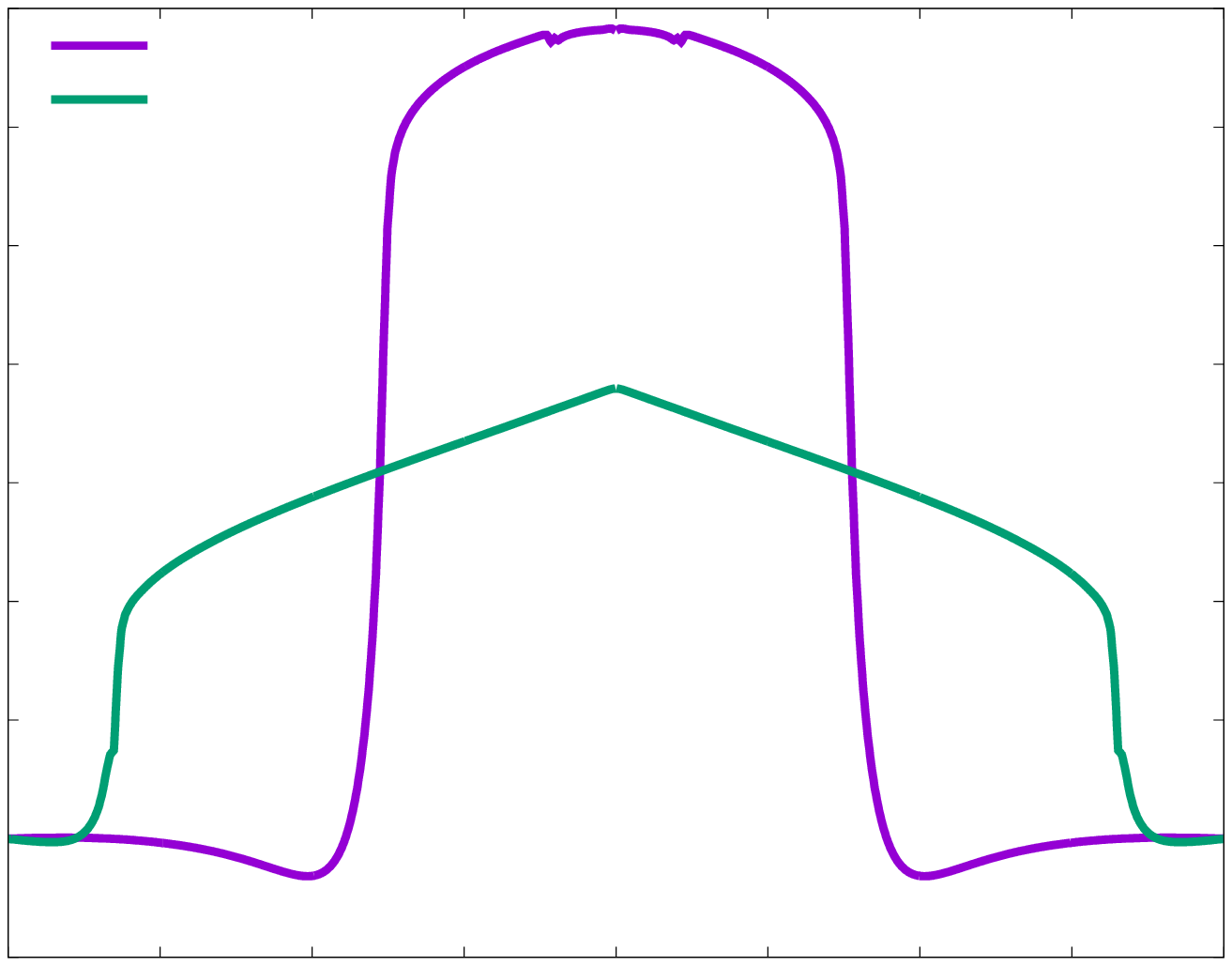}}%
    \gplfronttext
  \end{picture}%
  \endgroup
  }
  \end{minipage} \hspace{-5mm}  
  \begin{minipage}{0.34\textwidth}
  \resizebox{!}{3.9cm}{
\begingroup
  \makeatletter
  \providecommand\color[2][]{%
    \GenericError{(gnuplot) \space\space\space\@spaces}{%
      Package color not loaded in conjunction with
      terminal option `colourtext'%
    }{See the gnuplot documentation for explanation.%
    }{Either use 'blacktext' in gnuplot or load the package
      color.sty in LaTeX.}%
    \renewcommand\color[2][]{}%
  }%
  \providecommand\includegraphics[2][]{%
    \GenericError{(gnuplot) \space\space\space\@spaces}{%
      Package graphicx or graphics not loaded%
    }{See the gnuplot documentation for explanation.%
    }{The gnuplot epslatex terminal needs graphicx.sty or graphics.sty.}%
    \renewcommand\includegraphics[2][]{}%
  }%
  \providecommand\rotatebox[2]{#2}%
  \@ifundefined{ifGPcolor}{%
    \newif\ifGPcolor
    \GPcolortrue
  }{}%
  \@ifundefined{ifGPblacktext}{%
    \newif\ifGPblacktext
    \GPblacktexttrue
  }{}%
  \let\gplgaddtomacro\g@addto@macro
  \gdef\gplbacktext{}%
  \gdef\gplfronttext{}%
  \makeatother
  \ifGPblacktext
    \def\colorrgb#1{}%
    \def\colorgray#1{}%
  \else
    \ifGPcolor
      \def\colorrgb#1{\color[rgb]{#1}}%
      \def\colorgray#1{\color[gray]{#1}}%
      \expandafter\def\csname LTw\endcsname{\color{white}}%
      \expandafter\def\csname LTb\endcsname{\color{black}}%
      \expandafter\def\csname LTa\endcsname{\color{black}}%
      \expandafter\def\csname LT0\endcsname{\color[rgb]{1,0,0}}%
      \expandafter\def\csname LT1\endcsname{\color[rgb]{0,1,0}}%
      \expandafter\def\csname LT2\endcsname{\color[rgb]{0,0,1}}%
      \expandafter\def\csname LT3\endcsname{\color[rgb]{1,0,1}}%
      \expandafter\def\csname LT4\endcsname{\color[rgb]{0,1,1}}%
      \expandafter\def\csname LT5\endcsname{\color[rgb]{1,1,0}}%
      \expandafter\def\csname LT6\endcsname{\color[rgb]{0,0,0}}%
      \expandafter\def\csname LT7\endcsname{\color[rgb]{1,0.3,0}}%
      \expandafter\def\csname LT8\endcsname{\color[rgb]{0.5,0.5,0.5}}%
    \else
      \def\colorrgb#1{\color{black}}%
      \def\colorgray#1{\color[gray]{#1}}%
      \expandafter\def\csname LTw\endcsname{\color{white}}%
      \expandafter\def\csname LTb\endcsname{\color{black}}%
      \expandafter\def\csname LTa\endcsname{\color{black}}%
      \expandafter\def\csname LT0\endcsname{\color{black}}%
      \expandafter\def\csname LT1\endcsname{\color{black}}%
      \expandafter\def\csname LT2\endcsname{\color{black}}%
      \expandafter\def\csname LT3\endcsname{\color{black}}%
      \expandafter\def\csname LT4\endcsname{\color{black}}%
      \expandafter\def\csname LT5\endcsname{\color{black}}%
      \expandafter\def\csname LT6\endcsname{\color{black}}%
      \expandafter\def\csname LT7\endcsname{\color{black}}%
      \expandafter\def\csname LT8\endcsname{\color{black}}%
    \fi
  \fi
    \setlength{\unitlength}{0.0500bp}%
    \ifx\gptboxheight\undefined%
      \newlength{\gptboxheight}%
      \newlength{\gptboxwidth}%
      \newsavebox{\gptboxtext}%
    \fi%
    \setlength{\fboxrule}{0.5pt}%
    \setlength{\fboxsep}{1pt}%
\begin{picture}(9070.00,7086.00)%
    \gplgaddtomacro\gplbacktext{%
      \csname LTb\endcsname
      \put(682,734){\makebox(0,0)[r]{\strut{}\fsize 0}}%
      \put(682,1346){\makebox(0,0)[r]{\strut{}\fsize 2}}%
      \put(682,1957){\makebox(0,0)[r]{\strut{}\fsize 4}}%
      \put(682,2569){\makebox(0,0)[r]{\strut{}\fsize 6}}%
      \put(682,3181){\makebox(0,0)[r]{\strut{}\fsize 8}}%
      \put(682,3792){\makebox(0,0)[r]{\strut{}\fsize 10}}%
      \put(682,4404){\makebox(0,0)[r]{\strut{}\fsize 12}}%
      \put(682,5015){\makebox(0,0)[r]{\strut{}\fsize 14}}%
      \put(682,5627){\makebox(0,0)[r]{\strut{}\fsize 16}}%
      \put(682,6238){\makebox(0,0)[r]{\strut{}\fsize 18}}%
      \put(682,6850){\makebox(0,0)[r]{\strut{}\fsize 20}}%
      \put(814,499){\makebox(0,0){\strut{}\fsize $0$}}%
      \put(1796,499){\makebox(0,0){\strut{}\fsize $10$}}%
      \put(2779,499){\makebox(0,0){\strut{}\fsize $20$}}%
      \put(3761,499){\makebox(0,0){\strut{}\fsize $30$}}%
      \put(4744,499){\makebox(0,0){\strut{}\fsize $40$}}%
      \put(5726,499){\makebox(0,0){\strut{}\fsize $50$}}%
      \put(6708,499){\makebox(0,0){\strut{}\fsize $60$}}%
      \put(7691,499){\makebox(0,0){\strut{}\fsize $70$}}%
      \put(8673,499){\makebox(0,0){\strut{}\fsize $80$}}%
    }%
    \gplgaddtomacro\gplfronttext{%
      \csname LTb\endcsname
      \put(198,3784){\rotatebox{-270}{\makebox(0,0){\strut{}\fsize $w$/($1 \cdot 10^{-6}$ m)}}}%
      \put(4743,169){\makebox(0,0){\strut{}\fsize $x$/m}}%
      \csname LTb\endcsname
      \put(1801,6622){\makebox(0,0)[l]{\strut{}\fsize poro-elastic-plastic}}%
      \csname LTb\endcsname
      \put(1801,6292){\makebox(0,0)[l]{\strut{}\fsize poro-elastic}}%
    }%
    \gplbacktext
    \put(0,0){\includegraphics{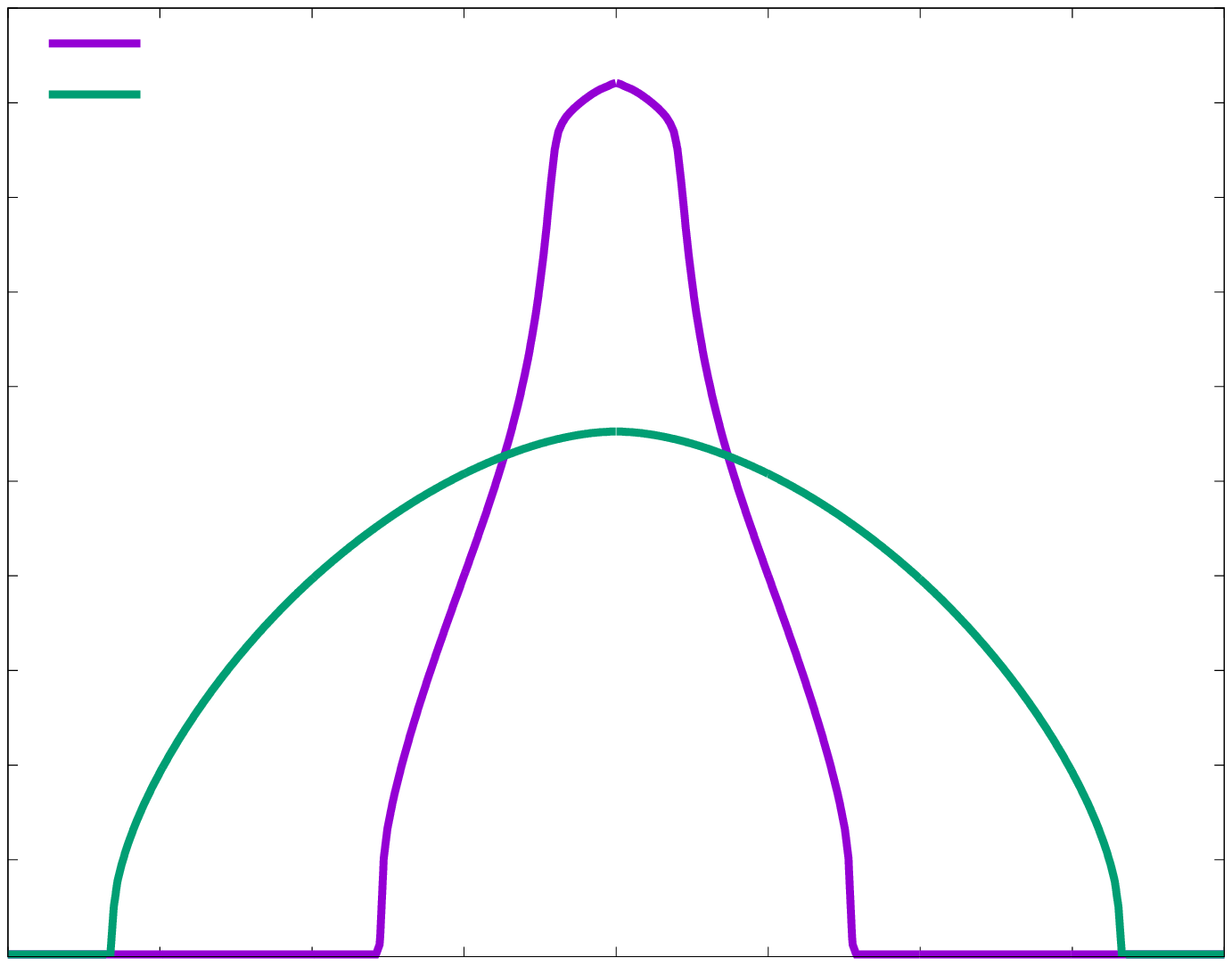}}%
    \gplfronttext
  \end{picture}%
\endgroup
}
  \end{minipage}
  \captionsetup{width=0.9\textwidth}
  \caption{Hydraulically induced ductile fracture: Change of fluid content $m$ and $m^p$, fluid pressure $p$ and fracture opening $w$ for porous-elastic and porous-elastic-plastic material over $x$ at $y = 40$~m.}
  \label{fig:elast_plastic_x}
\end{figure}

In Figure~\ref{fig:elast_plastic_x} we show the distribution of the change of fluid content as well as the pressure and the fracture-opening width along $x$ at $y=40$~m. These distributions refer to the final equilibrium state for both the porous-elastic and the porous-elastic-plastic material. The difference in the crack length for the two different materials and the positive change of the plastic fluid content as well as the negative change of the elastic fluid content in front of the fracture tip in the elastic-plastic material can again be observed.

\begin{figure}
  \begin{center}
  \newcommand{\fsize}[0]{\Large}
  \begin{minipage}{0.32\textwidth} \hspace{-0.05cm}\resizebox{!}{3.85cm}{
  \begingroup
  \makeatletter
  \providecommand\color[2][]{%
    \GenericError{(gnuplot) \space\space\space\@spaces}{%
      Package color not loaded in conjunction with
      terminal option `colourtext'%
    }{See the gnuplot documentation for explanation.%
    }{Either use 'blacktext' in gnuplot or load the package
      color.sty in LaTeX.}%
    \renewcommand\color[2][]{}%
  }%
  \providecommand\includegraphics[2][]{%
    \GenericError{(gnuplot) \space\space\space\@spaces}{%
      Package graphicx or graphics not loaded%
    }{See the gnuplot documentation for explanation.%
    }{The gnuplot epslatex terminal needs graphicx.sty or graphics.sty.}%
    \renewcommand\includegraphics[2][]{}%
  }%
  \providecommand\rotatebox[2]{#2}%
  \@ifundefined{ifGPcolor}{%
    \newif\ifGPcolor
    \GPcolortrue
  }{}%
  \@ifundefined{ifGPblacktext}{%
    \newif\ifGPblacktext
    \GPblacktexttrue
  }{}%
  \let\gplgaddtomacro\g@addto@macro
  \gdef\gplbacktext{}%
  \gdef\gplfronttext{}%
  \makeatother
  \ifGPblacktext
    \def\colorrgb#1{}%
    \def\colorgray#1{}%
  \else
    \ifGPcolor
      \def\colorrgb#1{\color[rgb]{#1}}%
      \def\colorgray#1{\color[gray]{#1}}%
      \expandafter\def\csname LTw\endcsname{\color{white}}%
      \expandafter\def\csname LTb\endcsname{\color{black}}%
      \expandafter\def\csname LTa\endcsname{\color{black}}%
      \expandafter\def\csname LT0\endcsname{\color[rgb]{1,0,0}}%
      \expandafter\def\csname LT1\endcsname{\color[rgb]{0,1,0}}%
      \expandafter\def\csname LT2\endcsname{\color[rgb]{0,0,1}}%
      \expandafter\def\csname LT3\endcsname{\color[rgb]{1,0,1}}%
      \expandafter\def\csname LT4\endcsname{\color[rgb]{0,1,1}}%
      \expandafter\def\csname LT5\endcsname{\color[rgb]{1,1,0}}%
      \expandafter\def\csname LT6\endcsname{\color[rgb]{0,0,0}}%
      \expandafter\def\csname LT7\endcsname{\color[rgb]{1,0.3,0}}%
      \expandafter\def\csname LT8\endcsname{\color[rgb]{0.5,0.5,0.5}}%
    \else
      \def\colorrgb#1{\color{black}}%
      \def\colorgray#1{\color[gray]{#1}}%
      \expandafter\def\csname LTw\endcsname{\color{white}}%
      \expandafter\def\csname LTb\endcsname{\color{black}}%
      \expandafter\def\csname LTa\endcsname{\color{black}}%
      \expandafter\def\csname LT0\endcsname{\color{black}}%
      \expandafter\def\csname LT1\endcsname{\color{black}}%
      \expandafter\def\csname LT2\endcsname{\color{black}}%
      \expandafter\def\csname LT3\endcsname{\color{black}}%
      \expandafter\def\csname LT4\endcsname{\color{black}}%
      \expandafter\def\csname LT5\endcsname{\color{black}}%
      \expandafter\def\csname LT6\endcsname{\color{black}}%
      \expandafter\def\csname LT7\endcsname{\color{black}}%
      \expandafter\def\csname LT8\endcsname{\color{black}}%
    \fi
  \fi
    \setlength{\unitlength}{0.0500bp}%
    \ifx\gptboxheight\undefined%
      \newlength{\gptboxheight}%
      \newlength{\gptboxwidth}%
      \newsavebox{\gptboxtext}%
    \fi%
    \setlength{\fboxrule}{0.5pt}%
    \setlength{\fboxsep}{1pt}%
  \begin{picture}(9070.00,7086.00)%
    \gplgaddtomacro\gplbacktext{%
      \csname LTb\endcsname
      \put(946,1050){\makebox(0,0)[r]{\strut{}\fsize $0$}}%
      \put(946,2189){\makebox(0,0)[r]{\strut{}\fsize $0.01$}}%
      \put(946,3329){\makebox(0,0)[r]{\strut{}\fsize $0.02$}}%
      \put(946,4468){\makebox(0,0)[r]{\strut{}\fsize $0.03$}}%
      \put(946,5608){\makebox(0,0)[r]{\strut{}\fsize $0.04$}}%
      \put(946,6747){\makebox(0,0)[r]{\strut{}\fsize $0.05$}}%
      \put(1078,602){\makebox(0,0){\strut{}\fsize $0$}}%
      \put(2027,602){\makebox(0,0){\strut{}\fsize $10$}}%
      \put(2977,602){\makebox(0,0){\strut{}\fsize $20$}}%
      \put(3926,602){\makebox(0,0){\strut{}\fsize $30$}}%
      \put(4876,602){\makebox(0,0){\strut{}\fsize $40$}}%
      \put(5825,602){\makebox(0,0){\strut{}\fsize $50$}}%
      \put(6774,602){\makebox(0,0){\strut{}\fsize $60$}}%
      \put(7724,602){\makebox(0,0){\strut{}\fsize $70$}}%
      \put(8673,602){\makebox(0,0){\strut{}\fsize $80$}}%
    }%
    \gplgaddtomacro\gplfronttext{%
      \csname LTb\endcsname
      \put(0,3784){\rotatebox{-270}{\makebox(0,0){\strut{}\fsize $m$/(kg/m$^3$)}}}%
      \put(4875,272){\makebox(0,0){\strut{}\fsize $x$/m}}%
      \csname LTb\endcsname
      \put(2065,6519){\makebox(0,0)[l]{\strut{}\fsize $m$}}%
      \csname LTb\endcsname
      \put(2065,6189){\makebox(0,0)[l]{\strut{}\fsize $m^p$}}%
    }%
    \gplbacktext
    \put(0,0){\includegraphics{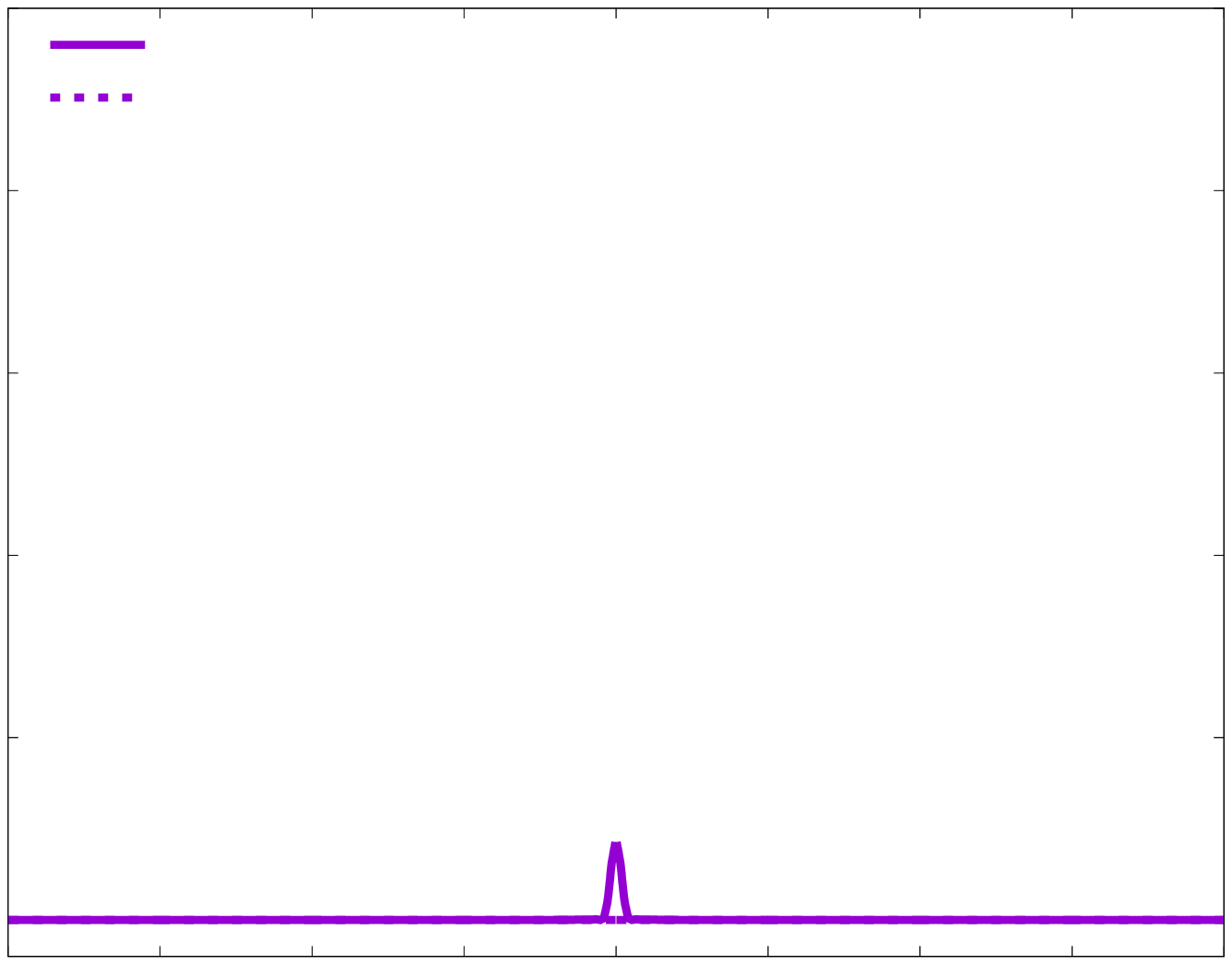}}%
    \gplfronttext
  \end{picture}%
  \endgroup
  }
  \end{minipage} \hspace{-5mm}
  \begin{minipage}{0.32\textwidth} \hspace{-0.05cm}\resizebox{!}{3.85cm}{
\begingroup
  \makeatletter
  \providecommand\color[2][]{%
    \GenericError{(gnuplot) \space\space\space\@spaces}{%
      Package color not loaded in conjunction with
      terminal option `colourtext'%
    }{See the gnuplot documentation for explanation.%
    }{Either use 'blacktext' in gnuplot or load the package
      color.sty in LaTeX.}%
    \renewcommand\color[2][]{}%
  }%
  \providecommand\includegraphics[2][]{%
    \GenericError{(gnuplot) \space\space\space\@spaces}{%
      Package graphicx or graphics not loaded%
    }{See the gnuplot documentation for explanation.%
    }{The gnuplot epslatex terminal needs graphicx.sty or graphics.sty.}%
    \renewcommand\includegraphics[2][]{}%
  }%
  \providecommand\rotatebox[2]{#2}%
  \@ifundefined{ifGPcolor}{%
    \newif\ifGPcolor
    \GPcolortrue
  }{}%
  \@ifundefined{ifGPblacktext}{%
    \newif\ifGPblacktext
    \GPblacktexttrue
  }{}%
  \let\gplgaddtomacro\g@addto@macro
  \gdef\gplbacktext{}%
  \gdef\gplfronttext{}%
  \makeatother
  \ifGPblacktext
    \def\colorrgb#1{}%
    \def\colorgray#1{}%
  \else
    \ifGPcolor
      \def\colorrgb#1{\color[rgb]{#1}}%
      \def\colorgray#1{\color[gray]{#1}}%
      \expandafter\def\csname LTw\endcsname{\color{white}}%
      \expandafter\def\csname LTb\endcsname{\color{black}}%
      \expandafter\def\csname LTa\endcsname{\color{black}}%
      \expandafter\def\csname LT0\endcsname{\color[rgb]{1,0,0}}%
      \expandafter\def\csname LT1\endcsname{\color[rgb]{0,1,0}}%
      \expandafter\def\csname LT2\endcsname{\color[rgb]{0,0,1}}%
      \expandafter\def\csname LT3\endcsname{\color[rgb]{1,0,1}}%
      \expandafter\def\csname LT4\endcsname{\color[rgb]{0,1,1}}%
      \expandafter\def\csname LT5\endcsname{\color[rgb]{1,1,0}}%
      \expandafter\def\csname LT6\endcsname{\color[rgb]{0,0,0}}%
      \expandafter\def\csname LT7\endcsname{\color[rgb]{1,0.3,0}}%
      \expandafter\def\csname LT8\endcsname{\color[rgb]{0.5,0.5,0.5}}%
    \else
      \def\colorrgb#1{\color{black}}%
      \def\colorgray#1{\color[gray]{#1}}%
      \expandafter\def\csname LTw\endcsname{\color{white}}%
      \expandafter\def\csname LTb\endcsname{\color{black}}%
      \expandafter\def\csname LTa\endcsname{\color{black}}%
      \expandafter\def\csname LT0\endcsname{\color{black}}%
      \expandafter\def\csname LT1\endcsname{\color{black}}%
      \expandafter\def\csname LT2\endcsname{\color{black}}%
      \expandafter\def\csname LT3\endcsname{\color{black}}%
      \expandafter\def\csname LT4\endcsname{\color{black}}%
      \expandafter\def\csname LT5\endcsname{\color{black}}%
      \expandafter\def\csname LT6\endcsname{\color{black}}%
      \expandafter\def\csname LT7\endcsname{\color{black}}%
      \expandafter\def\csname LT8\endcsname{\color{black}}%
    \fi
  \fi
    \setlength{\unitlength}{0.0500bp}%
    \ifx\gptboxheight\undefined%
      \newlength{\gptboxheight}%
      \newlength{\gptboxwidth}%
      \newsavebox{\gptboxtext}%
    \fi%
    \setlength{\fboxrule}{0.5pt}%
    \setlength{\fboxsep}{1pt}%
\begin{picture}(9070.00,7086.00)%
    \gplgaddtomacro\gplbacktext{%
      \csname LTb\endcsname
      \put(946,1050){\makebox(0,0)[r]{\strut{}\fsize $0$}}%
      \put(946,2189){\makebox(0,0)[r]{\strut{}\fsize $0.01$}}%
      \put(946,3329){\makebox(0,0)[r]{\strut{}\fsize $0.02$}}%
      \put(946,4468){\makebox(0,0)[r]{\strut{}\fsize $0.03$}}%
      \put(946,5608){\makebox(0,0)[r]{\strut{}\fsize $0.04$}}%
      \put(946,6747){\makebox(0,0)[r]{\strut{}\fsize $0.05$}}%
      \put(1078,602){\makebox(0,0){\strut{}\fsize $0$}}%
      \put(2027,602){\makebox(0,0){\strut{}\fsize $10$}}%
      \put(2977,602){\makebox(0,0){\strut{}\fsize $20$}}%
      \put(3926,602){\makebox(0,0){\strut{}\fsize $30$}}%
      \put(4876,602){\makebox(0,0){\strut{}\fsize $40$}}%
      \put(5825,602){\makebox(0,0){\strut{}\fsize $50$}}%
      \put(6774,602){\makebox(0,0){\strut{}\fsize $60$}}%
      \put(7724,602){\makebox(0,0){\strut{}\fsize $70$}}%
      \put(8673,602){\makebox(0,0){\strut{}\fsize $80$}}%
    }%
    \gplgaddtomacro\gplfronttext{%
      \csname LTb\endcsname
      \put(0,3784){\rotatebox{-270}{\makebox(0,0){\strut{}\fsize $m$/(kg/m$^3$)}}}%
      \put(4875,272){\makebox(0,0){\strut{}\fsize $x$/m}}%
      \csname LTb\endcsname
      \put(2065,6519){\makebox(0,0)[l]{\strut{}\fsize $m$}}%
      \csname LTb\endcsname
      \put(2065,6189){\makebox(0,0)[l]{\strut{}\fsize $m^p$}}%
    }%
    \gplbacktext
    \put(0,0){\includegraphics{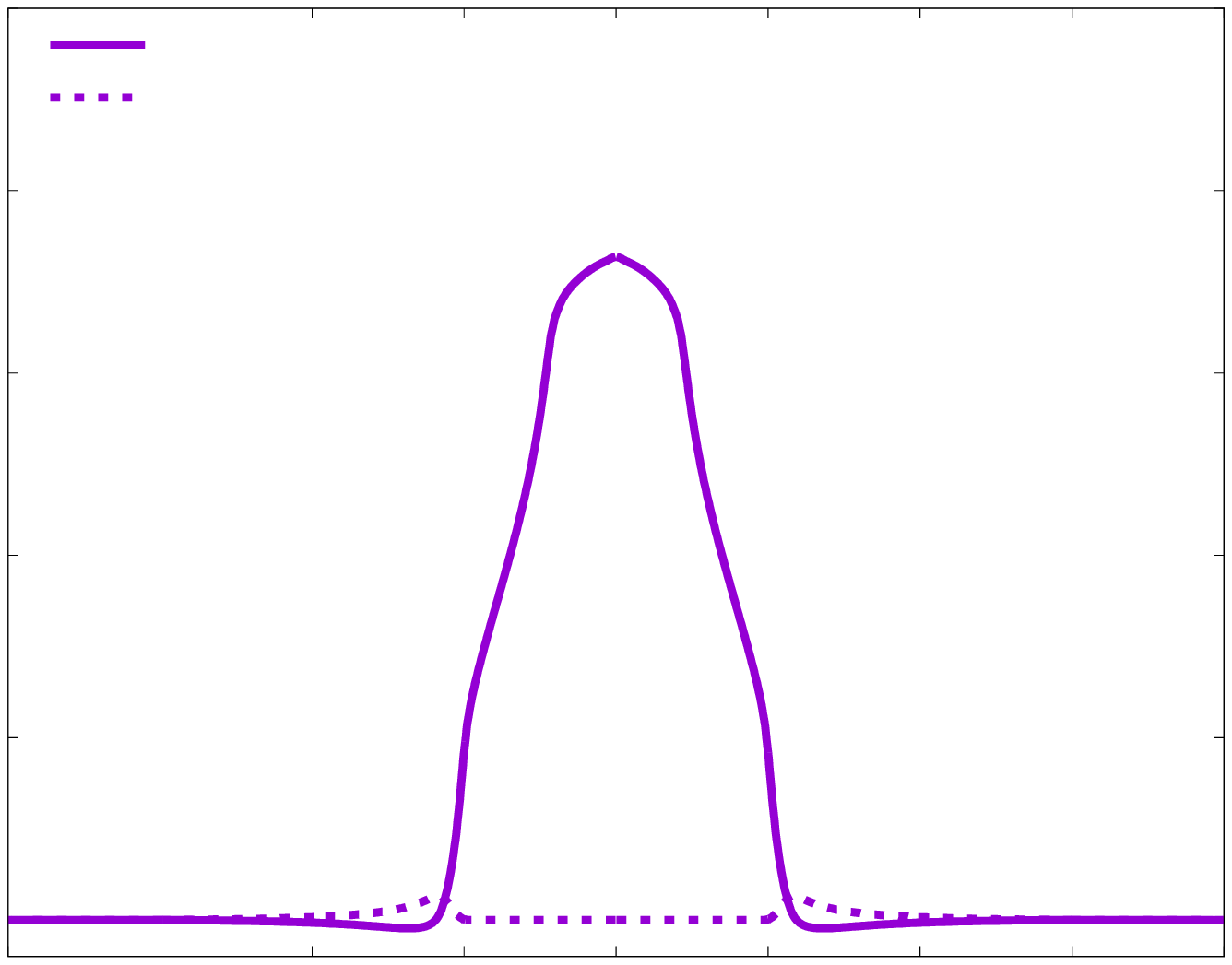}}%
    \gplfronttext
  \end{picture}%
  \endgroup
  }
  \end{minipage} \hspace{-5mm}
  \begin{minipage}{0.32\textwidth} \hspace{-0.05cm}\resizebox{!}{3.85cm}{
\begingroup
  \makeatletter
  \providecommand\color[2][]{%
    \GenericError{(gnuplot) \space\space\space\@spaces}{%
      Package color not loaded in conjunction with
      terminal option `colourtext'%
    }{See the gnuplot documentation for explanation.%
    }{Either use 'blacktext' in gnuplot or load the package
      color.sty in LaTeX.}%
    \renewcommand\color[2][]{}%
  }%
  \providecommand\includegraphics[2][]{%
    \GenericError{(gnuplot) \space\space\space\@spaces}{%
      Package graphicx or graphics not loaded%
    }{See the gnuplot documentation for explanation.%
    }{The gnuplot epslatex terminal needs graphicx.sty or graphics.sty.}%
    \renewcommand\includegraphics[2][]{}%
  }%
  \providecommand\rotatebox[2]{#2}%
  \@ifundefined{ifGPcolor}{%
    \newif\ifGPcolor
    \GPcolortrue
  }{}%
  \@ifundefined{ifGPblacktext}{%
    \newif\ifGPblacktext
    \GPblacktexttrue
  }{}%
  \let\gplgaddtomacro\g@addto@macro
  \gdef\gplbacktext{}%
  \gdef\gplfronttext{}%
  \makeatother
  \ifGPblacktext
    \def\colorrgb#1{}%
    \def\colorgray#1{}%
  \else
    \ifGPcolor
      \def\colorrgb#1{\color[rgb]{#1}}%
      \def\colorgray#1{\color[gray]{#1}}%
      \expandafter\def\csname LTw\endcsname{\color{white}}%
      \expandafter\def\csname LTb\endcsname{\color{black}}%
      \expandafter\def\csname LTa\endcsname{\color{black}}%
      \expandafter\def\csname LT0\endcsname{\color[rgb]{1,0,0}}%
      \expandafter\def\csname LT1\endcsname{\color[rgb]{0,1,0}}%
      \expandafter\def\csname LT2\endcsname{\color[rgb]{0,0,1}}%
      \expandafter\def\csname LT3\endcsname{\color[rgb]{1,0,1}}%
      \expandafter\def\csname LT4\endcsname{\color[rgb]{0,1,1}}%
      \expandafter\def\csname LT5\endcsname{\color[rgb]{1,1,0}}%
      \expandafter\def\csname LT6\endcsname{\color[rgb]{0,0,0}}%
      \expandafter\def\csname LT7\endcsname{\color[rgb]{1,0.3,0}}%
      \expandafter\def\csname LT8\endcsname{\color[rgb]{0.5,0.5,0.5}}%
    \else
      \def\colorrgb#1{\color{black}}%
      \def\colorgray#1{\color[gray]{#1}}%
      \expandafter\def\csname LTw\endcsname{\color{white}}%
      \expandafter\def\csname LTb\endcsname{\color{black}}%
      \expandafter\def\csname LTa\endcsname{\color{black}}%
      \expandafter\def\csname LT0\endcsname{\color{black}}%
      \expandafter\def\csname LT1\endcsname{\color{black}}%
      \expandafter\def\csname LT2\endcsname{\color{black}}%
      \expandafter\def\csname LT3\endcsname{\color{black}}%
      \expandafter\def\csname LT4\endcsname{\color{black}}%
      \expandafter\def\csname LT5\endcsname{\color{black}}%
      \expandafter\def\csname LT6\endcsname{\color{black}}%
      \expandafter\def\csname LT7\endcsname{\color{black}}%
      \expandafter\def\csname LT8\endcsname{\color{black}}%
    \fi
  \fi
    \setlength{\unitlength}{0.0500bp}%
    \ifx\gptboxheight\undefined%
      \newlength{\gptboxheight}%
      \newlength{\gptboxwidth}%
      \newsavebox{\gptboxtext}%
    \fi%
    \setlength{\fboxrule}{0.5pt}%
    \setlength{\fboxsep}{1pt}%
\begin{picture}(9070.00,7086.00)%
    \gplgaddtomacro\gplbacktext{%
      \csname LTb\endcsname
      \put(946,1050){\makebox(0,0)[r]{\strut{}\fsize $0$}}%
      \put(946,2189){\makebox(0,0)[r]{\strut{}\fsize $0.01$}}%
      \put(946,3329){\makebox(0,0)[r]{\strut{}\fsize $0.02$}}%
      \put(946,4468){\makebox(0,0)[r]{\strut{}\fsize $0.03$}}%
      \put(946,5608){\makebox(0,0)[r]{\strut{}\fsize $0.04$}}%
      \put(946,6747){\makebox(0,0)[r]{\strut{}\fsize $0.05$}}%
      \put(1078,602){\makebox(0,0){\strut{}\fsize $0$}}%
      \put(2027,602){\makebox(0,0){\strut{}\fsize $10$}}%
      \put(2977,602){\makebox(0,0){\strut{}\fsize $20$}}%
      \put(3926,602){\makebox(0,0){\strut{}\fsize $30$}}%
      \put(4876,602){\makebox(0,0){\strut{}\fsize $40$}}%
      \put(5825,602){\makebox(0,0){\strut{}\fsize $50$}}%
      \put(6774,602){\makebox(0,0){\strut{}\fsize $60$}}%
      \put(7724,602){\makebox(0,0){\strut{}\fsize $70$}}%
      \put(8673,602){\makebox(0,0){\strut{}\fsize $80$}}%
    }%
    \gplgaddtomacro\gplfronttext{%
      \csname LTb\endcsname
      \put(0,3784){\rotatebox{-270}{\makebox(0,0){\strut{}\fsize $m$/(kg/m$^3$)}}}%
      \put(4875,272){\makebox(0,0){\strut{}\fsize $x$/m}}%
      \csname LTb\endcsname
      \put(2065,6519){\makebox(0,0)[l]{\strut{}\fsize $m$}}%
      \csname LTb\endcsname
      \put(2065,6189){\makebox(0,0)[l]{\strut{}\fsize $m^p$}}%
    }%
    \gplbacktext
    \put(0,0){\includegraphics{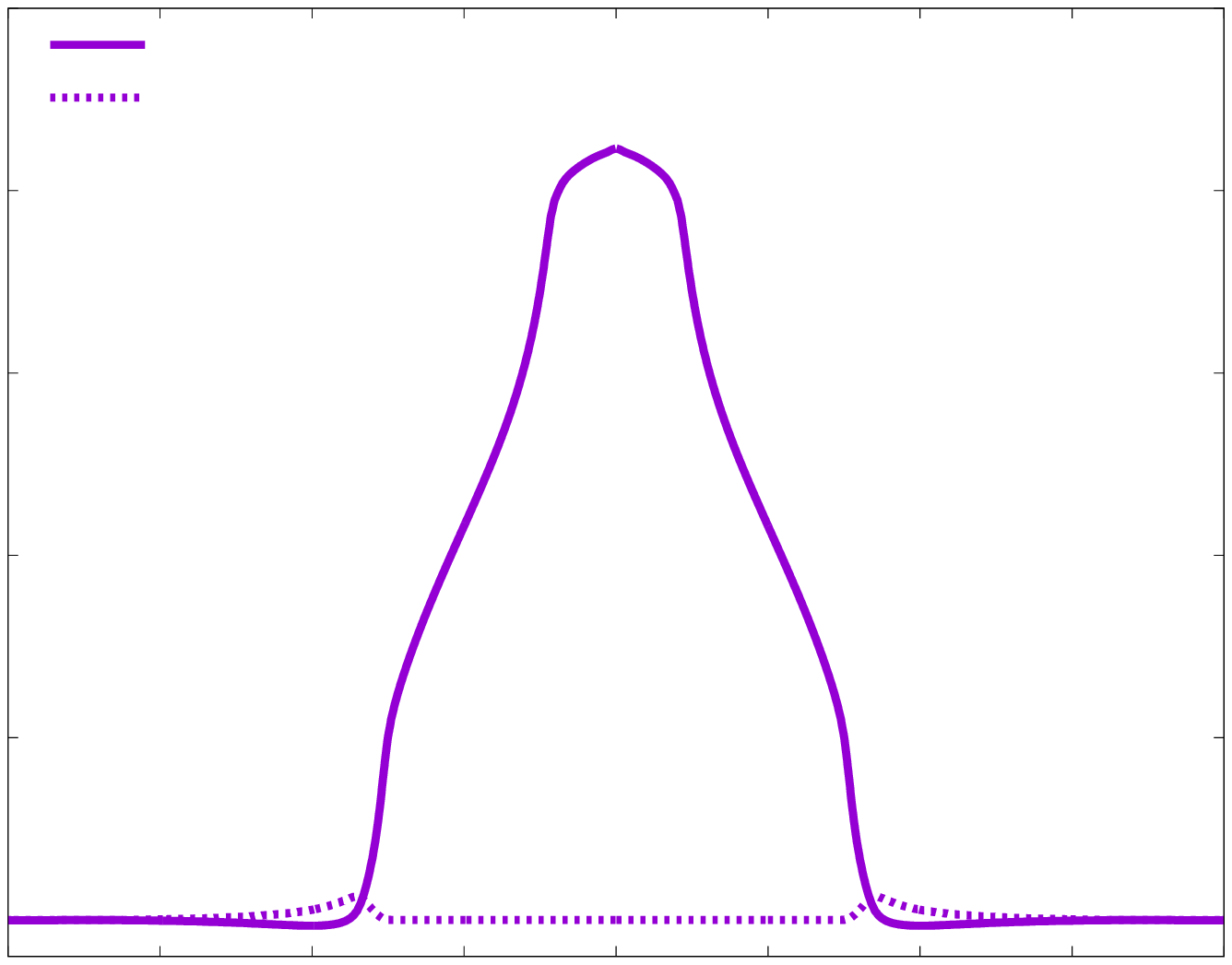}}%
    \gplfronttext
  \end{picture}%
  \endgroup
  }
  \end{minipage}
  \begin{minipage}{0.32\textwidth} \hspace{-0cm}\resizebox{!}{3.8cm}{
\begingroup
  \makeatletter
  \providecommand\color[2][]{%
    \GenericError{(gnuplot) \space\space\space\@spaces}{%
      Package color not loaded in conjunction with
      terminal option `colourtext'%
    }{See the gnuplot documentation for explanation.%
    }{Either use 'blacktext' in gnuplot or load the package
      color.sty in LaTeX.}%
    \renewcommand\color[2][]{}%
  }%
  \providecommand\includegraphics[2][]{%
    \GenericError{(gnuplot) \space\space\space\@spaces}{%
      Package graphicx or graphics not loaded%
    }{See the gnuplot documentation for explanation.%
    }{The gnuplot epslatex terminal needs graphicx.sty or graphics.sty.}%
    \renewcommand\includegraphics[2][]{}%
  }%
  \providecommand\rotatebox[2]{#2}%
  \@ifundefined{ifGPcolor}{%
    \newif\ifGPcolor
    \GPcolortrue
  }{}%
  \@ifundefined{ifGPblacktext}{%
    \newif\ifGPblacktext
    \GPblacktexttrue
  }{}%
  \let\gplgaddtomacro\g@addto@macro
  \gdef\gplbacktext{}%
  \gdef\gplfronttext{}%
  \makeatother
  \ifGPblacktext
    \def\colorrgb#1{}%
    \def\colorgray#1{}%
  \else
    \ifGPcolor
      \def\colorrgb#1{\color[rgb]{#1}}%
      \def\colorgray#1{\color[gray]{#1}}%
      \expandafter\def\csname LTw\endcsname{\color{white}}%
      \expandafter\def\csname LTb\endcsname{\color{black}}%
      \expandafter\def\csname LTa\endcsname{\color{black}}%
      \expandafter\def\csname LT0\endcsname{\color[rgb]{1,0,0}}%
      \expandafter\def\csname LT1\endcsname{\color[rgb]{0,1,0}}%
      \expandafter\def\csname LT2\endcsname{\color[rgb]{0,0,1}}%
      \expandafter\def\csname LT3\endcsname{\color[rgb]{1,0,1}}%
      \expandafter\def\csname LT4\endcsname{\color[rgb]{0,1,1}}%
      \expandafter\def\csname LT5\endcsname{\color[rgb]{1,1,0}}%
      \expandafter\def\csname LT6\endcsname{\color[rgb]{0,0,0}}%
      \expandafter\def\csname LT7\endcsname{\color[rgb]{1,0.3,0}}%
      \expandafter\def\csname LT8\endcsname{\color[rgb]{0.5,0.5,0.5}}%
    \else
      \def\colorrgb#1{\color{black}}%
      \def\colorgray#1{\color[gray]{#1}}%
      \expandafter\def\csname LTw\endcsname{\color{white}}%
      \expandafter\def\csname LTb\endcsname{\color{black}}%
      \expandafter\def\csname LTa\endcsname{\color{black}}%
      \expandafter\def\csname LT0\endcsname{\color{black}}%
      \expandafter\def\csname LT1\endcsname{\color{black}}%
      \expandafter\def\csname LT2\endcsname{\color{black}}%
      \expandafter\def\csname LT3\endcsname{\color{black}}%
      \expandafter\def\csname LT4\endcsname{\color{black}}%
      \expandafter\def\csname LT5\endcsname{\color{black}}%
      \expandafter\def\csname LT6\endcsname{\color{black}}%
      \expandafter\def\csname LT7\endcsname{\color{black}}%
      \expandafter\def\csname LT8\endcsname{\color{black}}%
    \fi
  \fi
    \setlength{\unitlength}{0.0500bp}%
    \ifx\gptboxheight\undefined%
      \newlength{\gptboxheight}%
      \newlength{\gptboxwidth}%
      \newsavebox{\gptboxtext}%
    \fi%
    \setlength{\fboxrule}{0.5pt}%
    \setlength{\fboxsep}{1pt}%
\begin{picture}(9070.00,7086.00)%
    \gplgaddtomacro\gplbacktext{%
      \csname LTb\endcsname
      \put(1078,874){\makebox(0,0)[r]{\strut{}\fsize $-0.02$}}%
      \put(1078,1602){\makebox(0,0)[r]{\strut{}\fsize $0$}}%
      \put(1078,2329){\makebox(0,0)[r]{\strut{}\fsize $0.02$}}%
      \put(1078,3057){\makebox(0,0)[r]{\strut{}\fsize $0.04$}}%
      \put(1078,3785){\makebox(0,0)[r]{\strut{}\fsize $0.06$}}%
      \put(1078,4512){\makebox(0,0)[r]{\strut{}\fsize $0.08$}}%
      \put(1078,5240){\makebox(0,0)[r]{\strut{}\fsize $0.1$}}%
      \put(1078,5967){\makebox(0,0)[r]{\strut{}\fsize $0.12$}}%
      \put(1078,6695){\makebox(0,0)[r]{\strut{}\fsize $0.14$}}%
      \put(1210,654){\makebox(0,0){\strut{}\fsize $0$}}%
      \put(2143,654){\makebox(0,0){\strut{}\fsize $10$}}%
      \put(3076,654){\makebox(0,0){\strut{}\fsize $20$}}%
      \put(4009,654){\makebox(0,0){\strut{}\fsize $30$}}%
      \put(4942,654){\makebox(0,0){\strut{}\fsize $40$}}%
      \put(5874,654){\makebox(0,0){\strut{}\fsize $50$}}%
      \put(6807,654){\makebox(0,0){\strut{}\fsize $60$}}%
      \put(7740,654){\makebox(0,0){\strut{}\fsize $70$}}%
      \put(8673,654){\makebox(0,0){\strut{}\fsize $80$}}%
    }%
    \gplgaddtomacro\gplfronttext{%
      \csname LTb\endcsname
      \put(198,3784){\rotatebox{-270}{\makebox(0,0){\strut{}\fsize $p$/MPa}}}%
      \put(4941,324){\makebox(0,0){\strut{}\fsize $x$/m}}%
    }%
    \gplbacktext
    \put(0,0){\includegraphics{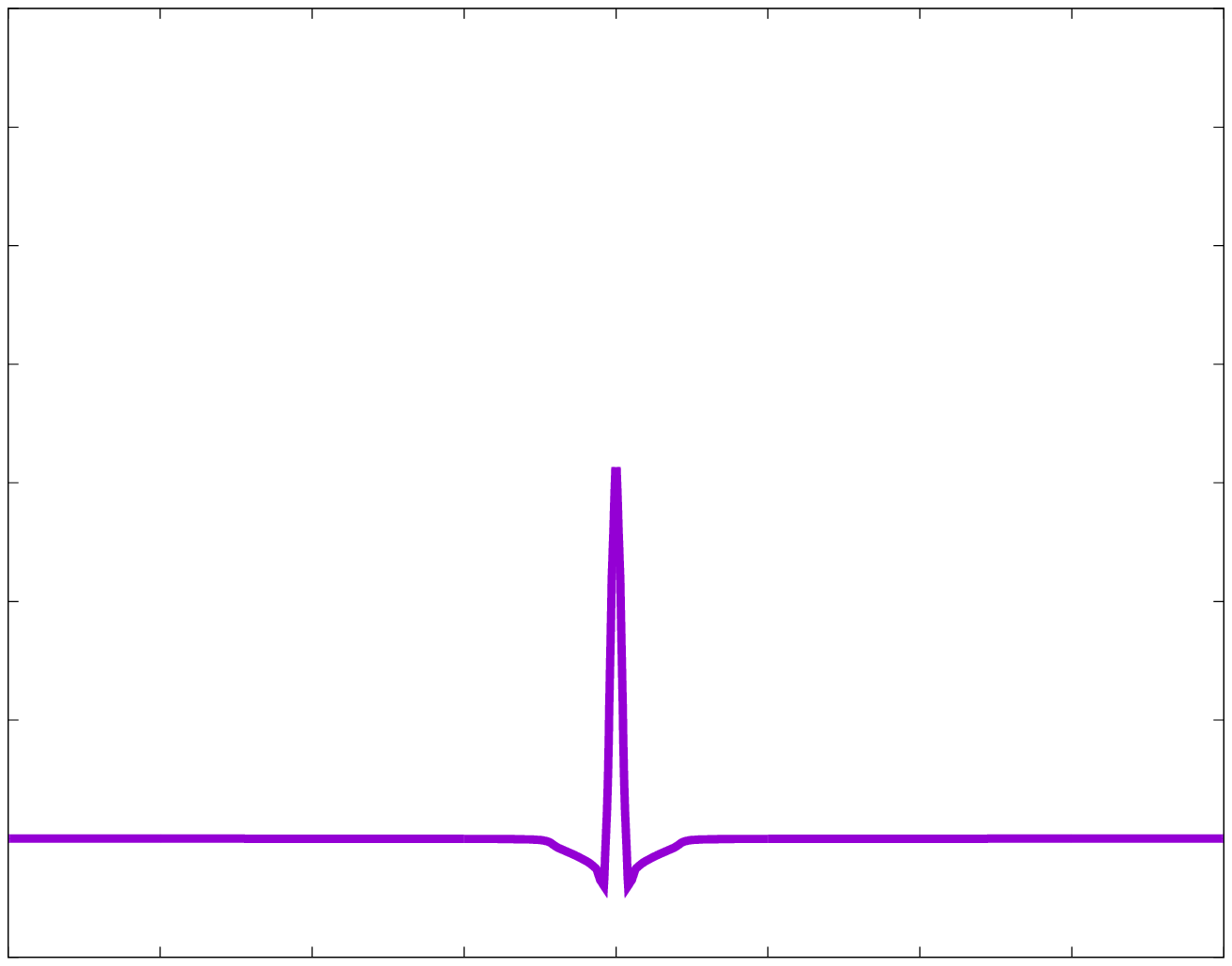}}%
    \gplfronttext
  \end{picture}%
  \endgroup
  }
  \end{minipage} \hspace{-5mm}
  \begin{minipage}{0.32\textwidth} \hspace{-0cm}\resizebox{!}{3.8cm}{
\begingroup
  \makeatletter
  \providecommand\color[2][]{%
    \GenericError{(gnuplot) \space\space\space\@spaces}{%
      Package color not loaded in conjunction with
      terminal option `colourtext'%
    }{See the gnuplot documentation for explanation.%
    }{Either use 'blacktext' in gnuplot or load the package
      color.sty in LaTeX.}%
    \renewcommand\color[2][]{}%
  }%
  \providecommand\includegraphics[2][]{%
    \GenericError{(gnuplot) \space\space\space\@spaces}{%
      Package graphicx or graphics not loaded%
    }{See the gnuplot documentation for explanation.%
    }{The gnuplot epslatex terminal needs graphicx.sty or graphics.sty.}%
    \renewcommand\includegraphics[2][]{}%
  }%
  \providecommand\rotatebox[2]{#2}%
  \@ifundefined{ifGPcolor}{%
    \newif\ifGPcolor
    \GPcolortrue
  }{}%
  \@ifundefined{ifGPblacktext}{%
    \newif\ifGPblacktext
    \GPblacktexttrue
  }{}%
  \let\gplgaddtomacro\g@addto@macro
  \gdef\gplbacktext{}%
  \gdef\gplfronttext{}%
  \makeatother
  \ifGPblacktext
    \def\colorrgb#1{}%
    \def\colorgray#1{}%
  \else
    \ifGPcolor
      \def\colorrgb#1{\color[rgb]{#1}}%
      \def\colorgray#1{\color[gray]{#1}}%
      \expandafter\def\csname LTw\endcsname{\color{white}}%
      \expandafter\def\csname LTb\endcsname{\color{black}}%
      \expandafter\def\csname LTa\endcsname{\color{black}}%
      \expandafter\def\csname LT0\endcsname{\color[rgb]{1,0,0}}%
      \expandafter\def\csname LT1\endcsname{\color[rgb]{0,1,0}}%
      \expandafter\def\csname LT2\endcsname{\color[rgb]{0,0,1}}%
      \expandafter\def\csname LT3\endcsname{\color[rgb]{1,0,1}}%
      \expandafter\def\csname LT4\endcsname{\color[rgb]{0,1,1}}%
      \expandafter\def\csname LT5\endcsname{\color[rgb]{1,1,0}}%
      \expandafter\def\csname LT6\endcsname{\color[rgb]{0,0,0}}%
      \expandafter\def\csname LT7\endcsname{\color[rgb]{1,0.3,0}}%
      \expandafter\def\csname LT8\endcsname{\color[rgb]{0.5,0.5,0.5}}%
    \else
      \def\colorrgb#1{\color{black}}%
      \def\colorgray#1{\color[gray]{#1}}%
      \expandafter\def\csname LTw\endcsname{\color{white}}%
      \expandafter\def\csname LTb\endcsname{\color{black}}%
      \expandafter\def\csname LTa\endcsname{\color{black}}%
      \expandafter\def\csname LT0\endcsname{\color{black}}%
      \expandafter\def\csname LT1\endcsname{\color{black}}%
      \expandafter\def\csname LT2\endcsname{\color{black}}%
      \expandafter\def\csname LT3\endcsname{\color{black}}%
      \expandafter\def\csname LT4\endcsname{\color{black}}%
      \expandafter\def\csname LT5\endcsname{\color{black}}%
      \expandafter\def\csname LT6\endcsname{\color{black}}%
      \expandafter\def\csname LT7\endcsname{\color{black}}%
      \expandafter\def\csname LT8\endcsname{\color{black}}%
    \fi
  \fi
    \setlength{\unitlength}{0.0500bp}%
    \ifx\gptboxheight\undefined%
      \newlength{\gptboxheight}%
      \newlength{\gptboxwidth}%
      \newsavebox{\gptboxtext}%
    \fi%
    \setlength{\fboxrule}{0.5pt}%
    \setlength{\fboxsep}{1pt}%
\begin{picture}(9070.00,7086.00)%
    \gplgaddtomacro\gplbacktext{%
      \csname LTb\endcsname
      \put(1078,874){\makebox(0,0)[r]{\strut{}\fsize $-0.02$}}%
      \put(1078,1602){\makebox(0,0)[r]{\strut{}\fsize $0$}}%
      \put(1078,2329){\makebox(0,0)[r]{\strut{}\fsize $0.02$}}%
      \put(1078,3057){\makebox(0,0)[r]{\strut{}\fsize $0.04$}}%
      \put(1078,3785){\makebox(0,0)[r]{\strut{}\fsize $0.06$}}%
      \put(1078,4512){\makebox(0,0)[r]{\strut{}\fsize $0.08$}}%
      \put(1078,5240){\makebox(0,0)[r]{\strut{}\fsize $0.1$}}%
      \put(1078,5967){\makebox(0,0)[r]{\strut{}\fsize $0.12$}}%
      \put(1078,6695){\makebox(0,0)[r]{\strut{}\fsize $0.14$}}%
      \put(1210,654){\makebox(0,0){\strut{}\fsize $0$}}%
      \put(2143,654){\makebox(0,0){\strut{}\fsize $10$}}%
      \put(3076,654){\makebox(0,0){\strut{}\fsize $20$}}%
      \put(4009,654){\makebox(0,0){\strut{}\fsize $30$}}%
      \put(4942,654){\makebox(0,0){\strut{}\fsize $40$}}%
      \put(5874,654){\makebox(0,0){\strut{}\fsize $50$}}%
      \put(6807,654){\makebox(0,0){\strut{}\fsize $60$}}%
      \put(7740,654){\makebox(0,0){\strut{}\fsize $70$}}%
      \put(8673,654){\makebox(0,0){\strut{}\fsize $80$}}%
    }%
    \gplgaddtomacro\gplfronttext{%
      \csname LTb\endcsname
      \put(198,3784){\rotatebox{-270}{\makebox(0,0){\strut{}\fsize $p$/MPa}}}%
      \put(4941,324){\makebox(0,0){\strut{}\fsize $x$/m}}%
    }%
    \gplbacktext
    \put(0,0){\includegraphics{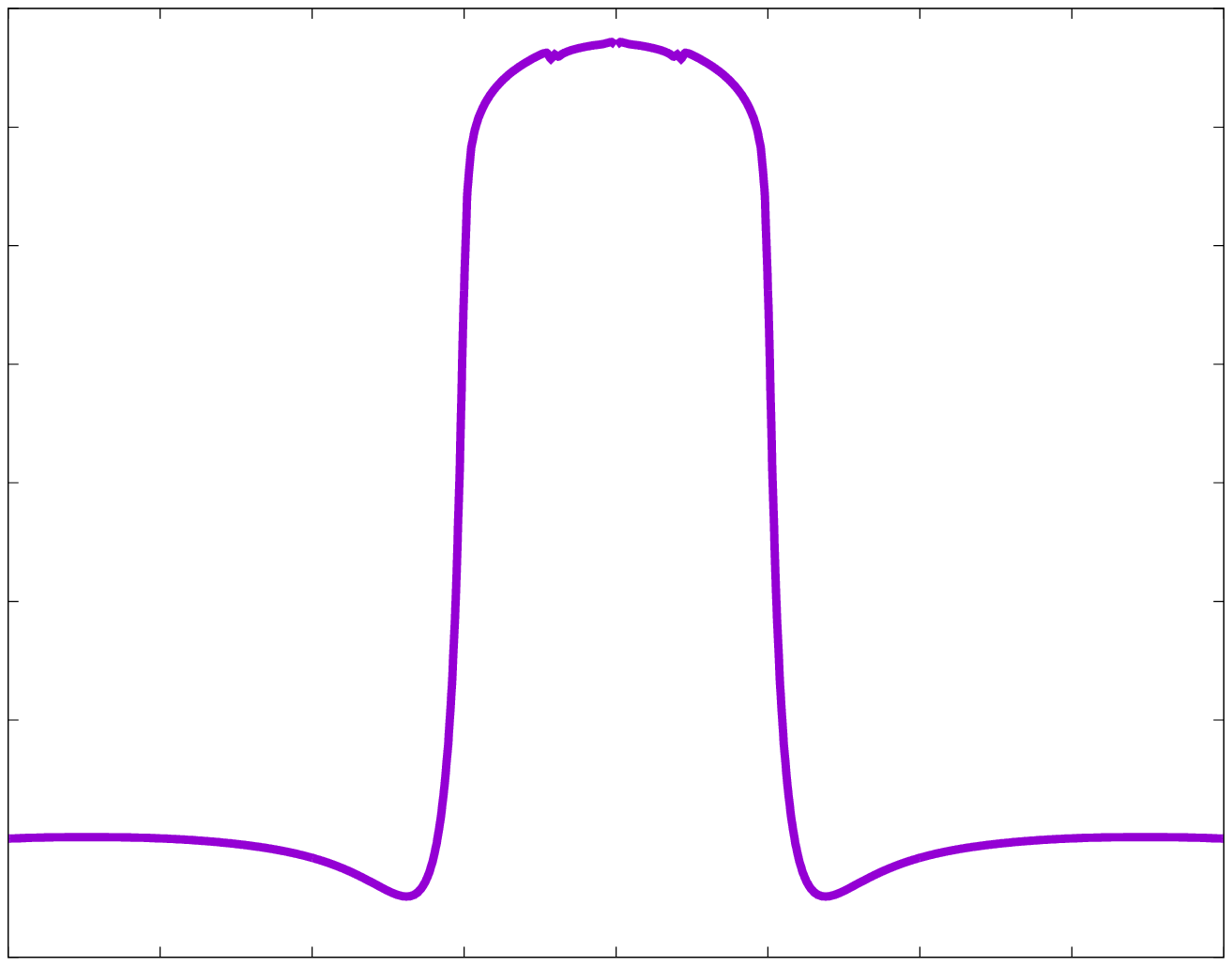}}%
    \gplfronttext
  \end{picture}%
  \endgroup
  }
  \end{minipage} \hspace{-5mm}
  \begin{minipage}{0.32\textwidth} \hspace{-0cm}\resizebox{!}{3.8cm}{
\begingroup
  \makeatletter
  \providecommand\color[2][]{%
    \GenericError{(gnuplot) \space\space\space\@spaces}{%
      Package color not loaded in conjunction with
      terminal option `colourtext'%
    }{See the gnuplot documentation for explanation.%
    }{Either use 'blacktext' in gnuplot or load the package
      color.sty in LaTeX.}%
    \renewcommand\color[2][]{}%
  }%
  \providecommand\includegraphics[2][]{%
    \GenericError{(gnuplot) \space\space\space\@spaces}{%
      Package graphicx or graphics not loaded%
    }{See the gnuplot documentation for explanation.%
    }{The gnuplot epslatex terminal needs graphicx.sty or graphics.sty.}%
    \renewcommand\includegraphics[2][]{}%
  }%
  \providecommand\rotatebox[2]{#2}%
  \@ifundefined{ifGPcolor}{%
    \newif\ifGPcolor
    \GPcolortrue
  }{}%
  \@ifundefined{ifGPblacktext}{%
    \newif\ifGPblacktext
    \GPblacktexttrue
  }{}%
  \let\gplgaddtomacro\g@addto@macro
  \gdef\gplbacktext{}%
  \gdef\gplfronttext{}%
  \makeatother
  \ifGPblacktext
    \def\colorrgb#1{}%
    \def\colorgray#1{}%
  \else
    \ifGPcolor
      \def\colorrgb#1{\color[rgb]{#1}}%
      \def\colorgray#1{\color[gray]{#1}}%
      \expandafter\def\csname LTw\endcsname{\color{white}}%
      \expandafter\def\csname LTb\endcsname{\color{black}}%
      \expandafter\def\csname LTa\endcsname{\color{black}}%
      \expandafter\def\csname LT0\endcsname{\color[rgb]{1,0,0}}%
      \expandafter\def\csname LT1\endcsname{\color[rgb]{0,1,0}}%
      \expandafter\def\csname LT2\endcsname{\color[rgb]{0,0,1}}%
      \expandafter\def\csname LT3\endcsname{\color[rgb]{1,0,1}}%
      \expandafter\def\csname LT4\endcsname{\color[rgb]{0,1,1}}%
      \expandafter\def\csname LT5\endcsname{\color[rgb]{1,1,0}}%
      \expandafter\def\csname LT6\endcsname{\color[rgb]{0,0,0}}%
      \expandafter\def\csname LT7\endcsname{\color[rgb]{1,0.3,0}}%
      \expandafter\def\csname LT8\endcsname{\color[rgb]{0.5,0.5,0.5}}%
    \else
      \def\colorrgb#1{\color{black}}%
      \def\colorgray#1{\color[gray]{#1}}%
      \expandafter\def\csname LTw\endcsname{\color{white}}%
      \expandafter\def\csname LTb\endcsname{\color{black}}%
      \expandafter\def\csname LTa\endcsname{\color{black}}%
      \expandafter\def\csname LT0\endcsname{\color{black}}%
      \expandafter\def\csname LT1\endcsname{\color{black}}%
      \expandafter\def\csname LT2\endcsname{\color{black}}%
      \expandafter\def\csname LT3\endcsname{\color{black}}%
      \expandafter\def\csname LT4\endcsname{\color{black}}%
      \expandafter\def\csname LT5\endcsname{\color{black}}%
      \expandafter\def\csname LT6\endcsname{\color{black}}%
      \expandafter\def\csname LT7\endcsname{\color{black}}%
      \expandafter\def\csname LT8\endcsname{\color{black}}%
    \fi
  \fi
    \setlength{\unitlength}{0.0500bp}%
    \ifx\gptboxheight\undefined%
      \newlength{\gptboxheight}%
      \newlength{\gptboxwidth}%
      \newsavebox{\gptboxtext}%
    \fi%
    \setlength{\fboxrule}{0.5pt}%
    \setlength{\fboxsep}{1pt}%
\begin{picture}(9070.00,7086.00)%
    \gplgaddtomacro\gplbacktext{%
      \csname LTb\endcsname
      \put(1078,874){\makebox(0,0)[r]{\strut{}\fsize $-0.02$}}%
      \put(1078,1602){\makebox(0,0)[r]{\strut{}\fsize $0$}}%
      \put(1078,2329){\makebox(0,0)[r]{\strut{}\fsize $0.02$}}%
      \put(1078,3057){\makebox(0,0)[r]{\strut{}\fsize $0.04$}}%
      \put(1078,3785){\makebox(0,0)[r]{\strut{}\fsize $0.06$}}%
      \put(1078,4512){\makebox(0,0)[r]{\strut{}\fsize $0.08$}}%
      \put(1078,5240){\makebox(0,0)[r]{\strut{}\fsize $0.1$}}%
      \put(1078,5967){\makebox(0,0)[r]{\strut{}\fsize $0.12$}}%
      \put(1078,6695){\makebox(0,0)[r]{\strut{}\fsize $0.14$}}%
      \put(1210,654){\makebox(0,0){\strut{}\fsize $0$}}%
      \put(2143,654){\makebox(0,0){\strut{}\fsize $10$}}%
      \put(3076,654){\makebox(0,0){\strut{}\fsize $20$}}%
      \put(4009,654){\makebox(0,0){\strut{}\fsize $30$}}%
      \put(4942,654){\makebox(0,0){\strut{}\fsize $40$}}%
      \put(5874,654){\makebox(0,0){\strut{}\fsize $50$}}%
      \put(6807,654){\makebox(0,0){\strut{}\fsize $60$}}%
      \put(7740,654){\makebox(0,0){\strut{}\fsize $70$}}%
      \put(8673,654){\makebox(0,0){\strut{}\fsize $80$}}%
    }%
    \gplgaddtomacro\gplfronttext{%
      \csname LTb\endcsname
      \put(198,3784){\rotatebox{-270}{\makebox(0,0){\strut{}\fsize $p$/MPa}}}%
      \put(4941,324){\makebox(0,0){\strut{}\fsize $x$/m}}%
    }%
    \gplbacktext
    \put(0,0){\includegraphics{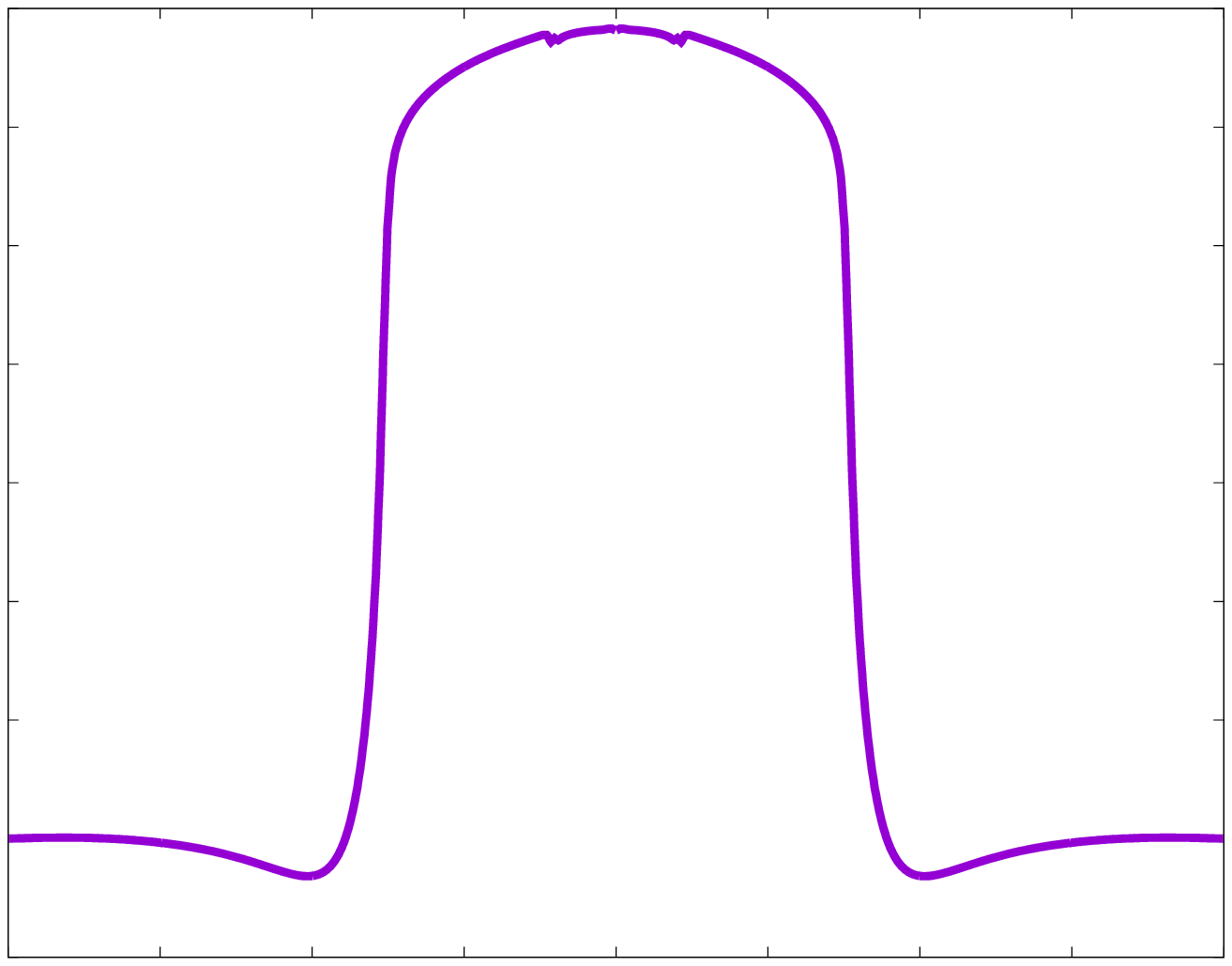}}%
    \gplfronttext
  \end{picture}%
  \endgroup
  }
  \end{minipage}
  \begin{minipage}{0.32\textwidth} \hspace{1mm}\resizebox{!}{3.8cm}{
\begingroup
  \makeatletter
  \providecommand\color[2][]{%
    \GenericError{(gnuplot) \space\space\space\@spaces}{%
      Package color not loaded in conjunction with
      terminal option `colourtext'%
    }{See the gnuplot documentation for explanation.%
    }{Either use 'blacktext' in gnuplot or load the package
      color.sty in LaTeX.}%
    \renewcommand\color[2][]{}%
  }%
  \providecommand\includegraphics[2][]{%
    \GenericError{(gnuplot) \space\space\space\@spaces}{%
      Package graphicx or graphics not loaded%
    }{See the gnuplot documentation for explanation.%
    }{The gnuplot epslatex terminal needs graphicx.sty or graphics.sty.}%
    \renewcommand\includegraphics[2][]{}%
  }%
  \providecommand\rotatebox[2]{#2}%
  \@ifundefined{ifGPcolor}{%
    \newif\ifGPcolor
    \GPcolortrue
  }{}%
  \@ifundefined{ifGPblacktext}{%
    \newif\ifGPblacktext
    \GPblacktexttrue
  }{}%
  \let\gplgaddtomacro\g@addto@macro
  \gdef\gplbacktext{}%
  \gdef\gplfronttext{}%
  \makeatother
  \ifGPblacktext
    \def\colorrgb#1{}%
    \def\colorgray#1{}%
  \else
    \ifGPcolor
      \def\colorrgb#1{\color[rgb]{#1}}%
      \def\colorgray#1{\color[gray]{#1}}%
      \expandafter\def\csname LTw\endcsname{\color{white}}%
      \expandafter\def\csname LTb\endcsname{\color{black}}%
      \expandafter\def\csname LTa\endcsname{\color{black}}%
      \expandafter\def\csname LT0\endcsname{\color[rgb]{1,0,0}}%
      \expandafter\def\csname LT1\endcsname{\color[rgb]{0,1,0}}%
      \expandafter\def\csname LT2\endcsname{\color[rgb]{0,0,1}}%
      \expandafter\def\csname LT3\endcsname{\color[rgb]{1,0,1}}%
      \expandafter\def\csname LT4\endcsname{\color[rgb]{0,1,1}}%
      \expandafter\def\csname LT5\endcsname{\color[rgb]{1,1,0}}%
      \expandafter\def\csname LT6\endcsname{\color[rgb]{0,0,0}}%
      \expandafter\def\csname LT7\endcsname{\color[rgb]{1,0.3,0}}%
      \expandafter\def\csname LT8\endcsname{\color[rgb]{0.5,0.5,0.5}}%
    \else
      \def\colorrgb#1{\color{black}}%
      \def\colorgray#1{\color[gray]{#1}}%
      \expandafter\def\csname LTw\endcsname{\color{white}}%
      \expandafter\def\csname LTb\endcsname{\color{black}}%
      \expandafter\def\csname LTa\endcsname{\color{black}}%
      \expandafter\def\csname LT0\endcsname{\color{black}}%
      \expandafter\def\csname LT1\endcsname{\color{black}}%
      \expandafter\def\csname LT2\endcsname{\color{black}}%
      \expandafter\def\csname LT3\endcsname{\color{black}}%
      \expandafter\def\csname LT4\endcsname{\color{black}}%
      \expandafter\def\csname LT5\endcsname{\color{black}}%
      \expandafter\def\csname LT6\endcsname{\color{black}}%
      \expandafter\def\csname LT7\endcsname{\color{black}}%
      \expandafter\def\csname LT8\endcsname{\color{black}}%
    \fi
  \fi
    \setlength{\unitlength}{0.0500bp}%
    \ifx\gptboxheight\undefined%
      \newlength{\gptboxheight}%
      \newlength{\gptboxwidth}%
      \newsavebox{\gptboxtext}%
    \fi%
    \setlength{\fboxrule}{0.5pt}%
    \setlength{\fboxsep}{1pt}%
\begin{picture}(9070.00,7086.00)%
    \gplgaddtomacro\gplbacktext{%
      \csname LTb\endcsname
      \put(682,719){\makebox(0,0)[r]{\strut{}\fsize 0}}%
      \put(682,2252){\makebox(0,0)[r]{\strut{}\fsize 5}}%
      \put(682,3785){\makebox(0,0)[r]{\strut{}\fsize 10}}%
      \put(682,5317){\makebox(0,0)[r]{\strut{}\fsize 15}}%
      \put(682,6850){\makebox(0,0)[r]{\strut{}\fsize 20}}%
      \put(814,499){\makebox(0,0){\strut{}\fsize $0$}}%
      \put(1796,499){\makebox(0,0){\strut{}\fsize $10$}}%
      \put(2779,499){\makebox(0,0){\strut{}\fsize $20$}}%
      \put(3761,499){\makebox(0,0){\strut{}\fsize $30$}}%
      \put(4744,499){\makebox(0,0){\strut{}\fsize $40$}}%
      \put(5726,499){\makebox(0,0){\strut{}\fsize $50$}}%
      \put(6708,499){\makebox(0,0){\strut{}\fsize $60$}}%
      \put(7691,499){\makebox(0,0){\strut{}\fsize $70$}}%
      \put(8673,499){\makebox(0,0){\strut{}\fsize $80$}}%
    }%
    \gplgaddtomacro\gplfronttext{%
      \csname LTb\endcsname
      \put(198,3784){\rotatebox{-270}{\makebox(0,0){\strut{}\fsize $w$/($1 \cdot 10^{-6}$ m)}}}%
      \put(4743,169){\makebox(0,0){\strut{}\fsize $x$/m}}%
    }%
    \gplbacktext
    \put(0,0){\includegraphics{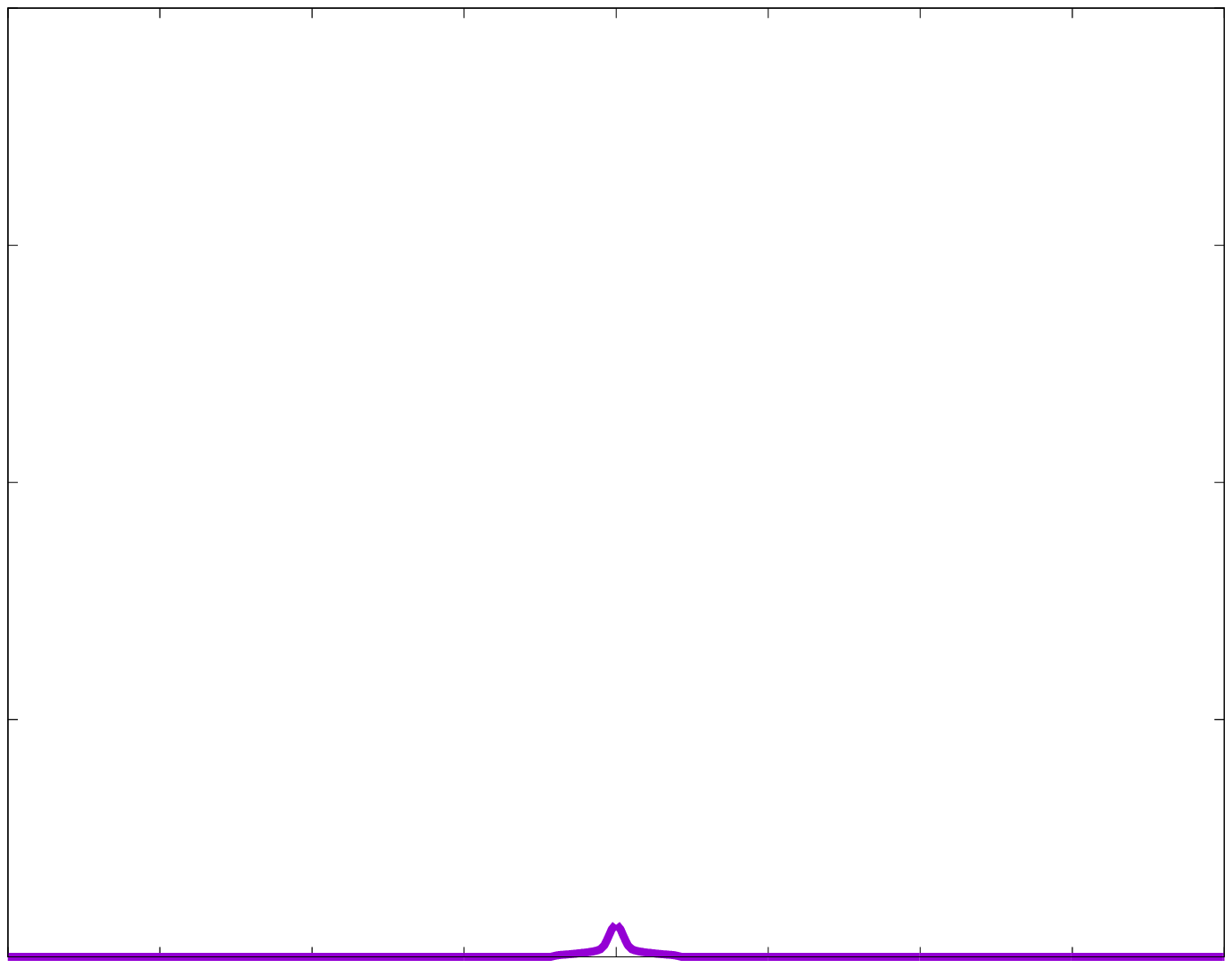}}%
    \gplfronttext
  \end{picture}%
  \endgroup
  }
  \end{minipage} \hspace{-5mm}
  \begin{minipage}{0.32\textwidth} \hspace{1mm}\resizebox{!}{3.8cm}{
\begingroup
  \makeatletter
  \providecommand\color[2][]{%
    \GenericError{(gnuplot) \space\space\space\@spaces}{%
      Package color not loaded in conjunction with
      terminal option `colourtext'%
    }{See the gnuplot documentation for explanation.%
    }{Either use 'blacktext' in gnuplot or load the package
      color.sty in LaTeX.}%
    \renewcommand\color[2][]{}%
  }%
  \providecommand\includegraphics[2][]{%
    \GenericError{(gnuplot) \space\space\space\@spaces}{%
      Package graphicx or graphics not loaded%
    }{See the gnuplot documentation for explanation.%
    }{The gnuplot epslatex terminal needs graphicx.sty or graphics.sty.}%
    \renewcommand\includegraphics[2][]{}%
  }%
  \providecommand\rotatebox[2]{#2}%
  \@ifundefined{ifGPcolor}{%
    \newif\ifGPcolor
    \GPcolortrue
  }{}%
  \@ifundefined{ifGPblacktext}{%
    \newif\ifGPblacktext
    \GPblacktexttrue
  }{}%
  \let\gplgaddtomacro\g@addto@macro
  \gdef\gplbacktext{}%
  \gdef\gplfronttext{}%
  \makeatother
  \ifGPblacktext
    \def\colorrgb#1{}%
    \def\colorgray#1{}%
  \else
    \ifGPcolor
      \def\colorrgb#1{\color[rgb]{#1}}%
      \def\colorgray#1{\color[gray]{#1}}%
      \expandafter\def\csname LTw\endcsname{\color{white}}%
      \expandafter\def\csname LTb\endcsname{\color{black}}%
      \expandafter\def\csname LTa\endcsname{\color{black}}%
      \expandafter\def\csname LT0\endcsname{\color[rgb]{1,0,0}}%
      \expandafter\def\csname LT1\endcsname{\color[rgb]{0,1,0}}%
      \expandafter\def\csname LT2\endcsname{\color[rgb]{0,0,1}}%
      \expandafter\def\csname LT3\endcsname{\color[rgb]{1,0,1}}%
      \expandafter\def\csname LT4\endcsname{\color[rgb]{0,1,1}}%
      \expandafter\def\csname LT5\endcsname{\color[rgb]{1,1,0}}%
      \expandafter\def\csname LT6\endcsname{\color[rgb]{0,0,0}}%
      \expandafter\def\csname LT7\endcsname{\color[rgb]{1,0.3,0}}%
      \expandafter\def\csname LT8\endcsname{\color[rgb]{0.5,0.5,0.5}}%
    \else
      \def\colorrgb#1{\color{black}}%
      \def\colorgray#1{\color[gray]{#1}}%
      \expandafter\def\csname LTw\endcsname{\color{white}}%
      \expandafter\def\csname LTb\endcsname{\color{black}}%
      \expandafter\def\csname LTa\endcsname{\color{black}}%
      \expandafter\def\csname LT0\endcsname{\color{black}}%
      \expandafter\def\csname LT1\endcsname{\color{black}}%
      \expandafter\def\csname LT2\endcsname{\color{black}}%
      \expandafter\def\csname LT3\endcsname{\color{black}}%
      \expandafter\def\csname LT4\endcsname{\color{black}}%
      \expandafter\def\csname LT5\endcsname{\color{black}}%
      \expandafter\def\csname LT6\endcsname{\color{black}}%
      \expandafter\def\csname LT7\endcsname{\color{black}}%
      \expandafter\def\csname LT8\endcsname{\color{black}}%
    \fi
  \fi
    \setlength{\unitlength}{0.0500bp}%
    \ifx\gptboxheight\undefined%
      \newlength{\gptboxheight}%
      \newlength{\gptboxwidth}%
      \newsavebox{\gptboxtext}%
    \fi%
    \setlength{\fboxrule}{0.5pt}%
    \setlength{\fboxsep}{1pt}%
\begin{picture}(9070.00,7086.00)%
    \gplgaddtomacro\gplbacktext{%
      \csname LTb\endcsname
      \put(682,719){\makebox(0,0)[r]{\strut{}\fsize 0}}%
      \put(682,2252){\makebox(0,0)[r]{\strut{}\fsize 5}}%
      \put(682,3785){\makebox(0,0)[r]{\strut{}\fsize 10}}%
      \put(682,5317){\makebox(0,0)[r]{\strut{}\fsize 15}}%
      \put(682,6850){\makebox(0,0)[r]{\strut{}\fsize 20}}%
      \put(814,499){\makebox(0,0){\strut{}\fsize $0$}}%
      \put(1796,499){\makebox(0,0){\strut{}\fsize $10$}}%
      \put(2779,499){\makebox(0,0){\strut{}\fsize $20$}}%
      \put(3761,499){\makebox(0,0){\strut{}\fsize $30$}}%
      \put(4744,499){\makebox(0,0){\strut{}\fsize $40$}}%
      \put(5726,499){\makebox(0,0){\strut{}\fsize $50$}}%
      \put(6708,499){\makebox(0,0){\strut{}\fsize $60$}}%
      \put(7691,499){\makebox(0,0){\strut{}\fsize $70$}}%
      \put(8673,499){\makebox(0,0){\strut{}\fsize $80$}}%
    }%
    \gplgaddtomacro\gplfronttext{%
      \csname LTb\endcsname
      \put(198,3784){\rotatebox{-270}{\makebox(0,0){\strut{}\fsize $w$/($1 \cdot 10^{-6}$ m)}}}%
      \put(4743,169){\makebox(0,0){\strut{}\fsize $x$/m}}%
    }%
    \gplbacktext
    \put(0,0){\includegraphics{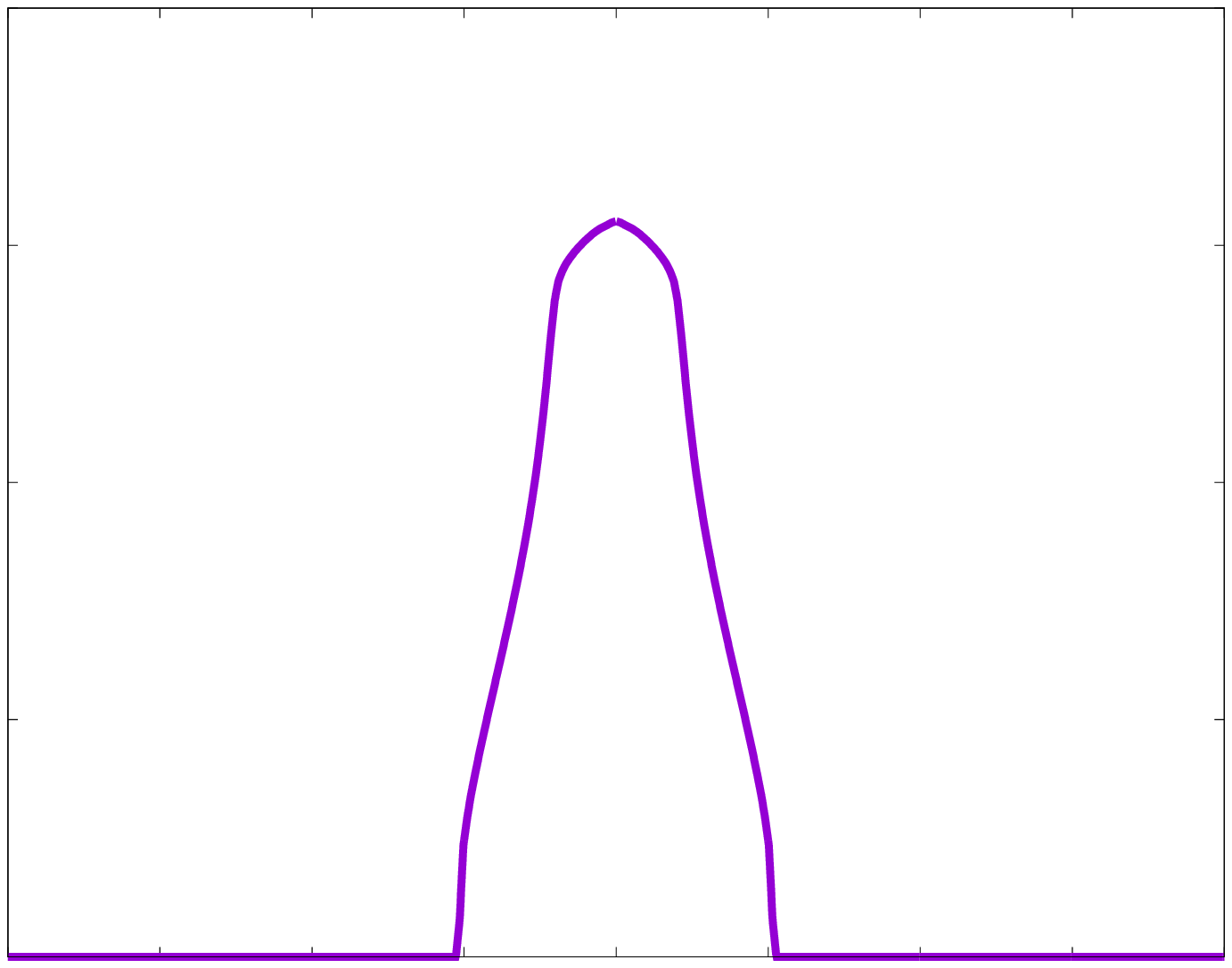}}%
    \gplfronttext
  \end{picture}%
  \endgroup
  }
  \end{minipage} \hspace{-5mm}
  \begin{minipage}{0.32\textwidth} \hspace{1mm}\resizebox{!}{3.8cm}{
\begingroup
  \makeatletter
  \providecommand\color[2][]{%
    \GenericError{(gnuplot) \space\space\space\@spaces}{%
      Package color not loaded in conjunction with
      terminal option `colourtext'%
    }{See the gnuplot documentation for explanation.%
    }{Either use 'blacktext' in gnuplot or load the package
      color.sty in LaTeX.}%
    \renewcommand\color[2][]{}%
  }%
  \providecommand\includegraphics[2][]{%
    \GenericError{(gnuplot) \space\space\space\@spaces}{%
      Package graphicx or graphics not loaded%
    }{See the gnuplot documentation for explanation.%
    }{The gnuplot epslatex terminal needs graphicx.sty or graphics.sty.}%
    \renewcommand\includegraphics[2][]{}%
  }%
  \providecommand\rotatebox[2]{#2}%
  \@ifundefined{ifGPcolor}{%
    \newif\ifGPcolor
    \GPcolortrue
  }{}%
  \@ifundefined{ifGPblacktext}{%
    \newif\ifGPblacktext
    \GPblacktexttrue
  }{}%
  \let\gplgaddtomacro\g@addto@macro
  \gdef\gplbacktext{}%
  \gdef\gplfronttext{}%
  \makeatother
  \ifGPblacktext
    \def\colorrgb#1{}%
    \def\colorgray#1{}%
  \else
    \ifGPcolor
      \def\colorrgb#1{\color[rgb]{#1}}%
      \def\colorgray#1{\color[gray]{#1}}%
      \expandafter\def\csname LTw\endcsname{\color{white}}%
      \expandafter\def\csname LTb\endcsname{\color{black}}%
      \expandafter\def\csname LTa\endcsname{\color{black}}%
      \expandafter\def\csname LT0\endcsname{\color[rgb]{1,0,0}}%
      \expandafter\def\csname LT1\endcsname{\color[rgb]{0,1,0}}%
      \expandafter\def\csname LT2\endcsname{\color[rgb]{0,0,1}}%
      \expandafter\def\csname LT3\endcsname{\color[rgb]{1,0,1}}%
      \expandafter\def\csname LT4\endcsname{\color[rgb]{0,1,1}}%
      \expandafter\def\csname LT5\endcsname{\color[rgb]{1,1,0}}%
      \expandafter\def\csname LT6\endcsname{\color[rgb]{0,0,0}}%
      \expandafter\def\csname LT7\endcsname{\color[rgb]{1,0.3,0}}%
      \expandafter\def\csname LT8\endcsname{\color[rgb]{0.5,0.5,0.5}}%
    \else
      \def\colorrgb#1{\color{black}}%
      \def\colorgray#1{\color[gray]{#1}}%
      \expandafter\def\csname LTw\endcsname{\color{white}}%
      \expandafter\def\csname LTb\endcsname{\color{black}}%
      \expandafter\def\csname LTa\endcsname{\color{black}}%
      \expandafter\def\csname LT0\endcsname{\color{black}}%
      \expandafter\def\csname LT1\endcsname{\color{black}}%
      \expandafter\def\csname LT2\endcsname{\color{black}}%
      \expandafter\def\csname LT3\endcsname{\color{black}}%
      \expandafter\def\csname LT4\endcsname{\color{black}}%
      \expandafter\def\csname LT5\endcsname{\color{black}}%
      \expandafter\def\csname LT6\endcsname{\color{black}}%
      \expandafter\def\csname LT7\endcsname{\color{black}}%
      \expandafter\def\csname LT8\endcsname{\color{black}}%
    \fi
  \fi
    \setlength{\unitlength}{0.0500bp}%
    \ifx\gptboxheight\undefined%
      \newlength{\gptboxheight}%
      \newlength{\gptboxwidth}%
      \newsavebox{\gptboxtext}%
    \fi%
    \setlength{\fboxrule}{0.5pt}%
    \setlength{\fboxsep}{1pt}%
\begin{picture}(9070.00,7086.00)%
    \gplgaddtomacro\gplbacktext{%
      \csname LTb\endcsname
      \put(682,719){\makebox(0,0)[r]{\strut{}\fsize 0}}%
      \put(682,2252){\makebox(0,0)[r]{\strut{}\fsize 5}}%
      \put(682,3785){\makebox(0,0)[r]{\strut{}\fsize 10}}%
      \put(682,5317){\makebox(0,0)[r]{\strut{}\fsize 15}}%
      \put(682,6850){\makebox(0,0)[r]{\strut{}\fsize 20}}%
      \put(814,499){\makebox(0,0){\strut{}\fsize $0$}}%
      \put(1796,499){\makebox(0,0){\strut{}\fsize $10$}}%
      \put(2779,499){\makebox(0,0){\strut{}\fsize $20$}}%
      \put(3761,499){\makebox(0,0){\strut{}\fsize $30$}}%
      \put(4744,499){\makebox(0,0){\strut{}\fsize $40$}}%
      \put(5726,499){\makebox(0,0){\strut{}\fsize $50$}}%
      \put(6708,499){\makebox(0,0){\strut{}\fsize $60$}}%
      \put(7691,499){\makebox(0,0){\strut{}\fsize $70$}}%
      \put(8673,499){\makebox(0,0){\strut{}\fsize $80$}}%
    }%
    \gplgaddtomacro\gplfronttext{%
      \csname LTb\endcsname
      \put(198,3784){\rotatebox{-270}{\makebox(0,0){\strut{}\fsize $w$/($1 \cdot 10^{-6}$ m)}}}%
      \put(4743,169){\makebox(0,0){\strut{}\fsize $x$/m}}%
    }%
    \gplbacktext
    \put(0,0){\includegraphics{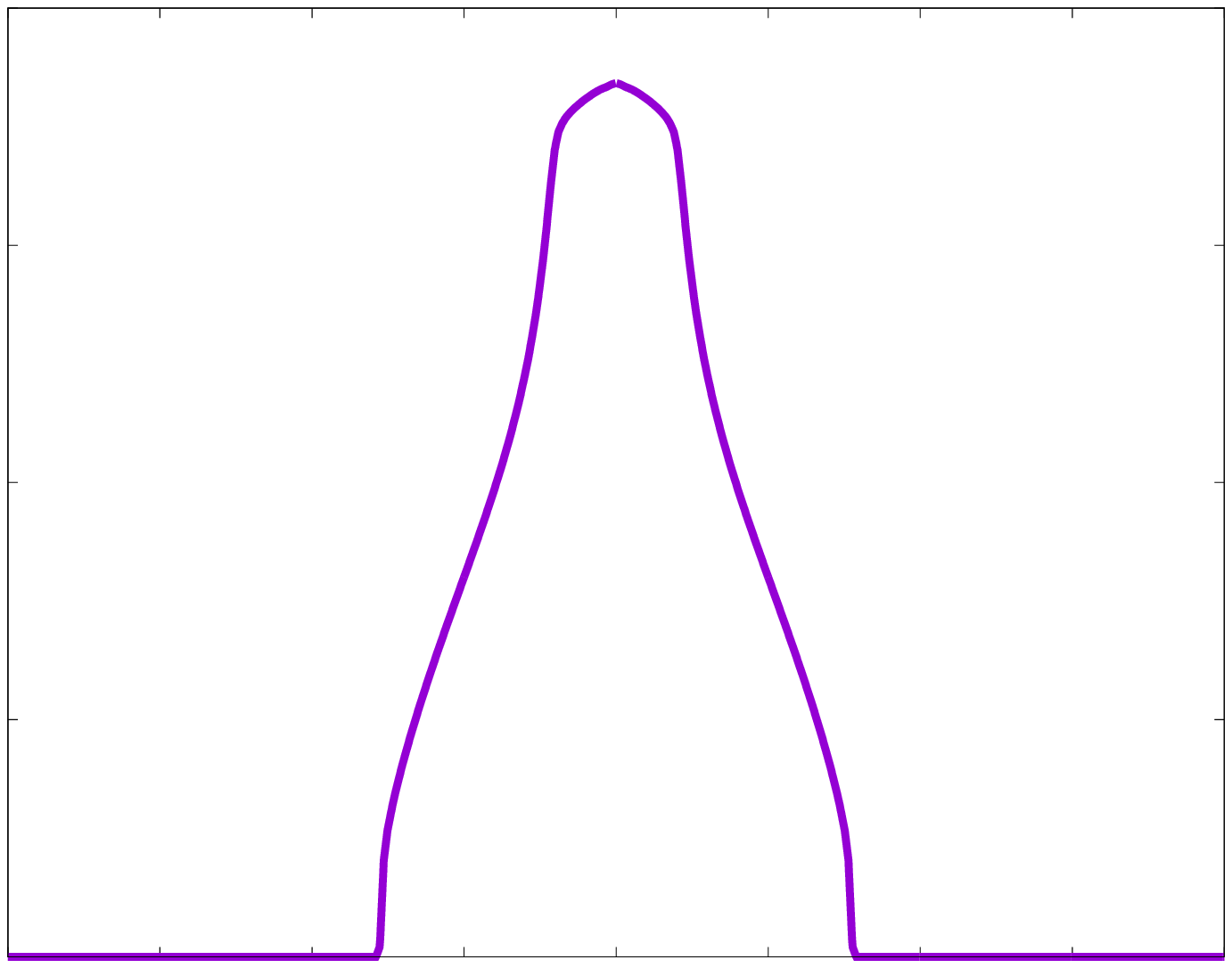}}%
    \gplfronttext
  \end{picture}%
  \endgroup
  }
  \end{minipage}
  \end{center}
  \vspace{-1mm}
  \hspace{9mm}
  \psfrag{a}    [c][c] {a)}
  \psfrag{b}    [c][c] {b)}
  \psfrag{c}    [c][c] {c)}
  \psfrag{pmax} [l][l] {$1$}
  \psfrag{p}    [l][l] {$d$}
  \psfrag{pmin} [l][l] {$0$}
  \includegraphics*[width=0.95\textwidth]{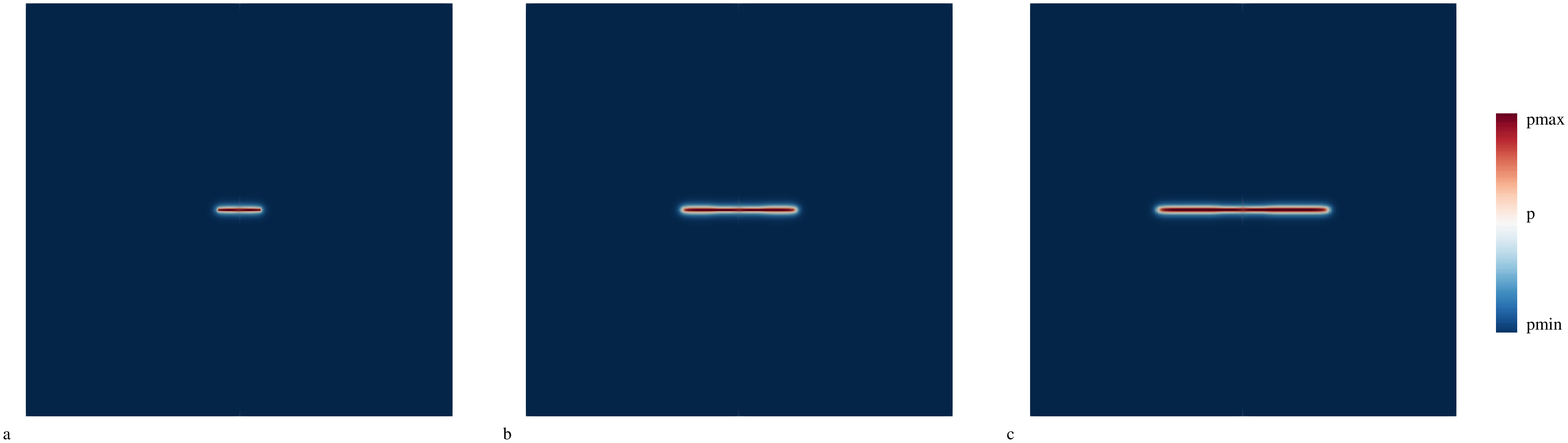}
  \captionsetup{width=0.9\textwidth}
  \caption{Hydraulically induced ductile fracture: Distribution of the fracture phase field and change of fluid content $m$ and $m^p$, fluid pressure $p$ and fracture opening $w$ over $x$ at $y=40$~m for three different time steps: a) $t=0$~s, b) $t=45$~s and c) $t=90$~s.}
 \label{fig:fracture_evolution}
\end{figure}

Finally, we depict a sequence of three snapshots in the course of the fracture evolution of the porous-elastic-plastic material in Figure~\ref{fig:fracture_evolution}. To be specific, we show the fracture phase field across the whole domain together with the change of fluid content, the pressure and the fracture-opening width at three different time steps along $x$ at $y=40$~m. The first time step is at $t=0$~s, the second time step is at $t=45$~s and the third time step is at the final state ($t=90$~s). As can be seen, all considered quantities are mainly concentrated in the center of the fracture. Such a concentration is less prominent in case of an elastic material, see Figure~\ref{fig:elast_plastic_x}.

\section{Conclusion}
\label{Section6}

A model for hydraulically induced fracturing in porous-elastic-plastic solids was developed in the present work. It incorporates a phase-field approach to fracture that is combined with a Drucker--Prager-type plasticity formulation and a Darcy--Biot-type fluid model. The model exploits a variational structure leading to a global minimization formulation. For this variational formulation it is crucial to introduce a plastic fluid content as an additional unknown yielding a constitutive fluid pressure in terms of only the elastic quantities. The global minimization structure demands the use of an $H(\div)$-conforming finite-element formulation, which has been implemented by means of Raviart--Thomas-type shape functions. The locking phenomenon of the plasticity formulation is overcome by using an enhanced-assumed-strain formulation for the deformation.

In the first numerical example a comparison of an undrained porous-elastic, an undrained porous-elastic-plastic and a drained elastic-plastic formulation was performed. Here, the different physical effects were investigated. It could be shown that the permeability of porous media can be considered as an additional hardening parameter. The second example shows the effect of the proposed modification of the fracture driving force. In the third example, a hydraulically induced crack in a porous-elastic and a porous-elastic-plastic medium was investigated. There it could be shown that neglecting the plastic effects underestimates the pressure inside the fracture and overestimates the fracture length.

{\bf Acknowledgments.}
This work was funded by the Deutsche Forschungsgemeinschaft (DFG, German Research Foundation) -- Project Number 327154368 -- SFB 1313. This funding is gratefully acknowledged.

\end{document}